\documentclass{article}

\usepackage[english]{babel}
\usepackage[utf8]{inputenc}
\usepackage[T1]{fontenc}

\usepackage[a4paper,top=3cm,bottom=2cm,left=3cm,right=3cm,marginparwidth=1.75cm]{geometry}

\usepackage{amssymb,amsmath,amsthm,enumitem}
\usepackage{graphicx}
\usepackage{authblk}
\usepackage{subcaption}
\usepackage{mathtools}
\usepackage{verbatim}
\usepackage{multirow}
\usepackage[justification=centering]{caption}
\usepackage{diagbox}
\usepackage[colorinlistoftodos]{todonotes}
\usepackage[colorlinks=true, allcolors=blue]{hyperref}
\usepackage{setspace}
\usepackage{pdfpages}
\usepackage{parskip}
\usepackage{csquotes}
\usepackage{apacite}
\usepackage{subcaption}
\usepackage[american]{circuitikz}
\usepackage{subfiles}
\usepackage{listings}
\usepackage{systeme}
\usepackage[sort&compress,square,comma,authoryear]{natbib}
\usepackage[bottom]{footmisc}
\usepackage{ marvosym }
\usepackage{algorithm} 
\usepackage{algpseudocode} 
\usepackage{tabularray-2021}
\usepackage{booktabs}
\hypersetup{
  colorlinks=false}
\usepackage{stackrel}
\usepackage{rsfso}
\usepackage{dsfont}
\usepackage{lipsum}
\usepackage{physics, color}
\usepackage[flushleft]{threeparttable}
\usepackage{soul} % for highlighting text
\usetikzlibrary{arrows.meta, fit,positioning}

\newcommand{\tikzscale}{0.75}
\newcommand{\tikznodesize}{0.5cm}

\newcommand{\tail}{\textnormal{t}}
\newcommand{\head}{\textnormal{h}}
\newcommand{\tin}{\textnormal{in}}
\newcommand{\tout}{\textnormal{out}}

\title{Uniformly sampling random directed hypergraphs with fixed degrees}
\author{Yanna J. Kraakman$^*$ and Clara Stegehuis}
\affil{Faculty of Electrical Engineering, Mathematics and Computer Science, University of Twente, Drienerlolaan 5, 7522 NB, Enschede, The Netherlands}
\affil{$^*$Corresponding author. Email: y.j.kraakman@utwente.nl}
\date{January 7, 2026}

\newtheorem{definition}{Definition}[section]
\newtheorem{theorem}{Theorem}
\newtheorem{lemma}{Lemma}[section]
\newtheorem{claim}{Claim}[section]
\newtheorem{property}{Property}[section]

\setcitestyle{authoryear,open={[},close={]}}

\begin{document}

\maketitle

\begin{abstract}
Many complex systems show non-pairwise interactions, which can be captured by hypergraphs. In this work, we propose an edge-swapping method to sample random directed hypergraphs with fixed vertex and hyperarc degrees, which can be applied to different classes of directed hypergraphs (containing self-loops, degenerate hyperarcs and/or multi-hyperarcs). We prove that this method indeed samples uniformly from the classes with self-loops and multi-hyperarcs, and that the method may not sample uniformly from classes without self-loops, or with self-loops and degenerate hyperarcs but without multi-hyperarcs. We present a partial result on the class with self-loops, but without degenerate hyperarcs or multi-hyperarcs.
\end{abstract}

\section{Introduction}
Many complex systems, from the spread of an epidemic to failures of the internet, can be modeled as mathematical graphs. Understanding these graphs is vital to predict the dynamics upon them, or to design better systems. A common method of analyzing a real-world graph is to compare the graph to a null model, a random graph model that matches the data in some features, while having a uniformly random structure with respect to other features. By comparing the graph to a random graph model, trivial features that can be expected based due to randomness can be distinguished from non-trivial features.

A popular random graph model is the random graph with a fixed degree sequence, sampled uniformly at random. A reason for fixing the degree sequence is that this strongly constrains important properties of the networks (\cite{newman2001}). Different methods exist to uniformly sample random graphs with fixed degrees. In the well-known \textit{stub-matching} sampling method, each edge in the original graph is cut in half, leaving two half-edges, or stubs. Then, these stubs are matched uniformly at random (\cite{bollobas1980}). By default, this method creates multigraphs, which may contain self-loops and multi-edges. However, often one is interested in sampling simple graphs only, in which case one has to reject the sample if it is non-simple and redo the procedure until a simple graph appears. If the variance of the vertex degrees is infinite, then the probability of obtaining a simple graph converges to zero for large graphs, so this rejection sampling is very time-consuming (\cite{hofstad2016}). Sometimes, stub-matching is referred to as the \textit{configuration model}, though other authors define the configuration model as the uniform distribution of graphs with a specified degree sequence. Stub-matching is then considered a method to sample from the configuration model (\cite{fosdick2018}, \cite{chodrow2019}). A different sampling method is \textit{edge-swapping} (\cite{petersen1891}), which can easily be tuned to always create simple graphs. This is a Markov chain Monte Carlo sampling method, during which edge pairs in the original graph swap heads sequentially. If a swap results in a self-loop or multi-edge, then this swap is undone. When the number of swaps tends to infinity, the obtained random simple graph converges in distribution to a uniform random simple graph, if the underlying Markov chain is regular, strongly connected and aperiodic. Depending on the graph class considered, it can be proved that the Markov chain has these properties for all initial graphs, or that there may be initial graphs for which it fails one of the properties. For undirected graphs, the edge-swapping method always samples uniformly from the classes of graphs that contain no self-loops, and the class that can contains both self-loops and multiple edges (proven by \cite{fosdick2018}, \cite{carstens2017}). For directed graphs, the edge-swapping method samples uniformly from the classes of graphs that can contain self-loops and have at least one vertex with in- or out-degree larger than 1. For the other classes, an edge-swap using three edges or a pre-sampling step may be necessary to assure uniformity of the method (proven by \cite{berger2010}, \cite{nishimura2018}). In such a pre-sampling step, for every triangle in the digraph an orientation is picked uniformly at random before applying the edge-swapping method.

Over the past decades, most research on graphs and null models has focused on graphs with pairwise connections. However, many complex systems show interactions that cannot be captured with only pairwise interactions (\cite{klamt2009}, \cite{battiston2021}, \cite{Bick2021}). This has caused a shift of focus from pairwise graphs to higher-order graphs, which can also describe interactions between groups of actors. One model for higher-order graphs is a simplicial complex, in which all group interactions are downwards closed. For this model, a null model has been proposed and studied (\cite{courtney2016}, \cite{young2017}). Here, we study a different model: a hypergraph, in which group structures need not be downwards closed. Hypergraphs are a generalization of pairwise graphs and model more general complex systems. As with graphs, hypergraphs can be characterized as undirected or directed. Most existing literature so far has focused on undirected hypergraphs, where groups interact together without different roles or directions. While this is realistic in examples such as social networks, there are many examples where directed hypergraphs are necessary instead. For instance, in graphs that model chemical reactions, the reaction $A+B\to C+D$ should be modeled by a directed hypergraph rather than a digraph or an undirected hypergraph, to capture that the reaction contains four reactants and that only $A$ and $B$ are necessary to start the reaction. Motivated by this example, we focus on directed hypergraphs.

Analyses of the edge-swapping sampling method on hypergraphs are scarce, in particular when the hypergraphs are directed. As far as we know, the only cases studied are undirected hypergraphs (\cite{chodrow2019}), and annotated hypergraphs (\cite{chodrow2020}), which are hypergraphs in which vertices can be of different types in every edge. In both analyses, the studied class of hypergraphs can contain multiple hyperedges, but no \textit{degenerate} hyperedges, which are hyperedges in which a vertex can appear more than once. The edge-swapping method samples uniformly from both these hypergraph classes. In some systems, however, degenerate hyperedges appear naturally and should thus appear in the null model. One example is the directed hypergraph that models chemical reactions, where the coefficients of molecules in a reaction can be interpreted as the multiplicity of a vertex in a hyperarc.

In this work, we define 16 classes of directed hypergraphs based on whether self-loops, degenerate hyperarcs and/or multiple hyperarcs are allowed and whether the directed hypergraphs are stub- or vertex-labeled. We then introduce an extension of the edge-swapping method to sample from these classes. This new method sequentially selects two random hyperarcs and mixes their heads and tails separately, while preserving the in- and out-degrees of the vertices and the tail- and head-degrees of the hyperarcs. We prove that for four classes the method always samples uniformly (Theorem \ref{thm:uniform_alldegrees}), and for ten classes it may sample with bias (Theorem \ref{thm:notuniform}). Analyzing the remaining two classes, where directed hypergraphs may contain self-loops, but no degenerate hyperarcs nor multiple hyperarcs, turned out to be the most difficult. Here, we prove that the introduced method samples uniformly from a special case of these classes, in which the directed hypergraphs have a degree sequence such that all hyperarcs have a tail of size 1 and a head of size 2 (corresponding to chemical reactions of the type $A\to B+C$), with two possible types of tails (Theorem \ref{thm:uniform_degrees(2,1)}). Lastly, we show that the results on vertex-labeled classes directly follow from the results on stub-labeled classes by introducing an acceptance probability for each edge-swap (Theorem \ref{thm:stub_to_vert}).

\paragraph{Overview of the paper.}
In Section \ref{section:model}, we introduce the directed hypergraph model and directed hypergraph spaces based on the classes. Section \ref{section:sampling} introduces our edge-swapping based sampling method. In Section \ref{section:results} we present our main results, specifying for which hypergraph spaces our swapping method samples uniformly, and we discuss these results in Section \ref{section:discussion}. Sections \ref{section:proofThm1+3}, \ref{section:proofThm2} and \ref{section:proofThm4} contain all proofs.

\section{Directed hypergraphs and their spaces}
\label{section:model}
We introduce our notation for directed hypergraphs and some of their characteristics in Section \ref{section:directed_hypergraphs}. We define all directed hypergraph spaces in Section \ref{section:configuration_models}.

\subsection{Directed hypergraphs}
\label{section:directed_hypergraphs}

We start by introducing notation for directed hypergraphs.
\\
\begin{definition}[Directed hypergraph]
    A directed hypergraph $H=(V,A)$ consists of a vertex set $V$ and a hyperarc multiset $A$. Each hyperarc $a \in A$ is defined by a pair $a^{\tail}$ and $a^{\head}$, which are both multisets of $V$. We write $a=(a^{\tail},a^{\head})$.
\end{definition}

We refer to $a^{\tail}$ and $a^{\head}$ as the tail resp. head of the hyperarc $a$. Now, we define self-loops, degenerate hyperarcs and multi-hyperarcs in the context of directed hypergraphs. These will be used in Section \ref{section:configuration_models} to define the hypergraph spaces.
\\
\begin{definition}[Self-loop, degenerate hyperarc, multi-hyperarcs]
    A hyperarc $a=(a^\tail,a^{\head})$ is called a \textit{self-loop} if $a^\tail=a^{\head}$ and a \textit{degenerate} hyperarc if $a^\tail$ or $a^\head$ contains at least one vertex with multiplicity $\geq 2$. Two hyperarcs $a_1=(a_1^\tail,a_1^{\head})$ and $a_2=(a_2^\tail,a_2^{\head})$ are called \textit{multi-hyperarcs} if $a_1^\tail=a_2^{\tail}$ and $a_1^\head=a_2^{\head}$.
\end{definition}

\begin{figure}[tbp]
\centering
\begin{subfigure}{0.55\textwidth}
    \centering
    \includegraphics[width=0.5\textwidth]{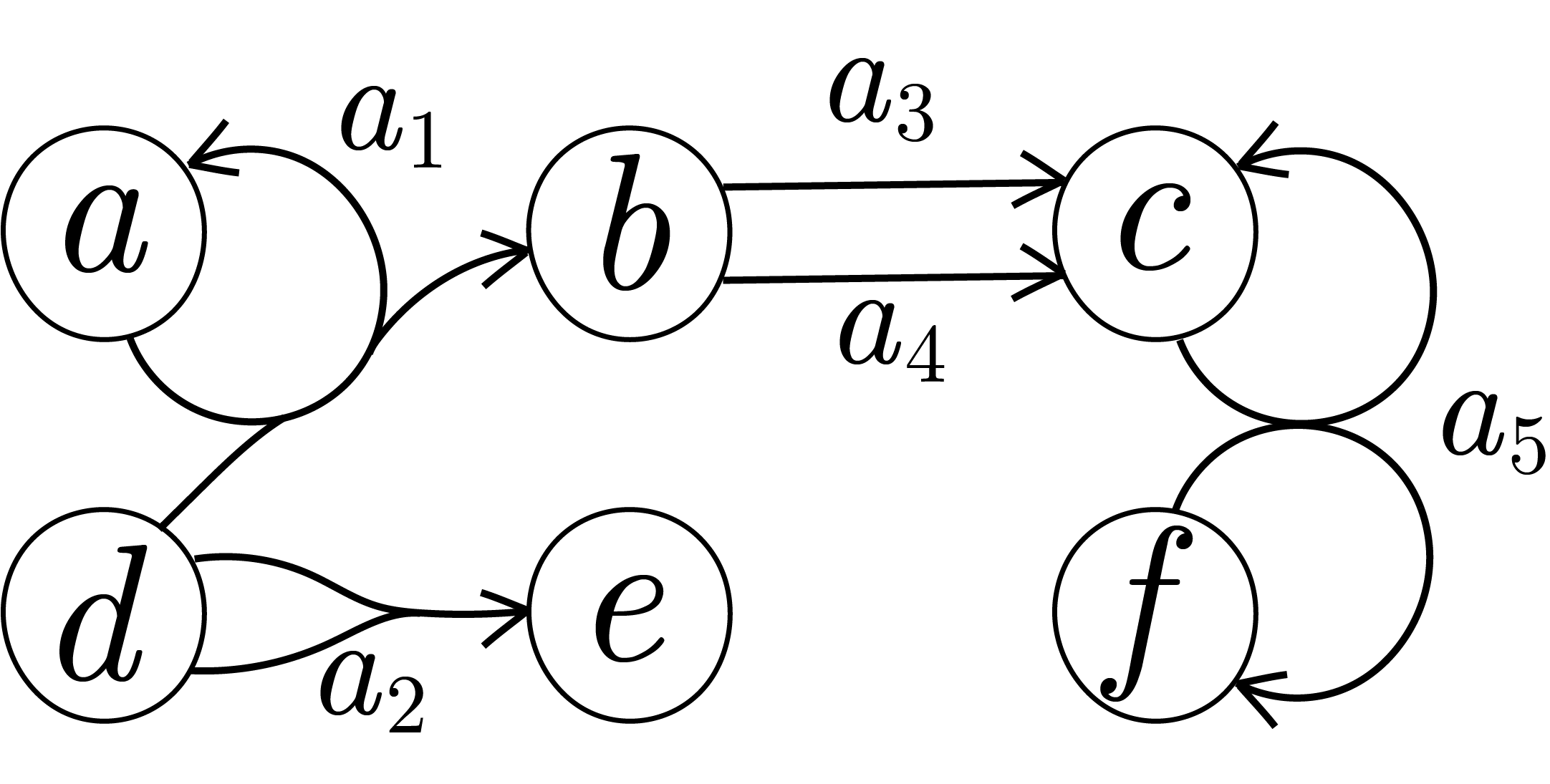}
    \caption{A representation of the directed hypergraph $H$.}
\end{subfigure}
\begin{subfigure}{0.35\textwidth}
    \centering
        \begin{tabular}{l}
            $\underline{A(H)}$ \\ 
            $a_1=(\{a,d\}, \{a,b\})$  \\
            $a_2=(\{d,d\}, \{e\})$ \\
            $a_3=(\{b\}, \{c\})$  \\
            $a_4=(\{b\}, \{c\})$ \\
            $a_5=(\{c,f\}, \{c,f\})$
        \end{tabular}
    \caption{The hyperarc set $A(H)$.}
\end{subfigure}
\caption{A directed hypergraph $H$ with vertex set $\{a,b,c,d,e,f\}$ and hyperarc set $\{a_1,a_2,a_3,a_4,a_5\}$, where $a_5$ is a self-loop, $a_2$ is a degenerate hyperarc and $a_3$ and $a_4$ are multi-hyperarcs. The vertex and hyperarc degree sequences of $H$ are $\vb*{d}_V=((1,1),(1,2),(3,1),(0,3),(1,0),(1,1))$ and $\vb*{\delta}_A=((2,2),(2,1),(1,1),(1,1),(2,2))$.}
\label{fig:example_directed_hypergraph}
\end{figure}

Figure \ref{fig:example_directed_hypergraph} shows an example of a directed hypergraph. 

Next, we define degrees of vertices and hyperarcs. These definitions are similar to other definitions of degrees in (directed) hypergraphs (\cite{ducournau2014}, \cite{chodrow2019}, \cite{ouvrard2020}), but unlike other definitions they take into account multiplicity of vertices in hyperedges.
\\
\begin{definition}[Vertex and hyperarc degrees]
    The in- and out-degree $d_v^{\tin}$ and $d_v^{\tout}$ of a vertex $v$ are the number of times the vertex appears in the heads resp. tails of arcs:
    \begin{align*}
        d_v^{\tin} &= \sum_{a \in A} m_{a^\head}(v)\\
        d_v^{\tout} &= \sum_{a \in A} m_{a^\tail}(v),
    \end{align*}
    where $m_s(v)$ denotes the multiplicity of element $v$ in multiset $s$. The tail- and head-degree $\delta_a^{\tail}$ and $\delta_a^{\head}$ of a hyperarc $a$ are the cardinalities of its tail resp. head:
    \begin{align*}
        \delta_a^{\tail} &= |a^{\tail}| = \sum_{v \in a^{\tail}} m_{a^\tail}(v) \\
        \delta_a^{\head} &= |a^{\head}| = \sum_{v \in a^{\head}} m_{a^\tail}(v).
    \end{align*}
\end{definition}

The degree sequence of a directed hypergraph consists of the vertex and hyperarc degree sequences.
\\
\begin{definition}[Degree sequence]
\label{def:degree_sequence}
    The vertex degree sequence $\vb*{d}_V$ of a directed hypergraph is a sequence of in- and out-degree pairs: $\vb*{d}_V= ((d_v^{\tin},d_v^{\tout}))_{v \in V}$. The hyperarc degree sequence $\vb*{\delta}_A$ of a directed hypergraph is a sequence of tail- and head-degree pairs: $\vb*{\delta}_A=((\delta_a^{\tail},\delta_a^{\head}))_{a \in A}$. The degree sequence $\vb*{d}$ of a directed hypergraph is the pair of its vertex and hyperarc degree sequences: $\vb*{d}=(\vb*{d}_V,\vb*{\delta}_A)$.
\end{definition}

\begin{figure}[tbp]
    \centering
    \includegraphics[width=0.4\textwidth]{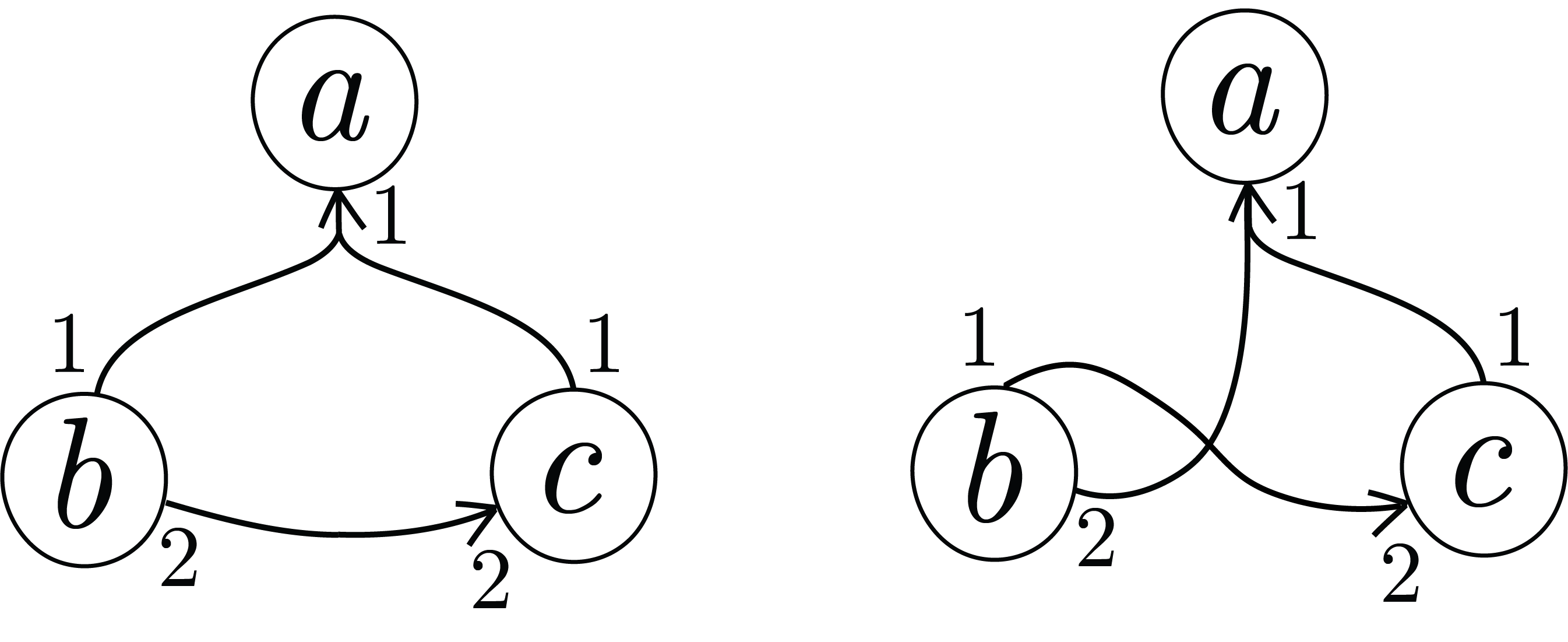}
    \caption{Two directed hypergraphs with their stubs labeled 1 or 2. The directed hypergraphs could be considered distinct (stub-labeled) or identical (vertex-labeled).}
    \label{fig:stub-labeled}
\end{figure}

Similarly to graphs, directed hypergraphs can be either vertex-labeled or stub-labeled. In stub-labeled directed hypergraphs, each connection between a vertex and a hyperarc, called a \textit{stub}, has a distinct label and meaning. In this setting, the two directed hypergraphs in Figure \ref{fig:stub-labeled}, with their stubs indicated by numbers, are considered different, even though they have the same vertex set $V=\{a,b,c\}$ and hyperarc set $A=\{(\{b,c\},\{a\}), (\{b\},\{c\})\}$. In vertex-labeled directed hypergraphs, the vertices have distinct labels, but the stubs are interchangeable. In that setting, the two directed hypergraphs in Figure \ref{fig:stub-labeled} are considered the same. When comparing given network data to a null model, these different notions of exchangeability may make a difference, and \cite{fosdick2018} notes that for most graphs that model real data, the vertex-labeled null model is more appropriate, with the exception of graphs modeling temporal data, or graphs in which the stubs are ordered. For a more detailed explanation of stub-labeled and vertex-labeled graphs, we refer to \cite{fosdick2018}.
\\
\begin{definition}[Stub-labeled and vertex-labeled directed hypergraphs]
    Let $v_1,v_2,\hdots,v_{d_v^{\tout}}$ be the out-stubs and let $v_{d_v^{\tout}+1},\hdots,v_{d_v^{\tout}+d_v^{\tin}}$ be the in-stubs of vertex $v$, for any $v \in V$. In a stub-labeled directed hypergraph, two directed hypergraphs $H_1$ and $H_2$ are equivalent if the sets of linked stubs in $H_1$ are the same as in $H_2$. In a vertex-labeled directed hypergraph, two directed hypergraphs $H_1$ and $H_2$ are equivalent if the sets of linked vertices in $H_1$ are the same as in $H_2$.
\end{definition}

When removing the stub labels on a stub-labeled directed hypergraph $H$, it becomes a vertex-labeled directed hypergraph $H'$. In particular, the directed hypergraph $H$ is one of the \textit{stub-labeled realizations} of $H'$.

While every stub in a stub-labeled directed hypergraph has a distinct label, the notions of self-loops, degenerate hyperarcs and multi-hyperarcs still apply; in their definitions, only the vertex labels are used. 

\subsection{Directed hypergraph spaces}
\label{section:configuration_models}

Directed hypergraphs may contain self-loops, degenerate hyperarcs and/or multi-hyperarcs. In the context of certain applications, however, some of these features may be nonsensical, such as a self-loop in a social network. We therefore identify different directed hypergraph spaces, that vary in whether or not the directed hypergraphs that they contain can have self-loops, degenerate hyperarcs and/or multi-hyperarcs and whether they are vertex- or stub-labeled.
\\
\begin{definition}[Directed hypergraph space]
    The directed hypergraph space $\mathcal{H}^{y}_{x} (\vb*{d})$, with $x \subseteq \{s,d,m\}$ and $y \in \{\textnormal{vert},\textnormal{stub}\}$, is a set containing all directed hypergraphs with degree sequence $\vb*{d}$ that are labeled as denoted by $y$ (vertex- or stub-labeled) and do not contain features $s$ (self-loops), $d$ (degenerate hyperarcs) or $m$ (multi-hyperarcs) that are not in $x$. 
\end{definition}

\begin{figure}[tbp]
    \centering
    \includegraphics[width=0.8\linewidth]{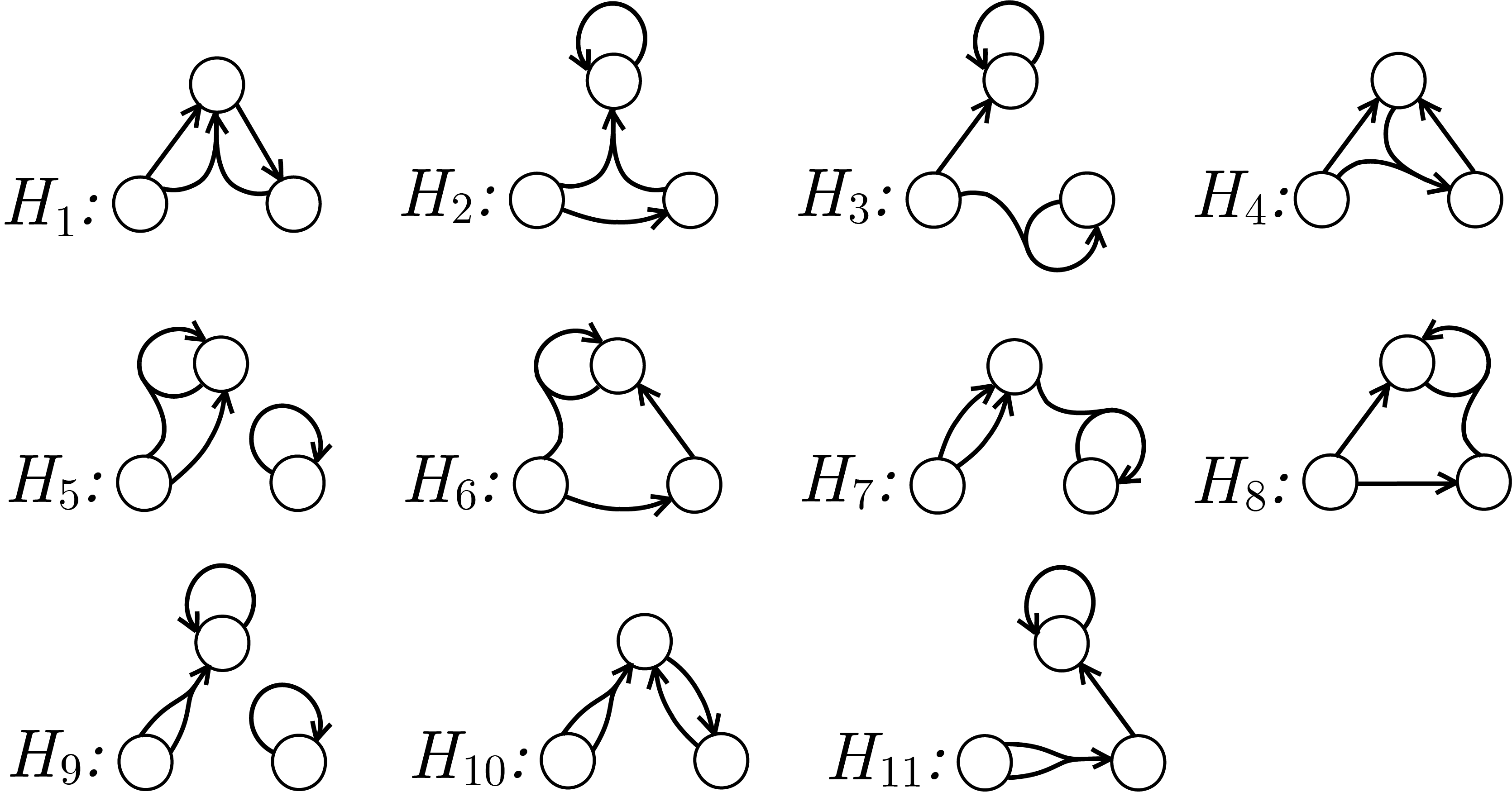}
    \caption{All vertex-labeled directed hypergraphs with $\vb*{d}_V=((2,1),(0,2),(1,1))$ and $\vb*{\delta}_A=((2,1),(1,1),(1,1))$. Some directed hypergraph spaces are $\mathcal{H}^{\textnormal{vert}}_{}(\vb*{d})=\{H_1,H_4,H_6,H_8\}$, $\mathcal{H}^{\textnormal{vert}}_{d}(\vb*{d})=\{H_1,H_4,H_6,H_8,H_{10}\}$, $\mathcal{H}^{\textnormal{vert}}_{s,m}(\vb*{d})=\{H_1,\hdots,H_8\}$ and $\mathcal{H}^{\textnormal{vert}}_{s,d,m}(\vb*{d})=\{H_1,\hdots,H_{11}\}$.}
    \label{fig:fixed_degrees}
\end{figure}

For any $\vb*{d}$, the least restrictive spaces are $\mathcal{H}^{\textnormal{vert}}_{s,d,m} (\vb*{d})$ and $\mathcal{H}^{\textnormal{stub}}_{s,d,m} (\vb*{d})$, and the most restrictive spaces are $\mathcal{H}^{\textnormal{vert}}_{} (\vb*{d})$ and $\mathcal{H}^{\textnormal{stub}}_{} (\vb*{d})$. Figure \ref{fig:fixed_degrees} shows all vertex-labeled directed hypergraphs with a specified degree sequence and illustrates certain spaces. These were found by enumerating all possible hyperarc multisets.

\section{Sampling from the hypergraph spaces} 
\label{section:sampling}
\begin{figure}[tbp]
    \centering
    \includegraphics[width=\linewidth]{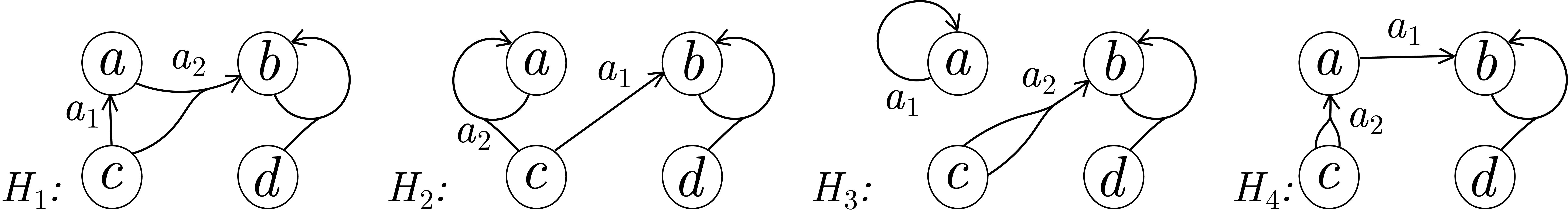}
    \caption{The double hyperarc shuffle $s(a_1,a_2|H_1)$ has state space $S(a_1,a_2|H_1)=\{H_1,H_2,H_3,H_4\}$.}
    \label{fig:double_hyperarc_shuffle}
\end{figure}

We introduce the \textit{double hyperarc shuffle} (Figure \ref{fig:double_hyperarc_shuffle}), which generalizes the edge-swap, and is an extension of the pairwise reshuffle defined by \cite{chodrow2019} for undirected hypergraphs. The double hyperarc shuffle takes two hyperarcs $a,b$ and `shuffles' the tails of $a$ and $b$ and the heads of $a$ and $b$ separately, while keeping the tail- and head-degrees of both hyperarcs fixed. By construction, the double hyperarc shuffle maintains both the vertex and the hyperarc degree sequence of a directed hypergraph. Let $A \uplus B$ denote the multiset union with multiplicity $m_{A \uplus B}(x)=m_A(x)+m_B(x)$. For example, if $A=\{a,b,b\}$ and $B=\{a,b,c\}$ then $A \uplus B=\{a,a,b,b,b,c\}$.
\\
\begin{definition}[Double hyperarc shuffle]
    Let $a,b$ be two distinct hyperarcs. The double hyperarc shuffle $s(a,b)$ is a random trial with the sample space
    \begin{align*}
        S(a,b) = \{ (\hat{a},\hat{b}): a^{\tail} \uplus b^{\tail} &= \hat{a}^{\tail} \uplus \hat{b}^{\tail}\\
        a^{\head} \uplus b^{\head} &= \hat{a}^{\head} \uplus \hat{b}^{\head}\\
        (\delta_a^{\tail},\delta_a^{\head}) &= (\delta_{\hat{a}}^{\tail},\delta_{\hat{a}}^{\head})\\
        (\delta_b^{\tail},\delta_b^{\head}) &= (\delta_{\hat{b}}^{\tail},\delta_{\hat{b}}^{\head}) \},
    \end{align*}
    where each outcome in the sample space is equally likely. Let $H$ be a directed hypergraph with $a,b \in A(H)$. Performing a double hyperarc shuffle on $H$ is denoted by $s(a,b|H)$ and is a random trial with the sample space 
    \begin{align*}
        S(a,b|H) = \{ H': A(H')=(A(H) \backslash \{a,b\}) \uplus \{\hat{a},\hat{b}\} \textnormal{ s.t. } (\hat{a},\hat{b}) \in S(a,b)\},
    \end{align*}
    where each outcome is equally likely. If the hyperarc $a$ appears twice in the hyperarc multiset, then $a=b$ is possible.
\end{definition}

Note that $S(a,b) = S(b,a)$. A double hyperarc shuffle $s(a,b|H)$ may result in a directed hypergraph $H' \in S(a,b|H)$ that is not in the considered directed hypergraph space $\mathcal{H}^{y}_{x}(\vb*{d})$, e.g., when $\mathcal{H}^{\textnormal{stub}}_{d,m}(\vb*{d})$ is considered and $H'$ contains a self-loop. Then, the double hyperarc shuffle is undone and the result is the directed hypergraph $H$. This is also known as a `switching and holding' mechanism (\cite{artzy-randrup2005}).

Performing the switching and holding mechanism can be considered a Markov chain, described by a random walk on a \textit{graph of graphs}, which consists of vertices that represent the directed hypergraphs in the considered space and arcs that represent the transitions possible by applying double hyperarc shuffles. In particular, let $\mathcal{G}(\mathcal{H}^y_x(\vb*{d})) = (\mathcal{V}, \mathcal{A})$ be the graph of graphs on the hypergraph space $\mathcal{H}^y_x(\vb*{d})$. Then, $\mathcal{V} = \mathcal{H}^y_x(\vb*{d})$. 
Furthermore, $(H,H')\in \mathcal{A}$ if and only if there exists a double hyperarc shuffle which transforms $H$ into $H'$. Every arc $(H,H') \in \mathcal{A}$ has an associated probability $p(H'|H)$ that equals the probability that the double hyperarc shuffle $s(a,b|H)$ with two uniformly at random picked hyperarcs $a,b \in A(H)$ results in $H'$.  

The outcome of performing $k$ double hyperarc shuffles sequentially, and undoing double hyperarc shuffles that create directed hypergraphs not in the considered directed hypergraph space, is called the \textit{double hyperarc shuffling method with $k$ steps}, and can be considered a random walk of length $k$ on $\mathcal{G}(\mathcal{H}^y_x(\vb*{d}))$.
\\
\begin{definition}[Double hyperarc shuffling method with $k$ steps]
    Given some directed hypergraph space $\mathcal{H}^y_x(\vb*{d})$ and an initial directed hypergraph $H_0 \in \mathcal{H}^y_x(\vb*{d})$, $\forall i \in \mathds{N}$ let $H_{i+1}$ be the outcome of $s(a_i,b_i|H_i)$, with $a_i,b_i \in A(H_i)$ picked uniformly at random. If the outcome of $s(a_i,b_i|H_i)$ is not an element of $\mathcal{H}^y_x(\vb*{d})$ then $H_{i+1}=H_i$. The double hyperarc shuffling method on the directed hypergraph space $\mathcal{H}^y_x(\vb*{d})$ is a random trial, of which the outcome is $H_k$.
\end{definition}

\section{Results}
\label{section:results}
In this section, we present our main results. As explained in Section \ref{section:configuration_models}, we distinguish 16 directed hypergraph spaces that one could sample from depending on whether self-loops, degenerate hyperarcs or multi-hyperarcs are allowed. From those, 8 are stub-labeled and 8 are their vertex-labeled equivalents. For each space, we show whether the double hyperarc shuffling method samples uniformly from it or not. The convergence that is used in Theorems \ref{thm:uniform_alldegrees}, \ref{thm:uniform_degrees(2,1)} and \ref{thm:stub_to_vert} is convergence in distribution. 

The following theorem establishes the spaces from which the double hyperarc shuffling method samples uniformly:
\\
\begin{theorem}
\label{thm:uniform_alldegrees}
    For $k \rightarrow \infty$, the double hyperarc shuffling method with $k$ steps on the space $\mathcal{H}^{\textnormal{stub}}_{s,d,m}(\vb*{d})$ samples uniformly from that space, and the double hyperarc shuffling method with $k$ steps on the space $\mathcal{H}^{\textnormal{stub}}_{s,m}(\vb*{d})$ samples uniformly from that space, for any degree sequence $\vb*{d}$.
\end{theorem}

Thus, according to Theorem \ref{thm:uniform_alldegrees}, the double hyperarc shuffling method on $\mathcal{H}^{\textnormal{stub}}_{s,d,m}(\vb*{d})$ samples uniformly this space that has no restrictions, other than the degree sequence of the directed hypergraphs: it allows self-loops, degenerate hyperarcs and multi-hyperarcs. This is the most general result, showing that the double hyperarc shuffling method can reshuffle any initial directed hypergraph into any other directed hypergraph with the same degree sequence. Secondly, the method on $\mathcal{H}^{\textnormal{stub}}_{s,m}(\vb*{d})$ samples uniformly from the space that has the restriction of no degenerate hyperarcs. These results generalize results on sampling directed multigraphs from a space with multigraphs that can contain self-loops and multi-arcs, which can be done using a similar edge-swapping method (\cite{carstens2017}).

We continue with a result on sampling from $\mathcal{H}^{\textnormal{stub}}_{s}(\vb*{d})$, which is the space with directed hypergraph that possibly contain self-loops, but no degenerate or multi-hyperarcs. We prove uniform sampling for a very restricted type of degree sequence.
\\
\begin{theorem}
\label{thm:uniform_degrees(2,1)}
    For $k \rightarrow \infty$, the double hyperarc shuffling method with $k$ steps on the space $\mathcal{H}^{\textnormal{stub}}_{s}(\vb*{d})$ samples uniformly from that space, for any degree sequence $\vb*{d}$ with $\forall a \in A: \delta_a^{\tail}=1, \delta_a^{\head}=2$ and $|\{v: d_v^{\textnormal{out}}>0\}|=2$.
\end{theorem}

Proving that the double hyperarc shuffling method on the space $\mathcal{H}^{\textnormal{stub}}_{s}(\vb*{d})$ samples uniformly from that space turns out to be a significantly greater challenge than proving the results of Theorem \ref{thm:uniform_alldegrees}. Theorem \ref{thm:uniform_degrees(2,1)} presents a result for the space with restrictions on $\vb*{d}$, in which the directed hypergraphs can be represented as colored graphs, allowing for graph theoretical techniques to be used in the proof. This result may be a setup for a more general theorem, as we explain in the discussion. The result generalizes the result on uniformly sampling directed graphs from a space with graphs that may contain self-loops and no multi-arcs, which can be done using a similar edge swapping method (\cite{rao1996}, \cite{berger2010}). By symmetry, the result of Theorem \ref{thm:uniform_degrees(2,1)} also holds for any degree sequence $\vb*{d}$ with $\delta_a^{\tail}=2, \delta_a^{\head}=1$ where $|\{v: d_v^{\textnormal{in}}>0\}|=2$. It remains an open question whether the theorem holds for any degree sequence $\vb*{d}$.

The next result shows that there are spaces from which the double hyperarc shuffling method does not sample uniformly for some degree sequences.
\\
\begin{theorem}
\label{thm:notuniform}
    There exist degree sequences $\vb*{d}_1, \vb*{d}_2$ such that the double hyperarc shuffling method with $k$ steps does not sample uniformly from the spaces $\mathcal{H}^{\textnormal{stub}}_{s,d}(\vb*{d}_1)$, $\mathcal{H}^{\textnormal{stub}}_{d,m}(\vb*{d}_2)$, $\mathcal{H}^{\textnormal{stub}}_{d}(\vb*{d}_2)$, $\mathcal{H}^{\textnormal{stub}}_{m}(\vb*{d}_2)$, and $\mathcal{H}^{\textnormal{stub}}_{}(\vb*{d}_2)$, for any $k$.
\end{theorem}

The proof of Theorem \ref{thm:notuniform} is by direct example of such degree sequences, and will follow in Section~\ref{section:proofThm1+3}. The results for $\mathcal{H}^{\textnormal{stub}}_{d,m}(\vb*{d}_2)$, $\mathcal{H}^{\textnormal{stub}}_{d}(\vb*{d}_2)$, $\mathcal{H}^{\textnormal{stub}}_{m}(\vb*{d}_2)$, and $\mathcal{H}^{\textnormal{stub}}_{}(\vb*{d}_2)$ follow directly from similar results for directed graphs: for hyperarc degree sequences $\vb*{d}_2$ with $\delta_a^{\tail}=\delta_a^{\head}=1$, $\mathcal{H}^{\textnormal{stub}}_{d,m}(\vb*{d}_2)$ and $\mathcal{H}^{\textnormal{stub}}_{m}(\vb*{d}_2)$ reduce to spaces of multigraphs without self-loops, and there exist vertex degree sequences for which our edge-swapping method cannot sample from those spaces (\cite{carstens2017}). Similarly, $\mathcal{H}^{\textnormal{stub}}_{d}(\vb*{d}_2)$ and $\mathcal{H}^{\textnormal{stub}}_{}(\vb*{d}_2)$ reduce to spaces of simple digraphs without self-loops, and our edge-swapping method cannot sample from those spaces, for some vertex degree sequences (\cite{rao1996}, \cite{berger2010}). The result for $\mathcal{H}^{\textnormal{stub}}_{s,d}(\vb*{d}_1)$ is remarkable: our edge-swapping method can sample from the space of stub-labeled simple graphs with possible self-loops, while this does not hold for the space of directed hypergraphs with self-loops and degenerate hyperarcs. Thus, for certain degree sequences, degenerate hyperarcs prevent uniform sampling from the directed hypergraph space using the double hyperarc shuffling method.

Lastly, we show that the results on stub-labeled spaces carry over to vertex-labeled spaces, when acceptance probabilities are used. We first introduce these acceptance probabilities. The acceptance probability of a double hyperarc shuffle $s(a,b|H)$ resulting in $H' \in S(a,b|H)$ in a vertex-labeled space is
\begin{align}
\label{eq:acceptance_prob}
    \alpha(H'|H) &= 
    (m_a m_b)^{-1} \prod_{v \in V} \Big( \binom{m_{\hat{a}^{\head}}(v) + m_{\hat{b}^{\head}}(v)}{m_{\hat{a}^{\head}}(v)} \binom{m_{\hat{a}^{\tail}}(v) + m_{\hat{b}^{\tail}}(v)}{m_{\hat{a}^{\tail}}(v)} \Big)^{-1},
\end{align}
where $\hat{a}$ and $\hat{b}$ are the hyperarcs in $H'$ that resulted from the shuffle of $a$ and $b$.
\\
\begin{theorem}
\label{thm:stub_to_vert}
    For $k \rightarrow \infty$, the double hyperarc shuffling method with $k$ steps on $\mathcal{H}^{\textnormal{stub}}_{x}(\vb*{d})$ with the acceptance probabilities as in \eqref{eq:acceptance_prob} samples uniformly from $\mathcal{H}^{\textnormal{vert}}_{x}(\vb*{d})$ if and only if the double hyperarc shuffling method with $k$ steps on $\mathcal{H}^{\textnormal{stub}}_{x}(\vb*{d})$ samples uniformly from $\mathcal{H}^{\textnormal{stub}}_{x}(\vb*{d})$.
\end{theorem}

This result generalizes similar results for undirected graphs (\cite{fosdick2018}) and undirected hypergraphs (\cite{chodrow2019}).

Theorems \ref{thm:uniform_alldegrees}, \ref{thm:uniform_degrees(2,1)} and \ref{thm:notuniform} are summarized in Table \ref{tab:results}. 

\begin{table}[h!]
\centering
\begin{threeparttable}
    \centering
    \caption{Sampling from the 8 stub-labeled directed hypergraph classes.}
    \begin{tabular}{c|c|c|c||c}
        Directed & Self-loops & Degenerate & multi- & Uniform sampling using double\\ 
        hypergraph space & & hyperarcs & hyperarcs & hyperarc shuffling method?\\ \hline
        $\mathcal{H}^{\textnormal{stub}}_{}(\vb*{d})$ & No & No & No & No\tnote{1}\\
        $\mathcal{H}^{\textnormal{stub}}_{d}(\vb*{d})$ & No & Yes & No & No\tnote{1}\\
        $\mathcal{H}^{\textnormal{stub}}_{m}(\vb*{d})$ & No & No & Yes & No\tnote{1}\\
        $\mathcal{H}^{\textnormal{stub}}_{d,m}(\vb*{d})$ & No & Yes & Yes & No\tnote{1}\\
        $\mathcal{H}^{\textnormal{stub}}_{s}(\vb*{d})$ & Yes & No & No & Sometimes\tnote{2}\\
        $\mathcal{H}^{\textnormal{stub}}_{s,d}(\vb*{d})$ & Yes & Yes & No & No\tnote{1}\\
        $\mathcal{H}^{\textnormal{stub}}_{s,m}(\vb*{d})$ & Yes & No & Yes & Yes\\
        $\mathcal{H}^{\textnormal{stub}}_{s,d,m}(\vb*{d})$ & Yes & Yes & Yes & Yes       
    \end{tabular}
     \begin{tablenotes}
  \item[1] The double hyperarc shuffling method does not sample uniformly from this space for at least one degree sequence.
  \item[2] Yes for degree sequences $\vb*{d}$ with $\delta_a^{\tail}=1, \delta_a^{\head}=2$, $|\{v: d_v^{\textnormal{out}}>0\}|=2$.
  \end{tablenotes}
    \label{tab:results}
\end{threeparttable}
\end{table}

\section{Discussion}
\label{section:discussion}

We now provide some discussion of our results.

\paragraph{Defining the degree sequence.} In this paper, we define the vertex degree sequence as pairs of in- and out-degrees of vertices (Definition \ref{def:degree_sequence}). The double hyperarc shuffle preserves these in- and out-degree pairs, similarly to the edge-swapping method for directed graphs that preserves the in- and out-degree pairs. Alternatively, one could define the vertex degree sequence as a multiset of in-degrees and a multiset of out-degrees. Sampling from such a space would require a different sampling method, that allows for vertices to swap their in-degree while fixing their out-degree. 

\paragraph{Sampling from the spaces of Theorem~\ref{thm:notuniform}} In Theorem \ref{thm:notuniform} we show that the double hyperarc shuffling method does not sample uniformly from several hypergraph spaces. However, the counterexamples describe rather small directed hypergraphs. It would be interesting to investigate whether adding extra conditions on the degree sequences may enable us to uniformly sample from these spaces, for example by requiring a minimum hyperarc degree. Other options to still be able to sample from these spaces could be to use a more complex version of the double hyperarc shuffle, for example by using three hyperarcs. This is also common in edge-swapping sampling methods on graphs (\cite{erdos2009}, \cite{nishimura2018}). Alternatively, one could try sampling from a more general space, and rejecting the sample if it does not meet the requirements.

\paragraph{Mixing times.} Theorems \ref{thm:uniform_alldegrees}, \ref{thm:uniform_degrees(2,1)} and \ref{thm:stub_to_vert} show that for $k \rightarrow \infty$, the double hyperarc shuffling method with $k$ steps samples uniformly from certain spaces. In practice, a relevant question is how large $k$ must be to obtain a sample that is approximately uniform. This is called the \textit{mixing time} (\cite{levin2008}). For several types of graphs, the mixing time is proven to be polynomial in the number of vertices (\cite{cooper2007}, \cite{Greenhill2011}, \cite{Erdos2019}). An open question is whether analogous bounds can be established for directed hypergraphs. There are also experimental approaches to analyzing mixing times for graphs \cite{dutta2024}. In the hypergraph setting, \cite{kraakman2025}, building on the preprint of this work, experimentally investigate the mixing time of the double hyperarc shuffling method.

\paragraph{Annotated hypergraphs.} As briefly mentioned in the introduction, \cite{chodrow2020} analyzed annotated hypergraphs, which are hypergraphs in which the vertices can have different roles in each hyperedge. Directed hypergraphs can be considered as annotated hypergraphs with two possible roles for vertices in a hyperedge: tail or head. \cite{chodrow2020} analyze the edge-swapping method for sampling from one class of annotated hypergraphs: one without degenerate hyperedges. It would be insightful to extend our analysis of sampling from the 16 spaces of directed hypergraphs to annotated hypergraphs. 

\paragraph{Generalizing Theorem \ref{thm:uniform_degrees(2,1)}.} For all theorems, the most challenging part in proving that the double hyperarc shuffling method samples uniformly from a space is proving that the underlying Markov chain is strongly connected. Generalizing the result of Theorem \ref{thm:uniform_degrees(2,1)} to a similar result without a restriction on the number of different tail vertices, $|\{v: d_v^{\textnormal{out}}>0\}|$, may be possible using a similar proof technique as we have used for Theorem \ref{thm:uniform_degrees(2,1)}. The difficulty lies in finding two sets of hyperarcs that allow us to construct a path in the Markov chain between states $H_1$ and $H_2$. We believe that generalizing to a similar result with a different restriction on the tail-degrees of the hyperarcs, $\delta_a^{\tail}=1$, may be possible using a multi-layer network approach of the proof of our result. One could, for example, fix $\forall a \in A: \delta_a^{\tail}=\delta_a^{\head}=2$, map the tails of the hyperarcs to edges in the first layer of a graph and the heads of the hyperarcs to edges in the second layer of the graph, and link both edges in the multi-layer graph. Generalizing to a similar result without any restriction on the degree sequence probably needs a different proof mechanism, since our method of finding the sets described above heavily relies on the head-degrees of the hyperarcs being 2.

\paragraph{Self-loops.} The used definition of a self-loop in a directed hypergraph requires the tail and head of a hyperarc to be equal. A less restrictive definition of a self-loop could be a hyperarc whose head and tail have a nonempty intersection. For this alternative definition of a self-loop, the results are the same, as Lemma \ref{lemma:stub_many_notconnected}, which results in Theorem \ref{thm:notuniform}, continues to hold. 

\section{Proofs of Theorems \ref{thm:uniform_alldegrees} and \ref{thm:notuniform}}
\label{section:proofThm1+3}
In this section, we present the proofs of Theorems \ref{thm:uniform_alldegrees}, and \ref{thm:notuniform}.

The double hyperarc shuffling method on a space $\mathcal{H}^{\textnormal{stub}}_x(\vb*{d})$ samples uniformly from that space if the underlying Markov chain has a uniform stationary distribution and this stationary distribution is unique. This is the case if for the graph of graphs $\mathcal{G}(\mathcal{H}^{\textnormal{stub}}_x(\vb*{d}))$ holds (\cite{levin2008})
\begin{enumerate}
    \item $\mathcal{G}(\mathcal{H}^{\textnormal{stub}}_x(\vb*{d}))$ is regular:
    \\ $\sum_{H' \in \mathcal{H}^{\textnormal{stub}}_x(\vb*{d})} p(H'|H)$ is equal for all $H \in \mathcal{H}^{\textnormal{stub}}_x(\vb*{d})$ and $\sum_{H' \in \mathcal{H}^{\textnormal{stub}}_x(\vb*{d})} p(H|H')$ is equal for all $H \in \mathcal{H}^{\textnormal{stub}}_x(\vb*{d})$.
    \item $\mathcal{G}(\mathcal{H}^{\textnormal{stub}}_x(\vb*{d}))$ is strongly connected:\\
    For every two vertices $v_1,v_2 \in \mathcal{G}(\mathcal{H}^{\textnormal{stub}}_x(\vb*{d}))$ there is a path $v_1 \rightarrow \cdots \rightarrow v_2$ in $\mathcal{G}(\mathcal{H}^{\textnormal{stub}}_x(\vb*{d}))$.
    \item $\mathcal{G}(\mathcal{H}^{\textnormal{stub}}_x(\vb*{d}))$ is aperiodic:\\
    The greatest common divisor of the length of the cycles in $\mathcal{G}(\mathcal{H}^{\textnormal{stub}}_x(\vb*{d}))$ is 1.
\end{enumerate}

For the directed hypergraph spaces in Theorem \ref{thm:uniform_alldegrees}, we prove the first and third property in Section \ref{section:proof_regularity_aperiodicity}. We prove the second, most difficult to prove property in Sections \ref{section:strong_connectivity_H_sdm} and \ref{section:strong_connectivity_H_sm}. For the directed hypergraph spaces in Theorem \ref{thm:notuniform}, counterexamples of the second property are given in Section \ref{section:no_connectivity}.

\subsection{Regularity and aperiodicity}
\label{section:proof_regularity_aperiodicity}
The following lemma shows that the graph of graphs of any directed hypergraph space that we consider is regular.
\\
\begin{lemma}
\label{lemma:regular}
    $\mathcal{G}(\mathcal{H}^{\textnormal{stub}}_x(\vb*{d}))$ is a regular graph, for any degree sequence $\vb*{d}$ and any $x \subseteq \{s,d,m\}$.
\end{lemma}
\begin{proof}
    By construction of the double hyperarc shuffle, $\sum_{H' \in \mathcal{H}} p(H'|H)$ is equal for all $H \in \mathcal{H}^{\textnormal{stub}}_x(\vb*{d})$:
    \begin{align*}
        \sum_{H' \in \mathcal{H}^{\textnormal{stub}}_x(\vb*{d})} p(H'|H) = 1.
    \end{align*}
    Next, we prove that $\sum_{H' \in \mathcal{H}^{\textnormal{stub}}_x(\vb*{d})} p(H|H')$ is equal for all $H \in \mathcal{H}^{\textnormal{stub}}_x(\vb*{d})$, by showing that each transition in $\mathcal{G}(\mathcal{H}^{\textnormal{stub}}_x(\vb*{d}))$ is bi-directional, i.e., that $p(H'|H)=p(H|H')$. To that end, observe that when applying double hyperarc shuffles in a stub-labeled directed hypergraph space, the probability of $s(a,b|H)$ resulting in any $H' \in S(a,b|H)$ is
      \begin{align*}
        p(H'|H) = 
             \binom{|A(H)|}{2}^{-1} \binom{\delta_a^{\tail}+\delta_b^{\tail}}{\delta_a^{\tail}}^{-1} \binom{\delta_a^{\head}+\delta_b^{\head}}{\delta_a^{\head}}^{-1},
    \end{align*}
    where $\binom{|A(H)|}{2}$ is the number of hyperarc pairs available for shuffling, and $\binom{\delta_a^{\tail}+\delta_b^{\tail}}{\delta_a^{\tail}}^{-1} \binom{\delta_a^{\head}+\delta_b^{\head}}{\delta_a^{\head}}^{-1}$ is the probability of a specific stub-labeled realization, given some hyperarc pair $a$ and $b$. Now, consider two directed hypergraphs $H,H'\in \mathcal{H}^{\textnormal{stub}}_x(\vb*{d})$ with $H' \in S(a,b|H)$ and $A(H')=(A(H) \backslash \{a,b\})\uplus\{\hat{a},\hat{b}\}$, where $a,b \in A(H)$ are two hyperarcs and $\hat{a},\hat{b}$ are the shuffled hyperarcs. 
    Then,
    \begin{align*}
        p(H'|H) &= \binom{|A(H)|}{2}^{-1} \binom{\delta_a^{\tail}+\delta_b^{\tail}}{\delta_a^{\tail}}^{-1} \binom{\delta_a^{\head}+\delta_b^{\head}}{\delta_a^{\head}}^{-1} \\
        &= \binom{|A(H')|}{2}^{-1} \binom{\delta_{\hat{a}}^{\tail}+\delta_{\hat{b}}^{\tail}}{\delta_{\hat{a}}^{\tail}}^{-1} \binom{\delta_{\hat{a}}^{\head}+\delta_{\hat{b}}^{\head}}{\delta_{\hat{a}}^{\head}}^{-1} \\
        &= p(H|H'),
    \end{align*}
        where the second equality follows from the degree preservation of the double hyperarc shuffle. This proves that each arc in $\mathcal{G}(\mathcal{H}^{\textnormal{stub}}_x(\vb*{d}))$ is bi-directional. In particular,
    \begin{align*}
        \sum_{H' \in \mathcal{H}^{\textnormal{stub}}_x(\vb*{d})} p(H|H') = \sum_{H' \in \mathcal{H}^{\textnormal{stub}}_x(\vb*{d})} p(H'|H) = 1,
    \end{align*}
    for every $H \in \mathcal{H}^{\textnormal{stub}}_x(\vb*{d})$.
\end{proof}

The following lemma shows that the graph of graphs of any directed hypergraph space that we consider is aperiodic.
\\
\begin{lemma}
\label{lemma:aperiodic}
    $\mathcal{G}(\mathcal{H}^{\textnormal{stub}}_x(\vb*{d}))$ is an aperiodic graph, for any degree sequence $\vb*{d}$ and any $x \subseteq \{s,d,m\}$.
\end{lemma}
\begin{proof}
    Since every double hyperarc shuffle can leave the directed hypergraph unchanged, i.e., $p(H|H)>0$ for all $H \in \mathcal{H}^{\textnormal{stub}}_x(\vb*{d})$, $\mathcal{G}(\mathcal{H}^{\textnormal{stub}}_x(\vb*{d}))$ contains self-loops. Therefore, $\mathcal{G}(\mathcal{H}^{\textnormal{stub}}_x(\vb*{d}))$ is aperiodic. 
\end{proof}

\subsection[Strong connectivity of G(H(sdm))]{Strong connectivity of $\mathcal{G}(\mathcal{H}^{\textnormal{stub}}_{s,d,m}(\vb*{d}))$}
\label{section:strong_connectivity_H_sdm}
Now, we show that the graph of graphs of $\mathcal{H}^{\textnormal{stub}}_{s,d,m}(\vb*{d})$ is a strongly connected graph.
\\
\begin{lemma}
\label{lemma:stub_sdm_connected}
    $\mathcal{G}(\mathcal{H}^{\textnormal{stub}}_{s,d,m}(\vb*{d}))$ is a strongly connected graph, for any degree sequence $\vb*{d}$.
\end{lemma}

Before showing the proof, we observe that there exists a path in $\mathcal{G}(\mathcal{H}^{\textnormal{stub}}_{s,d,m}(\vb*{d}))$ between any two stub-labeled realizations $H_1,H_2$ of some vertex-labeled directed hypergraph $H'$: Sequentially, the stubs can be rearranged using the double hyperarc shuffle. Thus, there exists a path between two stub-labeled directed hypergraphs $H_0$ and $H^*$ if and only if there exists a path between two stub-labeled directed hypergraphs $H'_0$ and $H^{\prime *}$, where $H_0$ and $H'_0$ are two stub-labeled realizations of the same vertex-labeled directed hypergraph, and $H^*$ and $H^{\prime *}$ are two stub-labeled realizations of the same vertex-labeled directed hypergraph. Therefore, we may neglect the stub labels in the proofs of connectivity.

\begin{proof}[Proof of Lemma \ref{lemma:stub_sdm_connected}]
Let $\mathcal{G}(\mathcal{H}^{\textnormal{stub}}_{s,d,m}(\vb*{d})) = (\mathcal{V},\mathcal{A})$. Let $H_0,H^* \in \mathcal{V}$ be two arbitrary vertices of $\mathcal{V}$, where $H_0=(V,A_0), H^*=(V,A^*)$. We show that there exists a path in $\mathcal{G}(\mathcal{H}^{\textnormal{stub}}_{s,d,m}(\vb*{d}))$ connecting $H_0$ to $H^*$.
    
Since the degree sequences of $H_0$ and $H^*$ are equal, for every hyperarc $a \in A_0$ there must exist a distinct hyperarc $\overline{a} \in A^*$ with $(\delta_a^{\tail},\delta_a^{\head})=(\delta_{\overline{a}}^{\tail},\delta_{\overline{a}}^{\head})$. We introduce a bijection $a \mapsto \overline{a}$ from $A_0$ to $A^*$ such that $\{\overline{a}: a \in A_0\} = A^*$. The arc $\overline{a}$ is considered the `target arc' of $a$: if $a$ is transformed to $\overline{a}$ for all $a \in A_0$ then $H_0$ is transformed to $H^*$. 

We quantify the difference between the two directed hypergraphs by counter $\mathcal{E}$:
\begin{align*}
    \mathcal{E}(H_0,H^*) = \sum_{a \in A(H_0)} (|\overline{a}^{\tail} \backslash a^{\tail}| + |\overline{a}^{\head} \backslash a^{\head}|).
\end{align*}
By construction, $\mathcal{E}(H_0,H^*) \in \mathds{N} \backslash \{1\}$. We prove the lemma by induction on the counter value.

\underline{Base cases}: $\mathcal{E}(H_0,H^*)=0$ and $\mathcal{E}(H_0,H^*)=2$\\
If $\mathcal{E}(H_0,H^*)=0$ then $A(H_0)=A(H^*)$, so $H_0=H^*$ and the states $H_0$ and $H^*$ are connected in $\mathcal{G}(\mathcal{H}^{\textnormal{stub}}_{s,d,m}(\vb*{d}))$, since every state is connected to itself. If $\mathcal{E}(H_0,H^*)=2$ then there must be two hyperarcs $a_1,a_2 \in A_0$ with target hyperarcs $\overline{a}_1$ and $\overline{a}_2$ and two vertices $v,w \in V$ for which holds $v \in a_1^* \cap \overline{a}_2^*, v \notin \overline{a}_1^* \cup a_2^*,w \notin a_1^* \cup \overline{a}_2^*, w \in \overline{a}_1^* \cap a_2^* $ for some $* \in \{\tail, \head\}$. W.l.o.g., let $*=\tail$. Now, the directed hypergraph $\Tilde{H} = (V,\Tilde{A}), A=(A(H_0) \backslash \{a_1,a_2\}) \uplus \{a_3,a_4\}$ with
\begin{align*}
    a_3 &= ((a_1^{\tail} \backslash \{v\}) \uplus \{w\}, a_1^{\head})\\
    a_4 &= ((a_2^{\tail} \backslash \{w\}) \uplus \{v\}, a_2^{\head}) 
\end{align*}
has $\Tilde{H} \in S(a_1,a_2|H_0)$ and $\Tilde{H} = H^*$, so $H_0$ and $H^*$ are adjacent in $\mathcal{G}(\mathcal{H}^{\textnormal{stub}}_{s,d,m}(\vb*{d}))$.

\underline{Induction hypothesis}: There exists a $k \in \mathds{N} \backslash \{0,1\}$ such that there exists a path from $H_0$ to $H^*$ in $\mathcal{G}(\mathcal{H}^{\textnormal{stub}}_{s,d,m}(\vb*{d}))$ if $\mathcal{E}(H_0,H^*) \leq k$.

\underline{To prove}: There exists a path from $H_0$ to $H^*$ in $\mathcal{G}(\mathcal{H}^{\textnormal{stub}}_{s,d,m}(\vb*{d}))$ if $\mathcal{E}(H_0,H^*) = k+1$.

Assume $\mathcal{E}(H_0,H^*) = k+1$. We prove the existence of a directed hypergraph $\Tilde{H} \in \mathcal{H}^{\textnormal{stub}}_{s,d,m}(\vb*{d})$ that is adjacent to $H_0$ in $\mathcal{G}(\mathcal{H}^{\textnormal{stub}}_{s,d,m}(\vb*{d}))$ and has $\mathcal{E}(\Tilde{H},H^*) \leq k$. 

Let $a_1 \in A(H_0)$ be some hyperarc with target $\overline{a}_1 \in A(H^*)$, where $a_1 \neq \overline{a}_1$. Since $a_1 \neq a_1^*$, there must exist a $v \in a_1^{\tail}$ with $v \notin a_1^{* \tail}$, and/or a $w \in a_1^{\head}$ with $w \notin a_1^{*\head}$. W.l.o.g., assume the former. Then, there must also exist a hyperarc $a_2 \in A(H_0)$ with target $\overline{a}_2 \in A(H^*)$ and $v \notin a_2^{\tail}, v \in a_2^{* \tail}$. Since $\delta_{a_2}^{\tail}=\delta_{\overline{a}_2}^{\tail}$, there also exists a $w \in a_2^{\tail}, w \notin a_2^{* \tail}$. 
Now, the directed hypergraph $\Tilde{H} = (V,\Tilde{A})$, $\Tilde{A}=(A(H_0) \backslash \{a_1,a_2\}) \uplus \{a_3,a_4\}$ with
    \begin{align*}
        a_3 &= ((a_1^{\tail} \backslash \{v\}) \uplus \{w\}, a_1^{\head})\\
        a_4 &= ((a_2^{\tail} \backslash \{w\} )\uplus \{v\}, a_2^{\head})        
    \end{align*}
has $\Tilde{H} \in S(a_1,a_2 | H_0)$. Since the double hyperarc shuffle preserves the degree sequence, we obtain that $\Tilde{H} \in \mathcal{H}^{\textnormal{stub}}_{s,d,m}(\vb*{d})$. 
In addition, 
\begin{align*}
    |a_3^{* \tail} \backslash a_3^{\tail}| &\leq |a_1^{* \tail} \backslash a_1^{\tail}|,\\
    |a_4^{* \tail} \backslash a_4^{\tail}| &= |a_2^{* \tail} \backslash a_2^{\tail}| -1
\end{align*}
and so
\begin{align*}
    \mathcal{E}(\Tilde{H},H^*) &= \sum_{a \in A(\Tilde{H})} (|a^{* \tail} \backslash a^{\tail}| + |a^{* \head} \backslash a^{\head}|) \\
    &= \sum_{a \in A(H_0) \backslash \{a_1,a_2\}} (|a^{* \tail} \backslash a^{\tail}| + |a^{* \head} \backslash a^{\head}|) + |a_3^{* \tail} \backslash a_3^{\tail}| + |a_3^{* \head} \backslash a_3^{\head}| + |a_4^{* \tail} \backslash a_4^{\tail}| + |a_4^{* \head} \backslash a_4^{\head}|\\
    &\leq \sum_{a \in A(H_0) \backslash \{a_1,a_2\}} (|a^{* \tail} \backslash a^{\tail}| + |a^{* \head} \backslash a^{\head}|) + |a_1^{* \tail} \backslash a_1^{\tail}| + |a_1^{* \head} \backslash a_1^{\head}| + |a_2^{* \tail} \backslash a_2^{\tail}|-1 + |a_2^{* \head} \backslash a_2^{\head}|\\
    &= \sum_{a \in A(H_0)} (|a^{* \tail} \backslash a^{\tail}| + |a^{* \head} \backslash a^{\head}|) -1\\
    &= k.
\end{align*}

By construction, $\Tilde{H} \in S(a_1,a_2|H_0)$, so $H_0$ and $\Tilde{H}$ are adjacent in $\mathcal{G}(\mathcal{H}^{\textnormal{stub}}_{s,d,m}(\vb*{d}))$. By the induction hypothesis, there exists a path from $\Tilde{H}$ to $H^*$ in $\mathcal{G}(\mathcal{H}^{\textnormal{stub}}_{s,d,m}(\vb*{d}))$.
\end{proof}

\subsection[Strong connectivity of G(H(sm))]{Strong connectivity of $\mathcal{G}(\mathcal{H}^{\textnormal{stub}}_{s,m}(\vb*{d}))$}
\label{section:strong_connectivity_H_sm}
Now, we show that the graph of graphs of $\mathcal{H}^{\textnormal{stub}}_{s,m}(\vb*{d})$ is a strongly connected graph.
\\
\begin{lemma}
\label{lemma:stub_sm_connected}
    $\mathcal{G}(\mathcal{H}^{\textnormal{stub}}_{s,m}(\vb*{d}))$ is a strongly connected graph, for any degree sequence $\vb*{d}$.
\end{lemma}

Although $\mathcal{G}(\mathcal{H}^{\textnormal{stub}}_{s,m}(\vb*{d})) \subseteq \mathcal{G}(\mathcal{H}^{\textnormal{stub}}_{s,d,m}(\vb*{d}))$, we cannot use Lemma \ref{lemma:stub_sdm_connected} to prove Lemma \ref{lemma:stub_sm_connected}. Lemma \ref{lemma:stub_sdm_connected} tells us that any two states $H,H' \in \mathcal{H}^{\textnormal{stub}}_{s,d,m}(\vb*{d})$ are connected via a path $H, H^{(1)}, H^{(2)}, \hdots, H^{(k)}, H'$ in $\mathcal{H}^{\textnormal{stub}}_{s,d,m}(\vb*{d})$. However, this path does not exist in $\mathcal{H}^{\textnormal{stub}}_{s,m}(\vb*{d})$ if any of the states $H^{(1)}, H^{(2)}, \hdots, H^{(k)}$ does not exist in  $\mathcal{H}^{\textnormal{stub}}_{s,m}(\vb*{d})$, i.e., if it represents a directed hypergraph with at least one degenerate hyperarc.

Before proving Lemma \ref{lemma:stub_sm_connected}, we introduce a mapping $f$ from directed hypergraphs to directed bipartite graphs, and we define structural zeros and alternating rectangles in adjacency matrices of such directed bipartite graphs. Later, we will use a theorem on sequentially transforming an adjacency matrix into another adjacency matrix with the same marginals.

Let $f$ be the following map.
\\
\begin{definition}[Map $f$]
    Let $H=(V,A)$ be a directed hypergraph. $f$ maps every hyperarc $(a^{\tail},a^{\head})$ to $\delta_a^{\tail}+\delta_a^{\head}$ arcs resulting in digraph $f(H)=G^{\tail} \cup G^{\head}$, where $G^{\tail}=(U^{\tout} \cup U^{\tail},A^{\tail})$ is a bipartite digraph representing the incidences of tails of hyperarcs in $H$ and $G^{\head}=(U^{\tin} \cup U^{\head},A^{\head})$ is a bipartite digraph representing the incidences of heads of hyperarcs in $H$. In particular,
    \begin{align*}
        U^{\tout} &= \{u^{\tout}_v \, : \, v \in V\},\\
        U^{\tin} &= \{u^{\tin}_v \, : \, v \in V\},\\
        U^{\tail} &= \{u^{\tail}_a \, : \, a \in A\},\\
        U^{\head} &= \{u^{\head}_a \, : \, a \in A\},\\
        A^t &= \{(u^{\tout}_v,u^{\tail}_a): v \in a^{\tail} \textnormal{ and } v \in V, a \in A\},\\
        A^h &= \{(u^{\tin}_v,u^{\head}_a): v \in a^{\head} \textnormal{ and } v \in V, a \in A\}.
    \end{align*}
\end{definition}

\begin{figure}[tbp]
\centering
\begin{subfigure}[t]{0.3\textwidth}
\centering
    \includegraphics[width=0.6\textwidth]{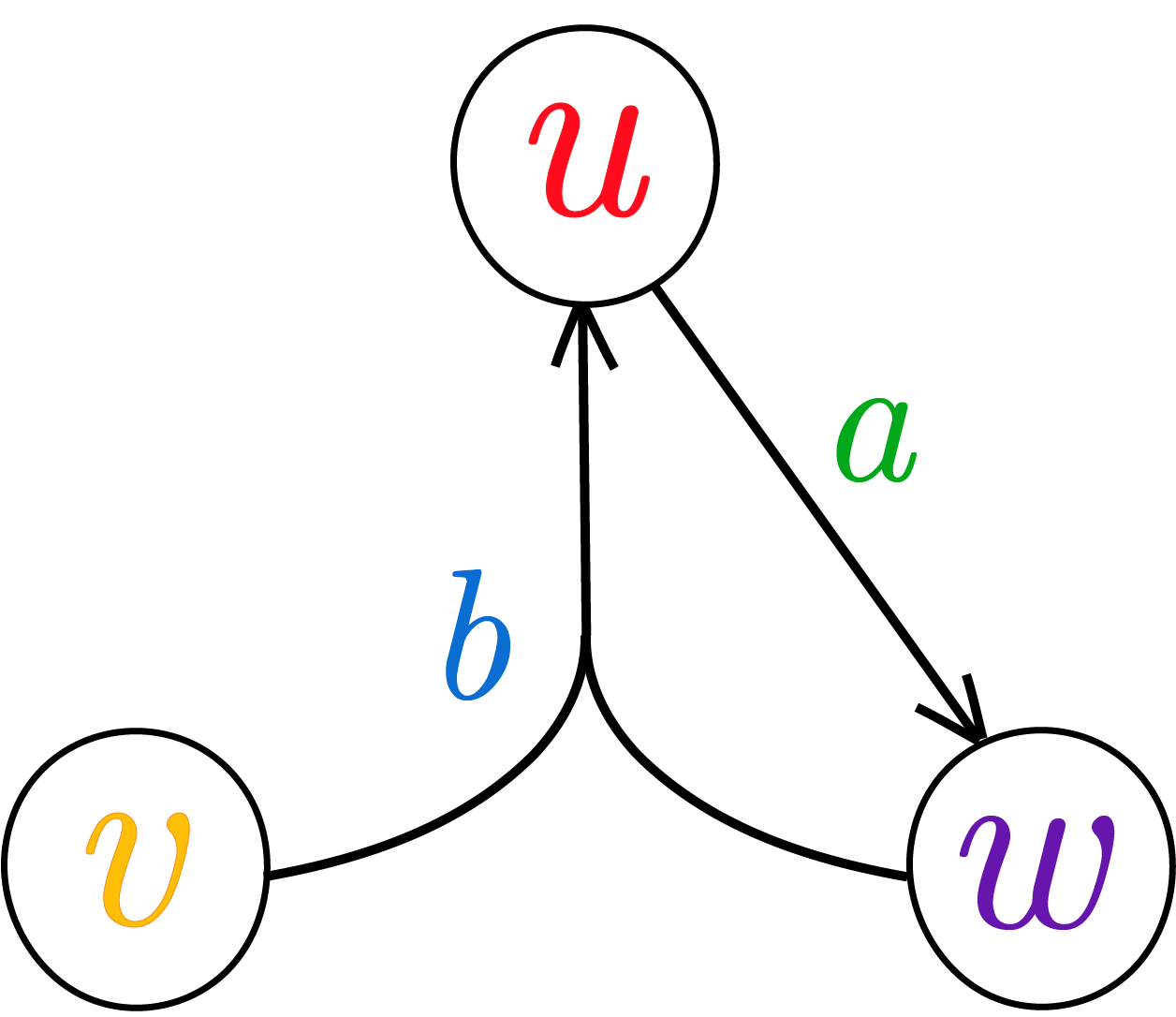}
    \caption{A directed hypergraph $H$.}
\end{subfigure}
\begin{subfigure}[t]{0.65\textwidth}
\centering
    \includegraphics[width=0.65\textwidth]{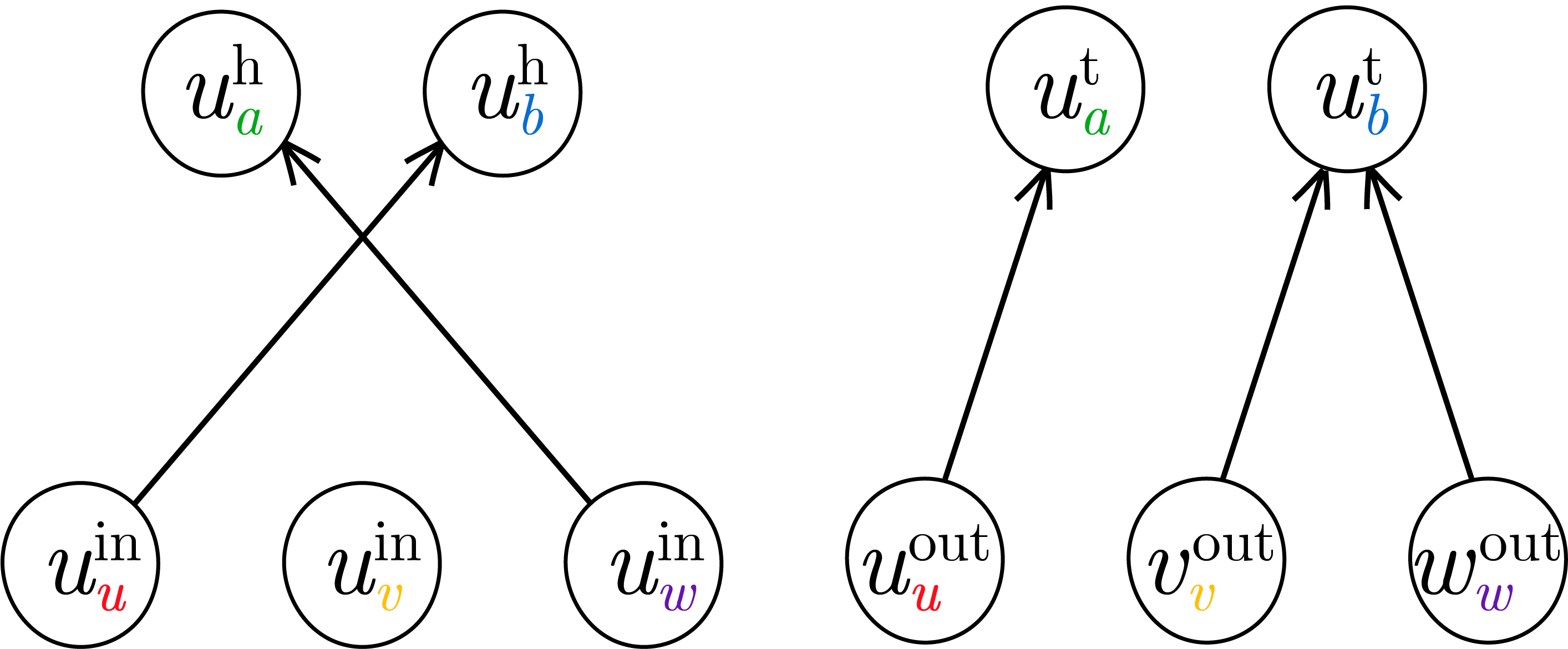}
    \caption{The digraph $f(H)$, where $U^{\tout}=\{u^{\tout}_u,u^{\tout}_v,u^{\tout}_w\}$, $U^{\tail}=\{u^{\tail}_a,u^{\tail}_b\}$, $U^{\tin}=\{u^{\tin}_v,u^{\tin}_v,u^{\tin}_w\}$, $U^{\head}=\{u^{\head}_a,u^{\head}_b\}$.}
\end{subfigure}        
\caption{A directed hypergraph $H \in \mathcal{H}^{\textnormal{stub}}_{s,m}(\vb*{d})$ and its mapping $f(H)$.}
\label{fig:mapping}
\end{figure}

The mapping $f$ is illustrated in Figure \ref{fig:mapping}. This mapping is one-to-one, so the inverse of $f$, $f^{-1}$, is well-defined.

For $f$ it holds that $f: \mathcal{H}^{\textnormal{stub}}_{x}(\vb*{d}) \rightarrow \mathcal{B}$, where $\mathcal{B}$ is the space of bipartite digraphs with vertex set $U^{\tin} \cup U^{\tout} \cup U^{\head} \cup U^{\tail}$ and arcs from $U^{\tin}$ to $U^{\head}$ and from $U^{\tout}$ to $U^{\tail}$. Moreover, $|U^{\tin}| = |U^{\tout}|$ and $|U^{\head}| = |U^{\tail}|$. A degenerate hyperarc in $H$ maps to a multi-arc in $f(H)$. Therefore, 
\begin{align}
\label{equivalence:H_sp-B*}
    H \in \mathcal{H}^{\textnormal{stub}}_{s,m}(\vb*{d}) \iff f(H) \in \mathcal{B}^*,
\end{align}
where $\mathcal{B}^*=\{B \in \mathcal{B}: B \textnormal{ contains no multi-arcs}\}$.

Next, we denote with $M_{\head}$ the adjacency matrix of $G^{\head}$ and with $M_{\tail}$ the adjacency matrix of $G^{\tail}$. Since graphs in $\mathcal{B}^*$ contain no multi-arcs, their adjacency matrices are (0,1)-matrices. Therefore, also $M_{\head}$ and $M_{\tail}$ are (0,1)-matrices. Moreover, by construction of $G^{\head}$, some entries of $M_{\head}$ always need to be zeros. An example is the entry that represents the connection of a vertex in $U^{\tout}$ to itself. The same holds for $M_{\tail}$. Such an entry is called a \textit{structural zero}.
\\
\begin{definition}[Structural zero]
    A structural zero in an adjacency matrix is a zero entry that is not allowed to be changed to a one, by construction of the graph represented by the adjacency matrix.
\end{definition}

In $M_{\head}$, the elements $ij$ with $i \in U^{\head}$ or $j \in U^{\tin}$ are the only structural zeros, since all arcs need to go from vertices in $U^{\tin}$ to vertices in $U^{\head}$ (by construction of $\mathcal{B}$).

Next, we introduce alternating rectangles in these adjacency matrix. An alternating rectangle represents two arcs $(i_1,j_1)$ and $(i_2,j_2)$ that can be shuffled to become $(i_1,j_2)$ and $(i_2,j_1)$, without creating multi-arcs, which later will be used to transform one directed hypergraph into another.
\\
\begin{definition}[Alternating rectangle (\cite{rao1996})]
    An alternating rectangle in the (0,1)-matrix $M$ is a set of four distinct cells of the type $\{i_1j_1,i_1j_2,i_2j_2,i_2j_1\}$ in $M$ such that the entries in the four cells are alternately zeros and ones as one goes around the rectangle in either direction and the zeros are not structural zeros.
\end{definition}

In $M_{\head}$, observe that any set of four distinct cells of the type $\{i_1j_1,i_1j_2,i_2j_2,i_2j_1\}$ such that the entries in the four cells are alternately zeros and ones is an alternating rectangle, as those zeros are never structural zeros.
\\
\begin{claim}
\label{claim:nostructural0s}
    In $M_{\head}$ and $M_{\tail}$, any set of four distinct cells of the type $\{i_1j_1,i_1j_2,i_2j_2,i_2j_1\}$ such that the entries in the four cells are alternately zeros and ones as one goes around the rectangle in either direction contains no structural zeros, and is therefore an alternating rectangle.
\end{claim}
\begin{proof}
    The proof is by contradiction. Let $\{i_1j_1,i_1j_2,i_2j_2,i_2j_1\}$ be such a set in $M_{\tail}$. W.l.o.g., let $i_1j_1=0$ be a structural zero. Then $i_1 \in U^{\head}$ or $j_1 \in U^{\tin}$. Moreover, $i_1j_2=1,i_2j_2=0,i_2j_1=1$.\\
    Case 1: $i_1 \in U^{\head}$. Then $i_1j_2=1$ is not possible, since arcs starting from vertices in $U^{\head}$ are not allowed in $\mathcal{B}$, by construction.\\
    Case 2: $j_1 \in U^{\tin}$. Then $i_2j_1=1$ is not possible, since arcs ending at vertices in $U^{\tin}$ are not allowed in $\mathcal{B}$, by construction.\\
    In conclusion, such a set cannot exist. By symmetry, the proof holds for $M_{\head}$ as well.
\end{proof}

Lastly, we introduce \textit{switching along an alternating rectangle}, which we will later relate to performing a double hyperarc shuffle.
\\
\begin{definition}[Switching along an alternating rectangle (\cite{rao1996})]
    Let $\{i_1j_1,i_1j_2,i_2j_2,i_2j_1\}$ be an alternating rectangle. Switching along this alternating rectangle is interchanging the zero and one entries. 
\end{definition}

Switching along an alternating rectangle transforms a graph in $\mathcal{B}^*$ into another graph in $\mathcal{B}^*$.
\\
\begin{claim}
\label{claim:switch_stays_in_B*}
    For any graph $G \in \mathcal{B}^*$ with adjacency matrix $M$, switching along an alternating rectangle in $M$ results in an adjacency matrix $M'$ that corresponds to a graph $G' \in \mathcal{B}^*$ with the same degree sequence as $G$.
\end{claim}
\begin{proof}
    Switching along an alternating rectangle in $M$ results in an adjacency matrix $M'$, which, by construction of the alternating rectangle, is a (0,1)-matrix. Therefore, the corresponding graph $G'$ contains no multi-arcs. Thus, by (\ref{equivalence:H_sp-B*}), $G' \in \mathcal{B}^*$. Moreover, the marginals of $M$ and $M'$ are equal, so $G$ and $G'$ have the same degree sequence.
\end{proof}

We now show that switching along an alternating rectangle is equivalent to performing a double hyperarc shuffle on the corresponding directed hypergraph.
\\
\begin{claim}
\label{claim:switching_equals_shuffle}
    Let $H \in \mathcal{H}^{\textnormal{stub}}_{s,m}(\vb*{d})$ be a directed hypergraph with mapping $f(H)=G^{\head} \cup G^{\tail}$, with adjacency matrices $M_{\head}$ and $M_{\tail}$ of $G^{\head}$ and $G^{\tail}$, respectively. Switching along an alternating rectangle $\{i_1j_1,i_1j_2,i_2j_2,i_2j_1\}$ in $M_{\head}$ ($M_{\tail}$) results in $M'_{\head}$ ($M'_{\tail}$), representing $G^{\prime \head}$ ($G^{\prime \tail}$). Let $H'=f^{-1}(G^{\prime \head} \cup G^{\tail})$ ($H'=f^{-1}(G^{ \head} \cup G^{\prime \tail})$). Then, $H' \in S(j_1,j_2|H)$ and $H' \in \mathcal{H}^{\textnormal{stub}}_{s,m}(\vb*{d})$.
\end{claim}
\begin{proof}
    By symmetry, we only need to prove the claim for $M_{\head}$. The adjacency matrix $M_{\head}$ denotes the incidences of the heads of the hyperarcs and the vertices in $H$. Since the alternating rectangle $\{i_1j_1,i_1j_2,i_2j_2,i_2j_1\}$ exists in $M_{\head}$, for the hyperarcs $j_1,j_2 \in A(H)$ holds $i_1 \in j_1^{\head}$, $i_2 \notin j_1^{\head}$, $i_1 \notin j_2^{\head}$ and $i_2 \in j_2^{\head}$. The switch along the alternating rectangle results in the adjacency matrix $M'_{\head}$, which denotes the incidences of the heads of the hyperarcs and the vertices in some other directed hypergraph $H'$, where $A(H')=(A(H) \backslash \{j_1,j_2\}) \uplus \{\hat{j}_1,\hat{j}_2\}$, where $\hat{j}_1 = (j_1^{\tail},(j_1^{\head} \backslash \{i_1\}) \uplus \{i_2\})$ and $\hat{j}_2 = (j_2^{\tail},(j_2^{\head} \backslash \{i_2\}) \uplus \{i_1\})$. This gives $H' \in S(j_1,j_2|H)$. By claim \ref{claim:switch_stays_in_B*}, $G^{\prime \head} \in \mathcal{B}^*$, so $H' \in \mathcal{H}^{\textnormal{stub}}_{s,m}(\vb*{d})$.
\end{proof}

Now we are ready to prove Lemma \ref{lemma:stub_sm_connected}. 

\begin{proof}[Proof of Lemma \ref{lemma:stub_sm_connected}]
Let $H_1=(V,A_1),H_2=(V,A_2) \in \mathcal{H}^{\textnormal{stub}}_{s,m}(\vb*{d})$ be two directed hypergraphs. We show that there exists a path from $H_1$ to $H_2$ in $\mathcal{G}(\mathcal{H}^{\textnormal{stub}}_{s,m}(\vb*{d}))$. The proof uses mapping $f$ and a theorem about strong connectivity of the graph of graphs of (0,1)-matrices.

Let $f(H_1) = G^{\head}_1 \cup G^{\tail}_1$, and $f(H_2) = G^{\head}_2 \cup G^{\tail}_2$. Consider the subgraphs $G^{\head}_1$ and $G^{\head}_2$. By (\ref{equivalence:H_sp-B*}), the adjacency matrices $M_{h,1}$ and $M_{h,2}$ of $G^{\head}_1$ and $G^{\head}_2$ are (0,1)-matrices. Since the degree sequences of $H_1$ and $H_2$ are equal, the marginals of $M_{h,1}$ and $M_{h,2}$ are equal. By Theorem 1 in \cite{rao1996}, there exists a sequence of switchings along alternating rectangles that transforms $M_{h,1}$ into $M_{h,2}$. By Claim \ref{claim:nostructural0s}, none of these switchings use structural zeros. By Claim \ref{claim:switching_equals_shuffle}, the sequence of switchings corresponds to a sequence of double hyperarc shuffles that transforms $H_1$ into $H'_1$, where $H_1'$ is the directed hypergraph with hyperarc set $A_1'$ in which the hyperarc heads equal the hyperarc heads in $A_2$ and the hyperarc tails equal the hyperarc tails in $A_1$. Therefore, there exists a path from $H_1$ to $H'_1$ in $\mathcal{G}(\mathcal{H}^{\textnormal{stub}}_{s,m}\vb*{d}))$.

In a similar way, consider $G^{\tail}_1$ and $G^{\tail}_2$. Theorem 1 in \cite{rao1996} again tells us that there exists a sequence of switchings along alternating rectangles that transforms the adjacency matrix of $G^{\tail}_1$ into the adjacency matrix of $G^{\tail}_2$. This sequence of double hyperarc shuffles corresponds to a path from $H'_1$ to $H_2$ in $\mathcal{G}(\mathcal{H}^{\textnormal{stub}}_{s,m}(\vb*{d}))$.  

In conclusion, this proves existence of a path in $\mathcal{G}(\mathcal{H}^{\textnormal{stub}}_{s,m}(\vb*{d}))$ from $H_1$ to $H_2$.
\end{proof}

We are now ready to prove Theorem \ref{thm:uniform_alldegrees}.
\begin{proof}[Proof of Theorem \ref{thm:uniform_alldegrees}]
    Since $\mathcal{G}(\mathcal{H}^{\textnormal{stub}}_{s,d,m}(\vb*{d}))$ and $\mathcal{G}(\mathcal{H}^{\textnormal{stub}}_{s,m}(\vb*{d}))$ are aperiodic (Lemma \ref{lemma:aperiodic}) and strongly connected (Lemmas \ref{lemma:stub_sdm_connected} and \ref{lemma:stub_sm_connected}, respectively), random walks on $\mathcal{G}(\mathcal{H}^{\textnormal{stub}}_{s,d,m}(\vb*{d}))$ and $\mathcal{G}(\mathcal{H}^{\textnormal{stub}}_{s,m}(\vb*{d}))$ are ergodic. Since $\mathcal{G}(\mathcal{H}^{\textnormal{stub}}_{s,d,m}(\vb*{d}))$ and $\mathcal{G}(\mathcal{H}^{\textnormal{stub}}_{s,m}(\vb*{d}))$ are also regular (Lemma \ref{lemma:regular}), both have a uniform stationary distribution.
\end{proof}

\subsection{Proof of Theorem \ref{thm:notuniform}}
\label{section:no_connectivity}
We show that $\mathcal{G}(\mathcal{H}^{\textnormal{stub}}_{}(\vb*{d}_2))$, $\mathcal{G}(\mathcal{H}^{\textnormal{stub}}_{d}(\vb*{d}_2))$, $\mathcal{G}(\mathcal{H}^{\textnormal{stub}}_{m}(\vb*{d}_2))$, $\mathcal{G}(\mathcal{H}^{\textnormal{stub}}_{d,m}(\vb*{d}_2))$, and $\mathcal{G}(\mathcal{H}^{\textnormal{stub}}_{s,d}(\vb*{d}_1))$ are not strongly connected, for some degree sequences $\vb*{d}_1,\vb*{d}_2$.
\\
\begin{lemma}
\label{lemma:stub_sd_notconnected}
   There exists a degree sequence $\vb*{d}_1$ such that $\mathcal{G}(\mathcal{H}^{\textnormal{stub}}_{s,d}(\vb*{d}_1))$ is not a strongly connected graph.
\end{lemma}
\begin{proof}
Consider the space $\mathcal{H}^{\textnormal{stub}}_{s,d}(\vb*{d}_1)$ with $\vb*{d}_1=(\vb*{d}_V,\vb*{\delta}_A)$ and $\vb*{d}_V=((0,2),(0,2),(0,2),(3,0))$, $\vb*{\delta}_A = \{(2,1),(2,1),(2,1)\}$. We present two directed hypergraphs $H_0,H^* \in \mathcal{H}^{\textnormal{stub}}_{s,d}(\vb*{d}_1)$ between which no path exists in $\mathcal{G}(\mathcal{H}^{\textnormal{stub}}_{s,d}(\vb*{d}_1))$. Let $H_0=(V,A_0)$ with $V=\{u,v,w,x\}$, \linebreak $A_0=\{(\{u,u\},\{x\}),(\{v,v\},\{x\}),(\{w,w\},\{x\})\}$ (Figure \ref{fig:stub_sd_notconnectedH}). Let $H^*=(V,A^*)$, with $A^*=\{(\{u,v\},\{x\}),(\{u,w\},\{x\}),(\{v,w\},\{x\})\}$ (Figure \ref{fig:stub_sd_notconnectedH'}).

\begin{figure}[h!]
\centering
\begin{subfigure}{0.4\textwidth}
    \centering
    \includegraphics[width=0.5\textwidth]{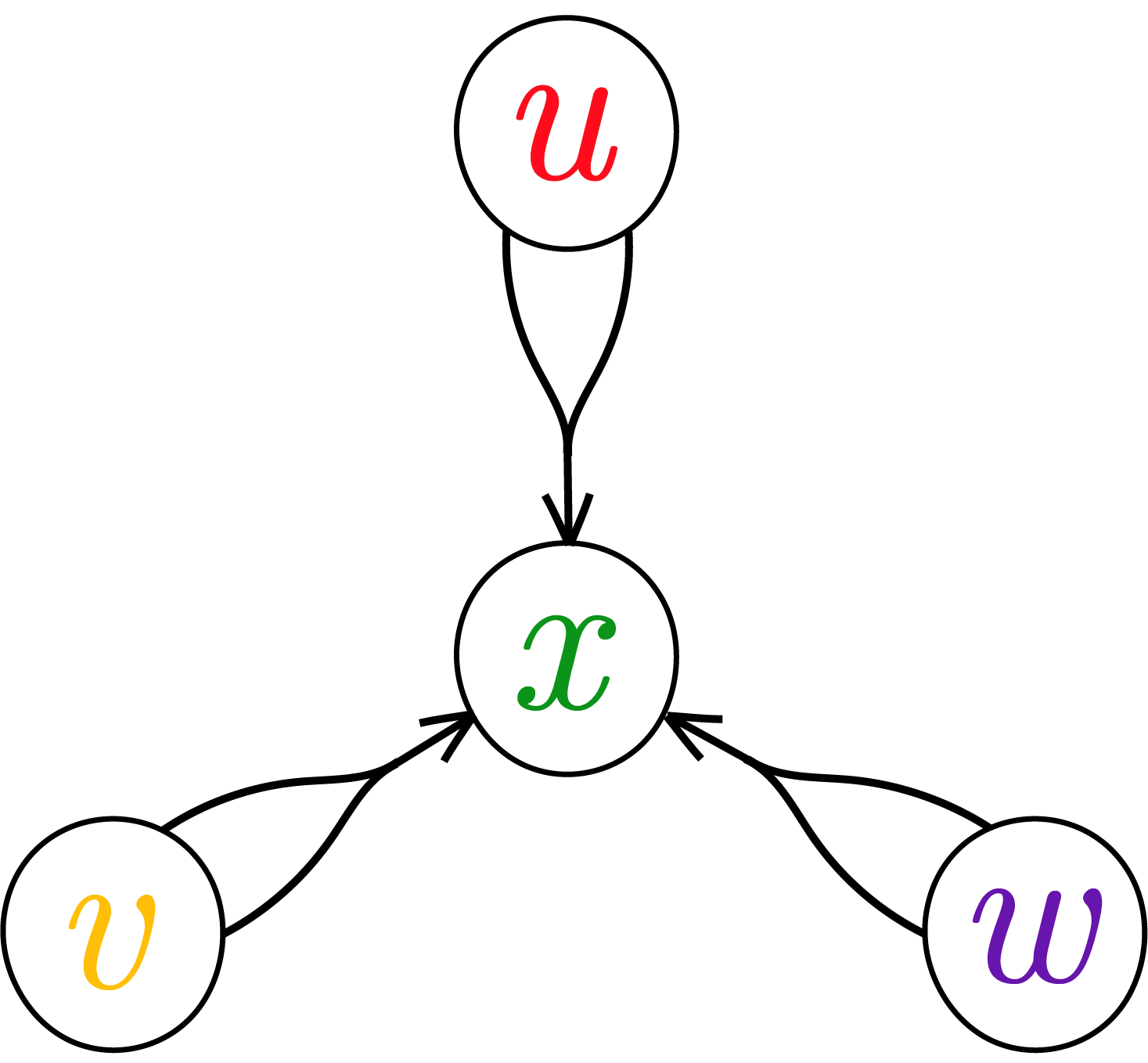}
    \caption{The directed hypergraph $H_0$.}
    \label{fig:stub_sd_notconnectedH}
\end{subfigure}
\hspace{1cm}
\begin{subfigure}{0.4\textwidth}
    \centering
    \includegraphics[width=0.5\textwidth]{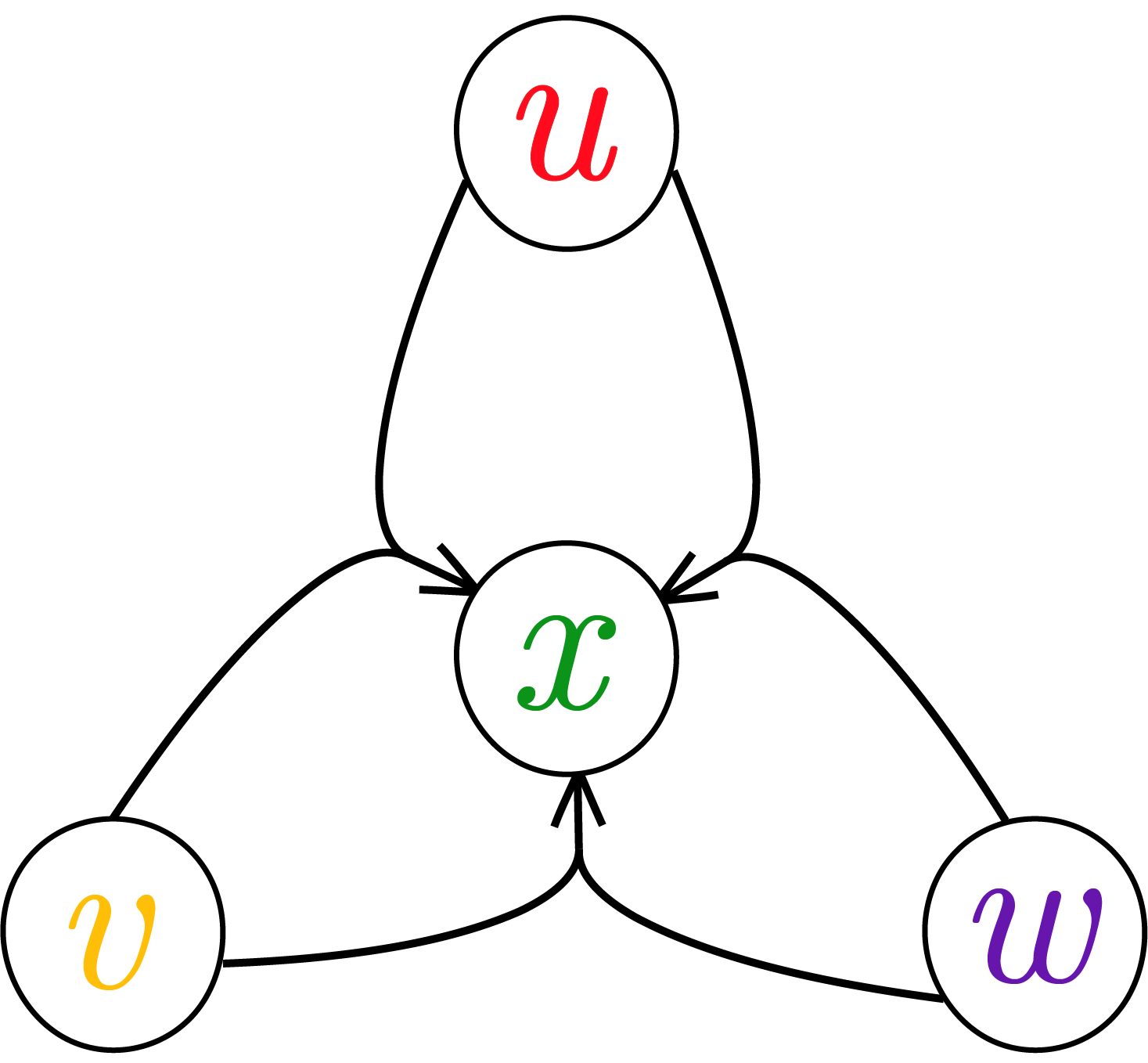}
    \caption{The directed hypergraph $H^*$.}
    \label{fig:stub_sd_notconnectedH'}
\end{subfigure}        
\caption{Two directed hypergraphs $H_0,H^* \in \mathcal{H}^{\textnormal{stub}}_{s,d}(\vb*{d}_1)$ between which there exists no path in $\mathcal{G}(\mathcal{H}^{\textnormal{stub}}_{s,d}(\vb*{d}_1))$.}
\end{figure}

We now show that there exists no path from $H_0$ to $H^*$ in $\mathcal{G}(\mathcal{H}^{\textnormal{stub}}_{s,d}(\vb*{d}_1))$. To that end, consider any two hyperarcs $a_1,a_2 \in A_0$, so $a_1=(\{m,m\},\{x\}),a_2=(\{n,n\},\{x\})$, with $m,n \in \{u,v,w\}$. The double hyperarc shuffle $s(a_1,a_2|H_0)$ can only result in the unchanged directed hypergraph $H_0$, or in a new directed hypergraph with the multi-hyperarcs pair $\hat{a}_1=\hat{a}_2=(\{m,n\},\{x\})$, which is not allowed. Therefore, $S(a_1,a_2|H_0)=\{H_0\}$ for all $a_1,a_2 \in A(H_0)$, so $H_0$ is an isolated vertex in $\mathcal{G}(\mathcal{H}^{\textnormal{stub}}_{s,d}(\vb*{d}_1))$. In particular, there exists no path from $H_0$ to $H^*$ in $\mathcal{G}(\mathcal{H}^{\textnormal{stub}}_{s,d}(\vb*{d}_1))$. 
\end{proof}
\begin{lemma}
\label{lemma:stub_many_notconnected}
   $\mathcal{G}(\mathcal{H}^{\textnormal{stub}}_{}(\vb*{d}_2))$, $\mathcal{G}(\mathcal{H}^{\textnormal{stub}}_{d}(\vb*{d}_2))$, $\mathcal{G}(\mathcal{H}^{\textnormal{stub}}_{m}(\vb*{d}_2))$ and $\mathcal{G}(\mathcal{H}^{\textnormal{stub}}_{d,m}(\vb*{d}_2))$ are no strongly connected graphs, for at least one degree sequence $\vb*{d}_2$.
\end{lemma}
\begin{proof}
Let $\vb*{d}_2=(\vb*{d}_V,\vb*{d}_A)$. If $\vb*{d}_A=((1,1))_{a \in A}$, then the spaces $\mathcal{H}^{\textnormal{stub}}_{d,m}(\vb*{d}_2)$ and $\mathcal{H}^{\textnormal{stub}}_{m}(\vb*{d}_2)$ reduce to the space of digraphs without self-loops and the spaces $\mathcal{H}^{\textnormal{stub}}_{d}(\vb*{d}_2)$ and $\mathcal{H}^{\textnormal{stub}}_{}(\vb*{d}_2)$ reduce to the space of simple digraphs. For the former space, there exists a degree sequence for which the graph of graphs is not strongly connected (\cite{carstens2017}), and for the latter space as well (\cite{rao1996}, \cite{berger2010}).
\end{proof}

We can now prove Theorem \ref{thm:notuniform}. 
\begin{proof}[Proof of Theorem \ref{thm:notuniform}]
    Since there exist degree sequences $\vb*{d}_1,\vb*{d}_2$ for which
    $\mathcal{G}(\mathcal{H}^{\textnormal{stub}}_{s,d}(\vb*{d}_1))$, $\mathcal{G}(\mathcal{H}^{\textnormal{stub}}_{d,m}(\vb*{d}_2))$, $\mathcal{G}(\mathcal{H}^{\textnormal{stub}}_{d}(\vb*{d}_2))$, $\mathcal{G}(\mathcal{H}^{\textnormal{stub}}_{m}(\vb*{d}_2))$ and $\mathcal{G}(\mathcal{H}^{\textnormal{stub}}_{}(\vb*{d}_2))$ are no strongly connected graphs (Lemmas \ref{lemma:stub_sd_notconnected} and \ref{lemma:stub_many_notconnected}), the double hyperarc shuffle method does not sample uniformly from the spaces $\mathcal{H}^{\textnormal{stub}}_{s,d}(\vb*{d}_1)$, $\mathcal{H}^{\textnormal{stub}}_{d,m}(\vb*{d}_2)$, $\mathcal{H}^{\textnormal{stub}}_{d}(\vb*{d}_2)$, $\mathcal{H}^{\textnormal{stub}}_{m}(\vb*{d}_2)$, and $\mathcal{H}^{\textnormal{stub}}_{}(\vb*{d}_2)$.
\end{proof}

\section{Proof of Theorem \ref{thm:uniform_degrees(2,1)}}
\label{section:proofThm2}
Lemmas \ref{lemma:regular} and \ref{lemma:aperiodic} also hold for $\mathcal{G}(\mathcal{H}^{\textnormal{stub}}_{s}(\vb*{d}))$. It remains to prove strong connectivity. Let
\begin{align*}
    D = \{\vb*{d}: \forall a \in A: \delta_a^{\tail}=1, \delta_a^{\head}=2, |\{v: d_v^{\tout}>0\}|=2\}
\end{align*}
denote all degree sequences that are allowed by Theorem \ref{thm:uniform_degrees(2,1)}.
\\
\begin{lemma}
\label{lemma:stub_s_connected}
    For any degree sequence $\vb*{d} \in D$, $\mathcal{G} (\mathcal{H}^{\textnormal{stub}}_{s}(\vb*{d}))$ is a strongly connected graph.
\end{lemma}

Let the two tail vertices be $\{v: d_v^{\textnormal{out}}>0\}= \{\tau,\sigma\}$. Before stating the proof of Lemma \ref{lemma:stub_s_connected}, we introduce the \textit{symmetric difference} between two directed hypergraphs, a mapping of the symmetric difference and \textit{minimal alternating cycles} in this mapped symmetric difference.

\begin{figure}[tbp]
\centering
\begin{subfigure}{0.4\textwidth}
        \centering
        \begin{tabular}{cc}
            $\underline{A(H_0)}$ & $\underline{A(H^*)}$ \\ 
            $(\tau, \{a,d\})$ & $(\tau, \{a,b\})$ \\
            $(\tau, \{b,c\})$ & $(\sigma, \{c,d\})$ \\
            $(\sigma, \{c,e\})$ & $(\tau, \{c,e\})$ \\
            $(\tau, \{a,c\})$ & $(\tau, \{a,c\})$
        \end{tabular}
    \caption{The hyperarc sets $A(H_0)$ and $A(H^*)$ of two directed hypergraphs $H_0$ and $H^*$ with equal degree sequence.}
    \label{fig:hyperarc sets A(H_0),A(H^*)}
\end{subfigure}
\begin{subfigure}{0.45\textwidth}
        \centering
        \begin{tabular}{c}
            $\underline{H_0 \Delta H^*}$ \\ 
            $(\tau, \{a,d\})$ \\
            $(\tau, \{b,c\})$  \\
            $(\sigma, \{c,e\})$  \\
            $(\tau, \{a,b\})$ \\
            $(\sigma, \{c,d\})$ \\
            $(\tau, \{c,e\})$ \\
        \end{tabular}
    \caption{The symmetric difference of $H_0$ and $H^*$.}
    \label{fig:symmetric difference}
\end{subfigure}\\ 
\hspace{1cm}
\begin{subfigure}{\linewidth}
    \center
    \newcommand{\offsetx}{1}
\newcommand{\offsety}{-0.2}
    
\begin{tikzpicture}[-,>=Stealth, line width=1.3pt,auto,
                    thick,main node/.style={circle,draw, minimum size=\tikznodesize}, inner sep=1pt]

                    %[-,>=Spurple, line width=1.3ptth,auto,node distance=2cm,
                    %thick,main node/.style={circle,draw,minimum size=0.8cm,font=\sffamily\Large\bfseries}]

  \node[main node] (a) at (0,0) {$a$};
  \node[main node] (b) at (2*\tikzscale,0) {$b$};
  \node[main node] (c) at (2*\tikzscale,-2*\tikzscale) {$c$};
  \node[main node] (d) at (0,-2*\tikzscale) {$d$};
  \node[main node] (e) at (4*\tikzscale,-2*\tikzscale) {$e$};

  \path[every node/.style={font=\sffamily\small}]
    (a) edge[purple,line width=1.3pt] node {} (b)
    (b) edge[cyan] node {} (c)
    (a) edge[cyan] node {} (d)
    (c) edge[purple, line width=1.3pt, dashed] node {} (d)
    (c) edge[purple, line width=1.3pt, bend right] node {} (e)
    (e) edge[cyan, bend right, dashed] node {} (c);

  % Legend title
  \node[anchor=south west, font=\bfseries\sffamily] at (5.0 + \offsetx,0.55+\offsety) {Legend};

  % Legend items
  \node at (5.9 + \offsetx,0.3 + \offsety) {: \textcolor{black}{$e \in f_{H_0,\tau}$}};
  \draw[cyan, line width=0.02cm] (4.6 + \offsetx,0.3 + \offsety) -- (5 + \offsetx,0.3 + \offsety);
  \node at (5.9 + \offsetx,-0.2 + \offsety) {: \textcolor{black}{$e \in f_{H_0,\sigma}$}};
  \draw[cyan,dashed, line width=0.02cm] (4.6 + \offsetx,-0.2 + \offsety) -- (5 + \offsetx,-0.2 + \offsety);
  \node at (5.9 + \offsetx,-0.7 + \offsety) {: \textcolor{black}{$e \in f_{H^*,\tau}$}};
  \draw[purple, line width=1.3pt] (4.6 + \offsetx,-0.7 + \offsety) -- (5 + \offsetx,-0.7 + \offsety);
  \node at (5.9 + \offsetx,-1.2 + \offsety) {: \textcolor{black}{$e \in f_{H^*,\sigma}$}};
  \draw[purple, line width=1.3pt, dashed] (4.6 + \offsetx,-1.2 + \offsety) -- (5 + \offsetx,-1.2 + \offsety);

  % Legend rectangle
  \node[draw, rectangle, fit={(4.7 + \offsetx,0.85 + \offsety) (6.7 + \offsetx,-1.4 + \offsety)}, inner sep=5pt] {};

\end{tikzpicture}
    \caption{The mapped symmetric difference of $H_0$ and $H^*$: $M(H_0 \Delta H^*)$, where $f_{H_0,\tau}(S) = \{\{a,d\}_{H_0,\tau}, \{b,c\}_{H_0,\tau}\},f_{H_0,\sigma}(S) = \{\{c,e\}_{H_0,\sigma}\}, f_{H^*,\tau}(S) = \{\{a,b\}_{H^*,\tau}, \{c,e\}_{H^*,\tau}\}$ and $f_{H^*,\sigma}(S) = \{\{c,d\}_{H^*,\sigma}\}$. Equation \ref{eq:counter} gives $\mathcal{E}(M(H_0 \Delta H^*))=4+\frac{1}{2}\cdot 2=5$. The legend in this figure is used in following figures as well.}
    \label{fig:mapped symmetric difference}
\end{subfigure}
\caption{Two directed hypergraphs $H_0,H^* \in \mathcal{H}^{\textnormal{stub}}_{s}(\vb*{d})$, their symmetric difference and the mapping of their symmetric difference.}
\end{figure}
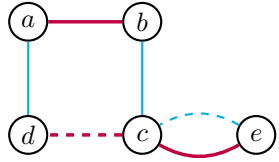
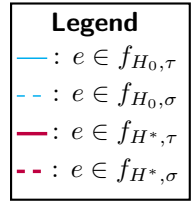

The core concept of the proof is the symmetric difference of two directed hypergraphs. Given two directed hypergraphs $H_0$ and $H^*$, the symmetric difference is defined as the set of hyperarcs that are either in $A(H_0)$ or in $A(H^*)$, but not in both.
\\
\begin{definition}[Symmetric difference]
$H_0 \Delta H^* = (A(H_0) \backslash A(H^*))\cup  (A(H^*) \backslash A(H_0))$
\end{definition}
Figure \ref{fig:symmetric difference} shows the symmetric difference of the directed hypergraphs in Figure \ref{fig:hyperarc sets A(H_0),A(H^*)}.
The symmetric difference of two directed hypergraphs can be mapped to a graph in the following way: 
\\
\begin{definition}[Map $M(H_0 \Delta H^*)$]
    Let $\vb*{d} \in D$ and $H_0,H^* \in \mathcal{H}^{\textnormal{stub}}_{s}(\vb*{d})$. Every hyperarc $(t,\{x,y\})$ in $H_0 \Delta H^*$ is mapped as follows:
    \begin{align*}
       (t,\{x,y\}) \mapsto_M 
        \begin{cases}
            \{x,y\}_{H_0,t} & \textnormal{ if } (t,\{x,y\}) \in A(H_0) \backslash A(H^*)\\
            \{x,y\}_{H^*,t} & \textnormal{ else}.
        \end{cases}
    \end{align*}
   The mapped hyperarc is interpreted as an edge $\{x,y\}$ with attributes $H_0$ and $t$ (resp. $H^*$ and $t$). 
\end{definition}

If $a \mapsto_M e$ then we may write $e=M(a)$ and $a=M^{-1}(e)$.

For an edge $e \in M(H_0 \Delta H^*)$, we refer to the first and second attribute by $e^1$ and $e^2$, respectively: when $e=\{x,y\}_{H_0,t}$, then $e^1=H_0, e^2=t$. We may also refer to an edge $\{a,b\}_{.,t}$ or $\{a,b\}_{H_0,.}$ in the symmetric difference if the first resp. second attribute is irrelevant.

The mapping $M(H_0 \Delta H^*)$ is illustrated in Figure \ref{fig:mapped symmetric difference}, where the first attribute is represented by a cyan edge (for $H_0$) or a purple edge (for $H^*$) and the second attribute is represented by a solid line (for tail $\tau$) or a dashed line (for tail $\sigma$). In the hypergraph space $\mathcal{H}^{\textnormal{stub}}_{s}(\vb*{d})$ no multi-hyperarcs are allowed. Therefore, the symmetric difference of two hypergraphs never contains multiple copies of the same hyperarc. Thus, the mapping contains no multi-edges with the same attributes. Moreover, since no degenerate hyperarcs are allowed, the mapping contains no self-loops. Lastly, the next property demonstrates that any two edges between the same vertices must have a different second attribute.
\\
\begin{property}
\label{prop:second_att}
    Let $\vb*{d} \in D$ and $H_0,H^* \in \mathcal{H}^{\textnormal{stub}}_{s}(\vb*{d})$. If $e_1,e_2\in M(H_0 \Delta H^*)$ with $e_1 = \{x,y\}_{.,s}$ and $e_2 = \{x,y\}_{.,t}$ then $s \neq t$, so that $e_1^2 \neq e_2^2$.
\end{property}
\begin{proof}
We prove this by contradiction. Suppose that Property \ref{prop:second_att} does not hold. Then one of the following two cases happens: (1) $e_1 = \{x,y\}_{H,t}$ and $e_2 = \{x,y\}_{H,t}$ for some $H \in \{H_0,H^*\}$ and $t \in \{\tau, \sigma\}$. Then, $H$ contains multi-hyperarcs, which is a contradiction with it being in $\mathcal{H}^{\textnormal{stub}}_{s}(\vb*{d})$. (2) $e_1 = \{x,y\}_{H_0,t}$ and $e_2 = \{x,y\}_{H^*,t}$. This would mean that $(\{t\},\{x,y\})$ appears in both $H_0$ and $H^*$ and should thus not appear in the symmetric difference.\end{proof}

The map $M(H_0 \Delta H^*)$, contains trails with edges alternating between $e^1=H_0$ and $e^1=H^*$, which are called \textit{alternating trails}. Recall that in a trail, each edge is distinct. A closed alternating trail is called an alternating circuit. In such an alternating circuit, each edge is used at most once, while vertices can be used more often. We use the term \textit{alternating cycle} to refer to an alternating circuit, in line with \cite{berger2010}.
\\

\begin{definition}[Alternating trail, alternating cycle]
A trail $v_1,e_1,v_2,e_2,\hdots,e_{n-1},v_n$ in $M(H_0 \Delta H^*)$ is an alternating trail of length $n$ if $e_i^1 \neq e_{i+1}^1$ for all $i=1,\hdots,{n-2}$. An alternating trail of even length with $v_n=v_1$ is called an alternating cycle.
\end{definition}

Since every vertex in $M(H_0 \Delta H^*)$ is adjacent to the same number of edges with first attribute $H_0$ as with first attribute $H^*$, $M(H_0 \Delta H^*)$ can be decomposed into edge-disjoint alternating cycles.
\\
\begin{lemma}
\label{lemma:decomposition_alt_cycles}
    Let $\vb*{d} \in D$. For any $H_0,H^* \in \mathcal{H}^{\textnormal{stub}}_{s}(\vb*{d})$, $M(H_0 \Delta H^*)$ decomposes into edge-disjoint alternating cycles.
\end{lemma}
\begin{proof}
    Since $H_0$ and $H^*$ have the same degree sequence, every vertex in $M(H_0 \Delta H^*)$ is adjacent to an equal number of edges with first attribute $H_0$ as with $H^*$. The edge-disjoint alternating cycles can be constructed as follows: start with some edge $e_1=\{a_1,a_2\}_{H_0,.}$, then add $e_2=\{a_2,a_3\}_{H^*,.}$, followed by $e_3 = \{a_3,a_4\}_{H_0,.}$. This process of adding edges with alternating first attributes continues until we reach an edge $e_{2n}=\{x,a_1\}$, which closes a circuit. If unused edges remain, the procedure is repeated with one of them, until all edges are part of an edge-disjoint alternating cycle. The process cannot terminate prematurely, since that would imply the existence of a vertex with a different number of adjacent edges with first attribute $H_0$ and $H^*$.
\end{proof}

There may be multiple ways to decompose some $M(H_0 \Delta H^*)$ into alternating cycles. For example, $M(H_0 \Delta H^*)$ in Figure \ref{fig:mapped symmetric difference} can be decomposed into two alternating cycles $a,b,c,d,a$ and $c,e,c$, or into one alternating cycle $a,b,c,e,c,d,a$. We therefore introduce the concept of minimal alternating cycles.
\\
\begin{definition}[Alternating 2-cycle]
A cycle in $M(H_0,H^*)$ is an alternating 2-cycle if it is an alternating cycle of length 2.\\
\end{definition}

\begin{definition}[Minimal alternating cycle]
    An alternating cycle in $M(H_0 \Delta H^*)$ is a minimal alternating cycle if i) it is an alternating 2-cycle or ii) it contains no subcycle that is an alternating cycle and it contains no edges that are part of an alternating 2-cycle.
\end{definition}

In $M(H_0 \Delta H^*)$ in Figure \ref{fig:mapped symmetric difference}, there is one alternating 2-cycle $c,e,c$, and there are two minimal alternating cycles $a,b,c,d,a$ and $c,e,c$. The alternating cycle $a,b,c,e,c,d,a$ is not a minimal alternating cycle, since it contains the alternating cycle $a,b,c,d,a$ as a subcycle, and also since it contains the edge $c,e$ which is in an alternating 2-cycle.

We now strengthen Lemma \ref{lemma:decomposition_alt_cycles} by showing that $M(H_0 \Delta H^*)$ decomposes into edge-disjoint minimal alternating cycles.
\\
\begin{lemma}
\label{lemma:decomposition_min_alt_cycles}
    Let $\vb*{d} \in D$. For any $H_0,H^* \in \mathcal{H}^{\textnormal{stub}}_{s}(\vb*{d})$, $M(H_0 \Delta H^*)$ decomposes into edge-disjoint minimal alternating cycles.
\end{lemma}
\begin{proof}
    The proof is by construction. First, choose each alternating 2-cycle as one minimal alternating cycle. This is possible, since any two alternating 2-cycles are edge disjoint: by Lemma \ref{prop:second_att}, there are at most two edges between any two vertices in $M(H_0 \Delta H^*)$. The remaining edges are not part of any 2-cycle and can be decomposed into edge-disjoint alternating cycles, by similar arguments as in the proof of Lemma \ref{lemma:decomposition_alt_cycles}. If such an alternating cycle is not minimal, then it must be because it contains a subcycle that is an alternating cycle. All alternating cycles that are not minimal can be split into their alternating subcycles, until all alternating cycles are minimal.
\end{proof}

In the proof of Lemma \ref{lemma:stub_s_connected}, it will be useful to have a method to check if two directed hypergraphs with degrees in $D$ and without degenerate hyperarcs or multi-hyperarcs have the same degree sequence, by means of the presented mapping. For this, we introduce the following sets of edges in $M(H_0 \Delta H^*)$: 
\begin{itemize}
    \item $f_{H_0}(M(H_0 \Delta H^*)) := \{e \, | \, e \in M(H_0 \Delta H^*) \textnormal{ and } e^1 = H_0\}$;
    \item $f_{H_0,t}(M(H_0 \Delta H^*)) = \{e \, | \, e \in M(H_0 \Delta H^*) \textnormal{ and } e^1 = H_0 \textnormal{ and } e^2 = t\}$.\\
\end{itemize}

\begin{claim}
\label{claim:H_in_H_s_equivalences}
    Let $\Tilde{\vb*{d}},\hat{\vb*{d}} \in D$, let $\Tilde{H} \in \mathcal{H}^{\textnormal{stub}}_{s}(\Tilde{\vb*{d}})$ and $\hat{H} \in \mathcal{H}^{\textnormal{stub}}_{s}(\hat{\vb*{d}})$. Then,
    \begin{align*}
        \Tilde{\vb*{d}}=\hat{\vb*{d}} \iff &M(\Tilde{H} \Delta \hat{H}) \textnormal{ decomposes into edge-disjoint alternating cycles and } \\
        &\forall t \in \{\tau, \sigma\}: |f_{\Tilde{H},t}(M(\Tilde{H}\Delta \hat{H}))| = |f_{\hat{H},t}(M(\Tilde{H}\Delta \hat{H}))|.
    \end{align*} 
\end{claim}
\begin{proof}
    $\Rightarrow$: By Lemma \ref{lemma:decomposition_alt_cycles}, $M(\Tilde{H} \Delta \hat{H})$ decomposes into edge-disjoint alternating cycles. In addition, if $|f_{\Tilde{H},t}(M(\Tilde{H} \Delta \hat{H}))| > |f_{\hat{H},t}(M(\Tilde{H} \Delta \hat{H}))|$ for some tail $t \in \{\tau, \sigma\}$, then the out-degree of $t$ in $\Tilde{H}$ would be larger than in $\hat{H}$, contradicting with $\Tilde{H}$ and $\hat{H}$ having equal degree sequences:
    \begin{align*}
        d_t^{\textnormal{out}}(\Tilde{H}) &= |\{a \, | \, a \in A(\Tilde{H}) \textnormal{ and } a=(t,.)\}|\\
        &= |\{ a \, | \, a \in A(\Tilde{H}) \cap A(\hat{H}) \textnormal{ and } a=(t,.)\}| + |f_{\Tilde{H},t}(M(\Tilde{H} \Delta \hat{H}))|\\
        &> |\{ a \, | \, a \in A(\Tilde{H}) \cap A(\hat{H}) \textnormal{ and } a=(t,.)\}| + |f_{\hat{H},t}(M(\Tilde{H} \Delta \hat{H}))|\\
        &= |\{a \, | \, a \in A(\hat{H}) \textnormal{ and } a=(t,.)\}|\\
        &= d_t^{\textnormal{out}}(\hat{H}).
    \end{align*}
    $\Leftarrow$: If $M(\Tilde{H} \Delta \hat{H})$ decomposes into edge-disjoint alternating cycles, then every vertex in $M(\Tilde{H} \Delta \hat{H})$ is adjacent to the same number of edges with first attribute $\Tilde{H}$ as the number of edges with first attribute $\hat{H}$. Therefore, the in-degrees of the vertices in $\Tilde{H}$ equal the in-degrees of the vertices in $\hat{H}$. Since $\forall t \in \{\tau, \sigma\}: |f_{\Tilde{H},t}(M(\Tilde{H}\Delta \hat{H}))| = |f_{\hat{H},t}(M(\Tilde{H}\Delta \hat{H}))|$, the out-degrees of the vertices $\tau$ and $\sigma$ in $\Tilde{H}$ equal the out-degrees of those vertices in $\hat{H}$. This establishes that $\Tilde{H}$ and $\hat{H}$ have the same vertex degree sequence. Since $\Tilde{\vb*{d}},\hat{\vb*{d}} \in D$, $\Tilde{H}$ and $\hat{H}$ have the same hyperarc degree sequence. Therefore, $\Tilde{\vb*{d}}=\hat{\vb*{d}}$.
\end{proof}

Lastly, we relate the symmetric difference to the double hyperarc shuffle in the following lemma.
\\
\begin{lemma}
\label{lemma:sym_diff_hypershuffle}
Let $\vb*{d} \in D$ and let $H_0, H^* \in \mathcal{H}^{\textnormal{stub}}_{s}(\vb*{d})$ with $H_0 \neq H^*$. Then, 
\begin{align*}
    \exists a,b \in A(H_0): H^* \in S(a,b|H_0) \iff |H_0 \Delta H^*| = 4.
\end{align*} 
\end{lemma}
\begin{proof}
    $\Rightarrow$: if $H^* \in S(a,b|H_0)$, then $A(H^*) = (A(H_0) \backslash \{a,b\}) \cup \{\hat{a},\hat{b}\}$, for some shuffled hyperarcs $\hat{a},\hat{b}$. In particular, $H_0 \Delta H^* = \{a,b,\hat{a},\hat{b}\}$, in which all hyperarcs are distinct. Else, $\{a,b\} = \{\hat{a}, \hat{b}\}$, which contradicts the assumption that $H_0 \neq H^*$. \\
    $\Leftarrow$: if $|H_0 \Delta H^*|=4$, then $|H_0 \Delta H^*|=\{a_1,a_2,a_3,a_4\}$, where $M(a_1),M(a_3) \in f_{H_0}(M(H_0 \Delta H^*))$ and $M(a_2),M(a_4) \in f_{H^*}(M(H_0 \Delta H^*))$. Now, $\{M(a_1),M(a_2),M(a_3),M(a_4)\}$ decomposes into alternating cycles (by Lemma \ref{lemma:decomposition_alt_cycles}). These could be either one minimal alternating cycle of length 4, two disconnected alternating 2-cycles or two alternating 2-cycles sharing one vertex. In the first case, consider four distinct vertices $v,w,x,y \in V$ and two vertices $s,t \in \{\tau, \sigma\}$, that are not necessarily distinct. Then,
    \begin{align*}
        M(a_1) &= \{v,w\}_{H_0,t}\\
        M(a_2) &= \{w,x\}_{H^*,t}\\
        M(a_3) &= \{x,y\}_{H_0,s}\\
        M(a_4) &= \{v,y\}_{H^*,s}.
    \end{align*}
    Now, 
    \begin{align*}
        a_2 &= (t,(a_1^{\head} \backslash \{v\}) \cup \{x\})\\
        a_4 &= (s,(a_3^{\head} \backslash \{x\}) \cup \{v\}),
    \end{align*}
    and so $H^* \in S(a_1,a_3|H_0)$. In the second case, consider four distinct vertices $v,w,x,y \in V$ and two distinct vertices $s,t \in \{\tau, \sigma\}$. Then,
    \begin{align*}
        M(a_1) &= \{v,w\}_{H_0,t}\\
        M(a_2) &= \{v,w\}_{H^*,s}\\
        M(a_3) &= \{x,y\}_{H_0,s}\\
        M(a_4) &= \{x,y\}_{H^*,t}.
    \end{align*}
    Now,
    \begin{align*}
        a_2 &= ((a_1^{\tail} \backslash \{t\}) \cup \{s\}, a_1^{\head})\\
        a_4 &= ((a_3^{\tail} \backslash \{s\}) \cup \{t\}, a_3^{\head})
    \end{align*}
    and so $H^* \in S(a_1,a_3|H_0)$. The third case equals the second case but with $v=x$.
\end{proof}

Let $\vb*{d} \in D$ and let $H_0,H^* \in \mathcal{H}^{\textnormal{stub}}_{s}(\vb*{d})$. The proof of Lemma \ref{lemma:stub_s_connected} is by induction on the number of edges in the symmetric difference $M(H_0 \Delta H^*)$. The edges are counted as follows: each edge that is not in an alternating 2-cycle counts as 1, and each edge in an alternating 2-cycle counts as $1/2$. This edge count is denoted by $\mathcal{E}(M(H_0 \Delta H^*))$:
\begin{align}
\label{eq:counter}
    \mathcal{E}(M(H_0 \Delta H^*)) = \, &|\{e: e \in M(H_0 \Delta H^*) \textnormal{ and } e \notin \textnormal{alt. 2-cycle}\}|\\
    &+ \frac{1}{2}|\{e: e \in M(H_0 \Delta H^*) \textnormal{ and } e \in \textnormal{alt. 2-cycle}\}|. \nonumber
\end{align}
Since each alternating 2-cycle contains exactly 2 edges, $\mathcal{E}(M(H_0 \Delta H^*)) \in \mathds{N}$. The counting mechanism is illustrated in Figure \ref{fig:mapped symmetric difference}.

The intuition behind this edge count is as follows. In the proof of Lemma \ref{lemma:stub_s_connected}, we show how to either decrease the number of edges in the symmetric difference, $M(H_0 \Delta H^*)$, or to increase the number of alternating 2-cycles in the symmetric difference while keeping the total number of edges fixed. The edge counter $\mathcal{E}(M(H_0 \Delta H^*))$ is defined in such a way that it decreases in both cases, such that mathematical induction can be applied. It is not possible to exclude the edges in alternating 2-cycles in the counter, as we explain in the proof of Lemma \ref{lemma:connected_when_multiple_cycles}.

\underline{Base case}: $\mathcal{E}(M(H_0 \Delta H^*))=0$.\\
Then, $M(H_0 \Delta H^*)$ contains no edges. Therefore, $H_0 = H^*$ and Lemma \ref{lemma:stub_s_connected} holds.

\underline{Induction hypothesis}: There exists a $k \in \mathds{N}$ such that $\mathcal{G}(\mathcal{H}^{\textnormal{stub}}_{s}(\vb*{d}))$ contains a path from $H_0$ to $H^*$ if $\mathcal{E}(M(H_0 \Delta H^*)) \leq k$.

\underline{To prove}: $\mathcal{G}(\mathcal{H}^{\textnormal{stub}}_{s}(\vb*{d}))$ contains a path from $H_0$ to $H^*$ if $\mathcal{E}(M(H_0 \Delta H^*))=k+1$.

We aim to prove existence of a directed hypergraph $\Tilde{H} \in \mathcal{H}^{\textnormal{stub}}_{s}(\vb*{d})$, with $\mathcal{E}(M(H_0 \Delta \Tilde{H})) \leq k$ and $\mathcal{E}(M(\Tilde{H} \Delta H^*)) \leq k$. Then, by the induction hypothesis, $\mathcal{G}(\mathcal{H}^{\textnormal{stub}}_{s}(\vb*{d}))$ contains a path from $H_0$ to $\Tilde{H}$ and from $\Tilde{H}$ to $H^*$.

The proof is split into two cases, which are covered in Lemmas \ref{lemma:connected_when_multiple_cycles} and \ref{lemma:connected_when_one_cycle}:
\begin{itemize}
    \item $M(H_0 \Delta H^*)$ contains multiple minimal alternating cycles (Lemma \ref{lemma:connected_when_multiple_cycles});
    \item $M(H_0 \Delta H^*)$ contains exactly one minimal alternating cycle (Lemma \ref{lemma:connected_when_one_cycle}).\\
\end{itemize}

\begin{lemma}
\label{lemma:connected_when_multiple_cycles}
Let $\vb*{d} \in D$, let $H_0,H^* \in \mathcal{H}^{\textnormal{stub}}_{s}(\vb*{d})$ and $\mathcal{E}(M(H_0 \Delta H^*)) = k+1$. If $M(H_0 \Delta H^*)$ consists of multiple minimal alternating cycles, then $\mathcal{G}(\mathcal{H}^{\textnormal{stub}}_{s}(\vb*{d}))$ contains a path from $H_0$ to $H^*$.
\end{lemma}
\begin{proof}
    Let $M(H_0 \Delta H^*)$ decompose into multiple minimal edge-disjoint alternating cycles $C_1,C_2,\hdots,C_n$ (Lemma \ref{lemma:decomposition_min_alt_cycles}). Then for all $i,j \in [n]$, $i \neq j$:
\begin{align}
\label{eq:length cycle leq k}
    \mathcal{E}(C_i) \leq \mathcal{E}(M(H_0 \Delta H^*)) - \mathcal{E}(C_j) \leq (k+1)-1 = k.
\end{align}
Note that for this inequality to hold, it is important that the edge count does not neglect edges in 2-cycles, as in that case $ \mathcal{E}(C_j) $ could equal zero, and we would obtain the inequality $\mathcal{E}(C_i) \leq k+1$, on which we cannot apply the induction hypothesis in cases 1 and 2 that follow. 

We would like to construct the directed hypergraph $\Tilde{H}$ with $M(H_0 \Delta \Tilde{H}) = C_1$. Then, $\Tilde{H} \in \mathcal{H}^{\textnormal{stub}}_{s}(\vb*{d})$ if and only if $f_{H_0,t}(C_1)=f_{H^*,t}(C_1)$ for all $t \in \{\tau,\sigma\}$ (Claim \ref{claim:H_in_H_s_equivalences}). Therefore, we distinguish the following two cases, depending on whether the number of tails equal to $\tau$ in the cycle $C_1$ are equal in $H_0$ and $H^*$ or not, illustrated in Figure \ref{multiple_cycles_cases}.
\begin{enumerate}
    \item $|f_{H_0,\tau}(C_1)|=|f_{H^*,\tau}(C_1)|$;
    \item $|f_{H_0,\tau}(C_1)|\neq|f_{H^*,\tau}(C_1)|$. \label{case:case_that_splits_in_7}
\end{enumerate}

\begin{figure}[tbp]
\centering
\begin{subfigure}[t]{0.4\textwidth}
        \centering
        \begin{tikzpicture}[-,>=Stealth,auto,
                    thick,main node/.style={circle,draw, minimum size=\tikznodesize}, inner sep=1pt]

  \node[main node] (a) at (0,0) {};
  \node[main node] (b) at (2*\tikzscale,0) {};
  \node[main node] (c) at (3*\tikzscale,-1.7*\tikzscale) {};
  \node[main node] (d) at (2*\tikzscale,-3.5*\tikzscale) {};
  \node[main node] (e) at (0,-3.5*\tikzscale) {};
  \node[main node] (f) at (-1*\tikzscale,-1.7*\tikzscale) {};

  \node at (1*\tikzscale,-1.7*\tikzscale) {$C_1$};

  \path[every node/.style={font=\sffamily\small}]
    (a) edge[cyan] node {} (b)
    (b) edge[purple, line width=1.3pt, dashed] node {} (c)
    (c) edge[cyan] node {} (d)
    (d) edge[purple, line width=1.3pt] node {} (e)
    (e) edge[cyan, dashed] node {} (f)
    (f) edge[purple, line width=1.3pt] node {} (a);

\end{tikzpicture}
    \caption{$M(H_0 \Delta H^*)$, where $|f_{H_0,\tau}(C_1)|=2=|f_{H^*,\tau}(C_1)|$.}
    \label{fig:multiple_cycles_case1}
\end{subfigure}
\begin{subfigure}[t]{0.4\textwidth}
        \centering
    \begin{tikzpicture}[-,>=Stealth,auto,
                    thick,main node/.style={circle,draw, minimum size=\tikznodesize}, inner sep=1pt]

  \node[main node] (a) at (0,0) {};
  \node[main node] (b) at (2*\tikzscale,0) {};
  \node[main node] (c) at (3*\tikzscale,-1.7*\tikzscale) {};
  \node[main node] (d) at (2*\tikzscale,-3.5*\tikzscale) {};
  \node[main node] (e) at (0,-3.5*\tikzscale) {};
  \node[main node] (f) at (-1*\tikzscale,-1.7*\tikzscale) {};

  \node at (1*\tikzscale,-1.7*\tikzscale) {$C_1$};

  \path[every node/.style={font=\sffamily\small}]
    (a) edge[cyan] node {} (b)
    (b) edge[purple, line width=1.3pt, dashed] node {} (c)
    (c) edge[cyan] node {} (d)
    (d) edge[purple, line width=1.3pt, dashed] node {} (e)
    (e) edge[cyan] node {} (f)
    (f) edge[purple, line width=1.3pt] node {} (a);

\end{tikzpicture}
    \caption{$M(H_0 \Delta H^*)$, where $|f_{H_0,\tau}(C_1)|=3\neq1=|f_{H^*,\tau}(C_1)|$. Here $\alpha=2$.}
    \label{fig:multiple_cycles_case2}
\end{subfigure}
\caption{Two cases for a cycle $C_1$: $|f_{H_0,\tau}(C_1)|=|f_{H^*,\tau}(C_1)|$ or $|f_{H_0,\tau}(C_1)| \neq |f_{H^*,\tau}(C_1)|$. Legend as in Figure \ref{fig:mapped symmetric difference}.}
\label{multiple_cycles_cases}
\end{figure}
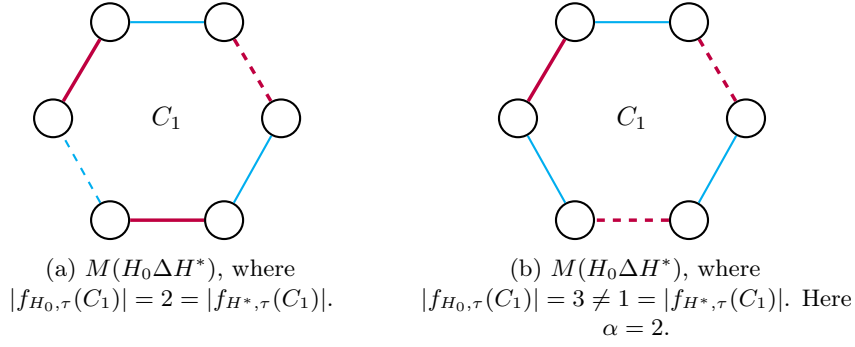

\underline{Case 1}: $|f_{H_0,\tau}(C_1)|=|f_{H^*,\tau}(C_1)|$ (Figure \ref{fig:multiple_cycles_case1}). \\
We construct the directed hypergraph $\Tilde{H}$ with 
\begin{align*}
    A(\Tilde{H}) &= (A(H_0) \backslash \{M^{-1}(e): e \in C_1 \textnormal{ and } e^1=H_0\}) \cup \{M^{-1}(e): e \in C_1 \textnormal{ and } e^1 = H^*\}.
\end{align*}
By Claim \ref{claim:H_in_H_s_equivalences}, $\Tilde{H} \in \mathcal{H}^{\textnormal{stub}}_{s}(\vb*{d})$, as $\Tilde{H}$ contains no degenerate or multi-hyperarcs by construction, $M(H_0 \Delta \Tilde{H})=C_1$ is an alternating cycle and 
\begin{align*}
|f_{\Tilde{H},\tau}(M(H_0 \Delta \Tilde{H})| = |f_{\Tilde{H},\tau}(C_1)| = |f_{H_0,\tau}(C_1)| = |f_{H_0,\tau}(M(H_0 \Delta \Tilde{H})|,
\end{align*}
where the second equality is by assumption of case 1. Moreover,
\begin{align*}
    \mathcal{E}(M(H_0 \Delta \Tilde{H})) = \mathcal{E}(C_1) \leq k
\end{align*}
by~\eqref{eq:length cycle leq k}, and
\begin{align*}
    \mathcal{E}(M(\Tilde{H} \Delta H^*))= \mathcal{E}(M(H_0 \Delta H^*)) - \mathcal{E}(C_1) \leq k+1 - 2 = k-1.
\end{align*}

By the induction hypothesis, $\mathcal{G}(\mathcal{H}^{\textnormal{stub}}_{s}(\vb*{d}))$ contains a path from $H_0$ to $\Tilde{H}$ and from $\Tilde{H}$ to $H^*$.

\underline{Case 2}: $|f_{H_0,\tau}(C_1)|\neq|f_{H^*,\tau}(C_1)|$ (Figure \ref{fig:multiple_cycles_case2}).\\
W.l.o.g., let $|f_{H_0,\tau}(C_1)|>|f_{H^*,\tau}(C_1)|$. Since $|f_{H_0,\tau}(M(H_0 \Delta H^*))|=|f_{H^*,\tau}(M(H_0 \Delta H^*))|$ (Claim \ref{claim:H_in_H_s_equivalences}), 
another cycle $C_2$ must exist with $|f_{H_0,\tau}(C_2)|<|f_{H^*,\tau}(C_2)|$ and thus $|f_{H_0,\sigma}(C_2)|>|f_{H^*,\sigma}(C_2)|$. This is illustrated in Figure \ref{fig:multiple_cycles_case2_preswap}. We now swap tails of edges in $f_{H_0,\tau}(C_1)$ with tails of edges in $f_{H_0,\sigma}(C_2)$ until $|f_{H_0,\tau}(C_1)|=|f_{H^*,\tau}(C_1)|$ and then construct $\Tilde{H}$ in a similar way as in case 1.

\begin{figure}[tbp]
\centering
\begin{subfigure}{\textwidth}
        \centering
        \begin{tikzpicture}[-,>=Stealth,auto,
                    thick,main node/.style={circle,draw, minimum size=\tikznodesize}, inner sep=1pt]

  \node[main node] (a) at (0,0) {};
  \node[main node] (b) at (2*\tikzscale,0) {};
  \node[main node] (c) at (3*\tikzscale,-1.7*\tikzscale) {};
  \node[main node] (d) at (2*\tikzscale,-3.5*\tikzscale) {};
  \node[main node] (e) at (0,-3.5*\tikzscale) {};
  \node[main node] (f) at (-1*\tikzscale,-1.7*\tikzscale) {};

  \node[font=\bfseries\sffamily] at (1*\tikzscale,-1.7*\tikzscale) {$C_1$};

  \node[main node] (g) at (6,1.1*\tikzscale) {};
  \node[main node] (h) at (2*\tikzscale +6,1.1*\tikzscale) {};
  \node[main node] (i) at (3.5*\tikzscale +6,-0.2*\tikzscale) {};
  \node[main node] (j) at (3.5*\tikzscale+6,-2.2*\tikzscale) {};
  \node[main node] (k) at (2*\tikzscale+6,-3.5*\tikzscale) {};
  \node[main node] (l) at (0+6,-3.5*\tikzscale) {};
  \node[main node] (m) at (-1.5*\tikzscale +6,-2.2*\tikzscale) {};
  \node[main node] (n) at (-1.5*\tikzscale +6,-0.2*\tikzscale) {};

  \node[font=\bfseries\sffamily] at (1*\tikzscale+6,-1.2*\tikzscale) {$C_2$};

  \path[every node/.style={font=\sffamily\small}]
    (a) edge[cyan] node {$\color{black} e_1$} (b)
    (b) edge[purple, line width=1.3pt, dashed] node {} (c)
    (c) edge[cyan] node {} (d)
    (d) edge[purple, line width=1.3pt, dashed] node {} (e)
    (e) edge[cyan] node {} (f)
    (f) edge[purple, line width=1.3pt] node {} (a)

    (g) edge[cyan, dashed] node {} (h)
    (h) edge[purple, line width=1.3pt] node {} (i)
    (i) edge[cyan, dashed] node {$\color{black} e_2$} (j)
    (j) edge[purple, line width=1.3pt] node {} (k)
    (k) edge[cyan, dashed] node {} (l)
    (l) edge[purple, line width=1.3pt] node {} (m)
    (m) edge[cyan] node {} (n)
    (n) edge[purple, line width=1.3pt] node {} (g);

\end{tikzpicture}
    \caption{$M(H_0 \Delta H^*)$.}
    \label{fig:multiple_cycles_case2_preswap}
\end{subfigure}\\
\vspace{0.5cm}
\begin{subfigure}{\textwidth}
        \centering
    \begin{tikzpicture}[-,>=Stealth,auto,
                    thick,main node/.style={circle,draw, minimum size=\tikznodesize}, inner sep=1pt]

  \node[main node] (a) at (0,0) {};
  \node[main node] (b) at (2*\tikzscale,0) {};
  \node[main node] (c) at (3*\tikzscale,-1.7*\tikzscale) {};
  \node[main node] (d) at (2*\tikzscale,-3.5*\tikzscale) {};
  \node[main node] (e) at (0,-3.5*\tikzscale) {};
  \node[main node] (f) at (-1*\tikzscale,-1.7*\tikzscale) {};

  \node[font=\bfseries\sffamily] at (1*\tikzscale,-1.7*\tikzscale) {$C_1$};

  \node[main node] (g) at (6,1.1*\tikzscale) {};
  \node[main node] (h) at (2*\tikzscale +6,1.1*\tikzscale) {};
  \node[main node] (i) at (3.5*\tikzscale +6,-0.2*\tikzscale) {};
  \node[main node] (j) at (3.5*\tikzscale+6,-2.2*\tikzscale) {};
  \node[main node] (k) at (2*\tikzscale+6,-3.5*\tikzscale) {};
  \node[main node] (l) at (0+6,-3.5*\tikzscale) {};
  \node[main node] (m) at (-1.5*\tikzscale +6,-2.2*\tikzscale) {};
  \node[main node] (n) at (-1.5*\tikzscale +6,-0.2*\tikzscale) {};

  \node[font=\bfseries\sffamily] at (1*\tikzscale+6,-1.2*\tikzscale) {$C_2$};

  \path[every node/.style={font=\sffamily\small}]
    (a) edge[cyan, dashed] node {$\color{black} e_1$} (b)
    (b) edge[purple, line width=1.3pt, dashed] node {} (c)
    (c) edge[cyan] node {} (d)
    (d) edge[purple, line width=1.3pt, dashed] node {} (e)
    (e) edge[cyan] node {} (f)
    (f) edge[purple, line width=1.3pt] node {} (a)

    (g) edge[cyan, dashed] node {} (h)
    (h) edge[purple, line width=1.3pt] node {} (i)
    (i) edge[cyan] node {$\color{black} e_2$} (j)
    (j) edge[purple, line width=1.3pt] node {} (k)
    (k) edge[cyan, dashed] node {} (l)
    (l) edge[purple, line width=1.3pt] node {} (m)
    (m) edge[cyan] node {} (n)
    (n) edge[purple, line width=1.3pt] node {} (g);

\end{tikzpicture}
    \caption{$M(H_1 \Delta H^*)$. Performing the tail swap (\ref{eq:tailswap}) with the edges $e_1$ and $e_2$ as in Figure \ref{fig:multiple_cycles_case2_preswap} results in the directed hypergraph $H_1$.}
    \label{fig:multiple_cycles_case2_afterswap}
\end{subfigure}
\caption{If $|f_{H_0,\tau}(C_1)| > |f_{H^*,\tau}(C_1)|$ then there must exist an alternating cycle $C_2$ with $|f_{H_0,\tau}(C_2)| < |f_{H^*,\tau}(C_2)|$. Legend as in Figure \ref{fig:mapped symmetric difference}.}
\end{figure}
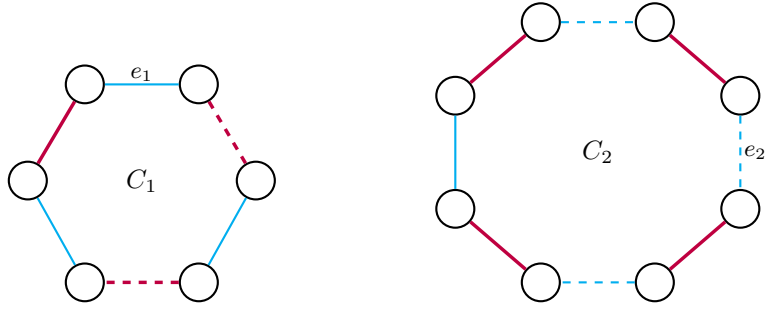
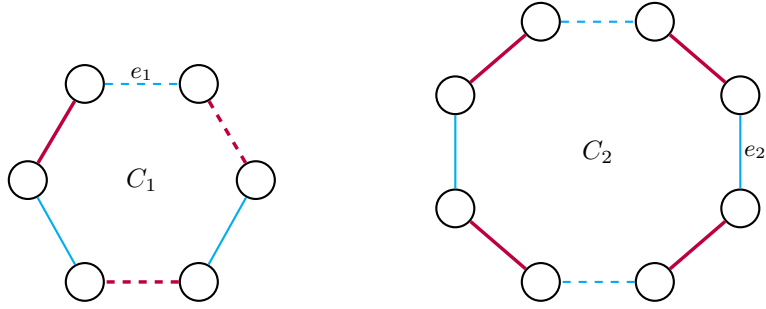

To be more precise, let 
\begin{align*}
    |f_{H_0,\tau}(C_1)|=|f_{H^*,\tau}(C_1)| + \alpha,
\end{align*}
where $\alpha \leq \frac{1}{2}|C_1|$.

Pick any edges $\{a,b\}_{H_0,\tau} \in C_1, \{c,d\}_{H_0,\sigma} \in C_2$. If $(\sigma, \{a,b\}) \notin A(H_0) \textnormal{ and } (\tau,\{c,d\}) \notin A(H_0)$, then perform the double hyperarc shuffle
\begin{align}\label{eq:tailswap}
    (\tau,\{a,b\}) \times (\sigma,\{c,d\}) &\rightarrow (\sigma,\{a,b\}),(\tau,\{c,d\})
\end{align}
in $H_0$. This creates the directed hypergraph $H_1$ with 
$H_0 \Delta H_1 = \{(\tau,\{a,b\}) , (\sigma,\{c,d\}), (\sigma,\{a,b\}),(\tau,\{c,d\})\}$. Moreover, $H_1 \Delta H^* = ((H_0 \Delta H^*) \backslash \{(\tau,\{a,b\}),(\sigma,\{c,d\})\}) \cup \{(\sigma,\{a,b\}),(\tau,\{c,d\})\}$ and thus
\begin{align*}
    \mathcal{E}(M(H_1 \Delta H^*)) &= \mathcal{E}(M(H_0 \Delta H^*)).
\end{align*}

In particular, we obtain 
\begin{align*}
    |f_{H_1,\tau}(C_1)|=|f_{H^*,\tau}(C_1)| + \alpha - 1.
\end{align*}

This tail swap is illustrated in Figure \ref{fig:multiple_cycles_case2_afterswap}. If one can perform $\alpha$ such swaps on $C_1$, using the cycle $C_2$ or possibly other cycles $C_2, C_3,\hdots,C_m$, $m \leq \alpha$, then this results in the directed hypergraph $H_{\alpha} \in \mathcal{H}^{\textnormal{stub}}_{s}(\vb*{d})$, with 
\begin{align*}
    \mathcal{E}(M(H_0 \Delta H_{\alpha})) &= 2 \alpha \leq \mathcal{E}(C_1) \leq k
\end{align*}
by~\eqref{eq:length cycle leq k}, and 
\begin{align*}
    \mathcal{E}(M(H_{\alpha} \Delta H^*)) &= \mathcal{E}(M(H_0 \Delta H^*)).
\end{align*}
If $C_1$ is an alternating 2-cycle, or any of the cycles $C_2,\hdots,C_m$ is an alternating 2-cycle, then this alternating 2-cycle disappears from $M(H_{\alpha} \Delta H^*)$ since both edges in the alternating 2-cycle now have the same tail (Property \ref{prop:second_att}) and we obtain $\mathcal{E}(M(H_{\alpha} \Delta H^*)) \leq \mathcal{E}(M(H_0 \Delta H^*)) - \frac{1}{2} \cdot 2$, so
\begin{align*}
    \mathcal{E}(M(H_{\alpha} \Delta H^*)) \leq \mathcal{E}(M(H_0 \Delta H^*)) - \frac{1}{2}\cdot 2 = k+1 - 1 = k.
\end{align*}
Now, let $\Tilde{H}= H_{\alpha}$. By construction, $\mathcal{G}(\mathcal{H}^{\textnormal{stub}}_{s}(\vb*{d}))$ contains a path from $H_0$ to $H_{\alpha}$. By the induction hypothesis, $\mathcal{G}(\mathcal{H}^{\textnormal{stub}}_{s}(\vb*{d}))$ contains a path from $\Tilde{H}$ to $H^*$.

If none of the cycles $C_1,C_2,\hdots,C_m$ are an alternating 2-cycle, then as in case 1, we can construct $\Tilde{H}$ using the cycle $C_1$:
\begin{align*}
    A(\Tilde{H}) &= (A(H_{\alpha}) \backslash \{M^{-1}(e): e \in C_1 \textnormal{ and } e^1=H_{\alpha}\}) \cup \{M^{-1}(e): e \in C_1 \textnormal{ and } e^1 = H^*\}.
\end{align*}

With similar arguments as in case 1 of the proof of Lemma \ref{lemma:connected_when_multiple_cycles}, we can show that $\Tilde{H} \in \mathcal{H}^{\textnormal{stub}}_{s}(\vb*{d})$, $\mathcal{E}(M(H_{\alpha} \Delta \Tilde{H})) \leq k$ and $\mathcal{E}(M(\Tilde{H} \Delta H^*)) \leq k-1$. By construction, $\mathcal{G}(\mathcal{H}^{\textnormal{stub}}_{s}(\vb*{d}))$ contains a path from $H_0$ to $H_{\alpha}$ and by the induction hypothesis, $\mathcal{G}(\mathcal{H}^{\textnormal{stub}}_{s}(\vb*{d}))$ contains a path from $H_{\alpha}$ to $H^*$.

If it is impossible to perform the tail swaps~\eqref{eq:tailswap} $\alpha$ times, then it must be that after performing $\alpha'<\alpha$ tail swaps, for all cycles $C_i \in \{C_2,\hdots,C_n\}$ with $|f_{H_0,\sigma}(C_i)|>|f_{H^*,\sigma}(C_i)|$, the suggested tail swap between an edge in $C_1$ and an edge in $C_i$ is not possible because they create multi-hyperarcs. Let $H_{\alpha'} \in \mathcal{H}^{\textnormal{stub}}_{s}(\vb*{d})$ be the directed hypergraph after the $\alpha'$ previous swaps. Then, all possible tail swaps must create multi-hyperarcs:
\begin{align*}
    &\forall \{a,b\}_{H_{\alpha'},\tau} \in C_1 \, \forall \{c,d\}_{H_{\alpha'},\sigma} \in C_i: (\sigma, \{a,b\}) \in A(H_{\alpha'}) \lor (\tau,\{c,d\}) \in A(H_{\alpha'}).
\end{align*}
In particular, one of the following two cases must hold:
\begin{enumerate}
    \item[(i)] $\forall \{a,b\}_{H_{\alpha'},\tau} \in C_1: (\sigma,\{a,b\}) \in A(H_{\alpha'})$; 
    \item[(ii)] $\forall \{c,d\}_{H_{\alpha'},\sigma} \in C_i: (\tau,\{c,d\}) \in A(H_{\alpha'})$.
\end{enumerate}

By symmetry, the two cases are considered identical. Therefore, we only discuss case (i), illustrated in Figure \ref{fig:multiple_cycles_case2_block_vanuit_C1}. Observe that:
\begin{itemize}
    \item By construction, $|f_{H_{\alpha'},\tau}(C_1)|=|f_{H^*,\tau}(C_1)| + \alpha-\alpha'$;
    \item $|C_1| \geq 4$. If instead $|C_1|=2$, then $C_1=(e_1,e_2)$ with $e_1^1 \neq e_2^1$ and $e_1^2 \neq e_2^2$ (Property \ref{prop:second_att})  
    W.l.o.g., let $e_1=\{a,b\}_{H_{\alpha'},\tau}$ and $e_2=\{a,b\}_{H^*,\sigma}$. By definition of the symmetric difference, $(\sigma, \{a,b\}) \notin A(H_{\alpha'})$. 
\end{itemize}

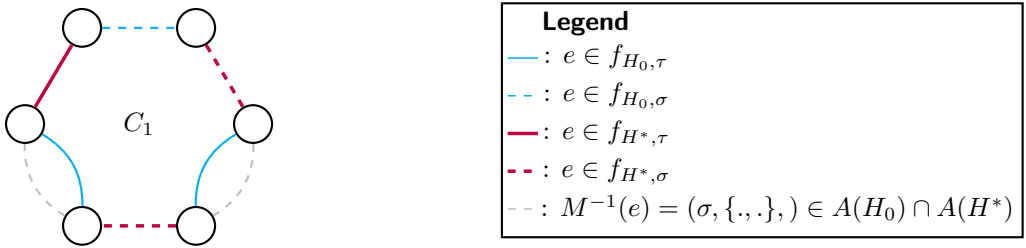
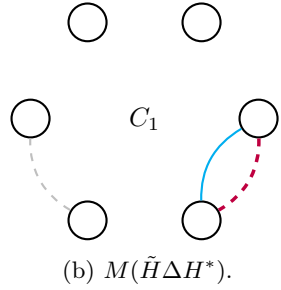
\begin{figure}[tbp]
\centering
\begin{subfigure}{\textwidth}
    \newcommand{\offsetx}{1}
\newcommand{\offsety}{-0.7}

\begin{tikzpicture}[-,>=Stealth,auto,
                    thick,main node/.style={circle,draw, minimum size=\tikznodesize}, inner sep=1pt]

  \node[main node] (a) at (0,0) {};
  \node[main node] (b) at (2*\tikzscale,0) {};
  \node[main node] (c) at (3*\tikzscale,-1.7*\tikzscale) {};
  \node[main node] (d) at (2*\tikzscale,-3.5*\tikzscale) {};
  \node[main node] (e) at (0,-3.5*\tikzscale) {};
  \node[main node] (f) at (-1*\tikzscale,-1.7*\tikzscale) {};

  \node[font=\bfseries\sffamily] at (1*\tikzscale,-1.7*\tikzscale) {$C_1$};

  \path[every node/.style={font=\sffamily\small}]
    (a) edge[cyan, dashed, line width = 3 pt] node {} (b)
    (b) edge[purple, dashed, line width = 3 pt] node {} (c)
    (c) edge[cyan, bend right] node {} (d)
    (c) edge[lightgray, dashed, bend left, line width = 3 pt] node {} (d)
    (d) edge[purple, dashed, line width = 3 pt] node {} (e)
    (e) edge[cyan, bend right, line width = 3 pt] node {} (f)
    (e) edge[lightgray, dashed, bend left] node {} (f)
    (f) edge[purple, line width = 3 pt] node {} (a);

% Legend title
  \node[anchor=south west, font=\bfseries\sffamily] at (5.0 + \offsetx,0.55+\offsety) {Legend};

  % Legend items
  \node at (5.9 + \offsetx,0.3 + \offsety) {: \textcolor{black}{$e \in f_{H_0,\tau}$}};
  \draw[cyan, line width=0.02cm] (4.6 + \offsetx,0.3 + \offsety) -- (5 + \offsetx,0.3 + \offsety);
  \node at (5.9 + \offsetx,-0.2 + \offsety) {: \textcolor{black}{$e \in f_{H_0,\sigma}$}};
  \draw[cyan,dashed, line width=0.02cm] (4.6 + \offsetx,-0.2 + \offsety) -- (5 + \offsetx,-0.2 + \offsety);
  \node at (5.9 + \offsetx,-0.7 + \offsety) {: \textcolor{black}{$e \in f_{H^*,\tau}$}};
  \draw[purple, line width=1.3pt] (4.6 + \offsetx,-0.7 + \offsety) -- (5 + \offsetx,-0.7 + \offsety);
  \node at (5.9 + \offsetx,-1.2 + \offsety) {: \textcolor{black}{$e \in f_{H^*,\sigma}$}};
  \draw[purple, line width=1.3pt, dashed] (4.6 + \offsetx,-1.2 + \offsety) -- (5 + \offsetx,-1.2 + \offsety);
  \node at (8.15+ \offsetx,-1.7 + \offsety) {: \textcolor{black}{$M^{-1}(e)=(\sigma,\{.,.\}) \in A(H_0) \cap A(H^*)$}};
  \draw[lightgray, dashed, line width=0.02cm] (4.6 + \offsetx,-1.7 + \offsety) -- (5 + \offsetx,-1.7 + \offsety);

  % Legend rectangle
  \node[draw, rectangle, fit={(4.7 + \offsetx,0.85 + \offsety) (11.2  + \offsetx,-1.9 + \offsety)}, inner sep=5pt] {};
  
\end{tikzpicture}
    \caption{$M(H_{\alpha'} \Delta H^*)$. For all edges in $f_{H_0,\tau}$ the directed hyperarc preimage of this edge, but with tail $\sigma$, is present in both $A(H_0)$ and $A(H^*)$. We obtain $|S|=1$. Let $S$ contain the bottom right cyan edge, then the artificial cycle $\Tilde{C}_1$ is indicated by the bold edges.}
    \label{fig:multiple_cycles_case2_block_vanuit_C1}
\end{subfigure}\\
\vspace{0.5cm}
\begin{subfigure}[t]{0.4\textwidth}
        \centering
    \newcommand{\offsetx}{1}
\newcommand{\offsety}{-0.7}

\begin{tikzpicture}[-,>=Stealth,auto,node distance=2cm,
                   [-,>=Stealth,auto,
                    thick,main node/.style={circle,draw, minimum size=\tikznodesize}, inner sep=1pt]

  \node[main node] (a) at (0,0) {};
  \node[main node] (b) at (2*\tikzscale,0) {};
  \node[main node] (c) at (3*\tikzscale,-1.7*\tikzscale) {};
  \node[main node] (d) at (2*\tikzscale,-3.5*\tikzscale) {};
  \node[main node] (e) at (0,-3.5*\tikzscale) {};
  \node[main node] (f) at (-1*\tikzscale,-1.7*\tikzscale) {};

  \node[font=\bfseries\sffamily] at (1*\tikzscale,-1.7*\tikzscale) {$C_1$};

  \path[every node/.style={font=\sffamily\small}]
    (c) edge[cyan, bend right] node {} (d)
    (c) edge[purple, line width=1.3pt, dashed, bend left] node {} (d)
    (e) edge[lightgray, dashed, bend left] node {} (f);
\end{tikzpicture}
    \caption{$M(\Tilde{H} \Delta H^*)$.}
    \label{fig:multiple_cycles_case2_finalswap}
\end{subfigure}
\caption{If no more tail swaps (\ref{eq:tailswap}) can be performed, then $\Tilde{H}$ can be constructed using the artificial cycle $\Tilde{C_1}$.}
\label{fig:construction_G_tilde_from_alpha'}
\end{figure}

We now create an artificial cycle $\Tilde{C_1}$, which may not be present in $M(H_{\alpha '} \Delta H^*)$, and construct $\Tilde{H}$ using this cycle, similarly as in case 1.

Let $S$ be a subset of $f_{H_{\alpha'},\tau}(C_1)$ such that $|S|=\alpha-\alpha'$. Let $\Tilde{S}$ be the set of edges in $S$, but with tails $\sigma$ instead of $\tau$:
\begin{align*}
    \Tilde{S} = \{\{e,f\}_{H_{\alpha'}, \sigma}: \{e,f\}_{H_{\alpha '}, \tau} \in S\}.
\end{align*}
We create the artificial cycle $\Tilde{C_1}$ which equals $C_1$, but with the edges in $S$ replaced by the edges in $\Tilde{S}$:
\begin{align*}
    \Tilde{C_{1}} = (C_1 \backslash S) \cup \Tilde{S}.
\end{align*}
$\Tilde{C}_1$ is illustrated in Figure \ref{fig:multiple_cycles_case2_block_vanuit_C1}. For $\Tilde{C_1}$ holds $|\{e \in \Tilde{C_1}: e^1=H_{\alpha'} \textnormal{ and } e^2 = \tau\}| = |\{e \in \Tilde{C_1}: e^1=H^* \textnormal{ and } e^2 = \tau\}|$. We construct $\Tilde{H}$ using $\Tilde{C_1}$:
\begin{align*}
    A(\Tilde{H}) &= (A(H_{\alpha'}) \backslash \{M^{-1}(e): e \in \Tilde{C_1} \textnormal{ and } e^1=H_{\alpha'}\}) \cup \{M^{-1}(e): e \in \Tilde{C_1} \textnormal{ and } e^1 = H^*\},
\end{align*}
which is illustrated in Figure \ref{fig:multiple_cycles_case2_finalswap}. With similar arguments as in case 1 of the proof of Lemma \ref{lemma:connected_when_multiple_cycles}, we can show that $\Tilde{H} \in \mathcal{H}^{\textnormal{stub}}_{s}(\vb*{d})$ and $\mathcal{E}(M(H_{\alpha'} \Delta \Tilde{H})) \leq k$. In addition,
\begin{align*}
    \mathcal{E}(M(\Tilde{H} \Delta H^*))&=  \mathcal{E}(M(H_{\alpha'} \Delta H^*)) - \frac{ 1}{2}|\Tilde{C_1}| - \Big(\frac{1}{2}|\Tilde{C_1}|-|S|\Big)  \\
    &= \mathcal{E}(M(H_{\alpha'} \Delta H^*)) - |\Tilde{C_1}| + |S|  \\
    &\leq \mathcal{E}(M(H_{\alpha'} \Delta H^*)) - \frac{1}{2} |\Tilde{C_1}|\\
    &\leq k+1 - 2 =k-1.
\end{align*}

The first equality is explained best by considering the example in Figure \ref{fig:construction_G_tilde_from_alpha'}, where $S$ contains only the bottom right cyan edge. The term $- \frac{ 1}{2}|\Tilde{C_1}|$ represents the edges in $f_{H^*}(\Tilde{C_1})$ in $M(H_{\alpha'} \Delta H^*)$ that are displayed in the example in purple. As shown in the example, all these edges have disappeared in $M(\Tilde{H} \Delta H^*)$. This is because these $\frac{1}{2}|\Tilde{C_1}|$ hyperarcs are not elements of $A(H_{\alpha' })$ and are elements of $A(\Tilde{H})$ and $A(H^*)$. The term $- (\frac{1}{2}|\Tilde{C_1}|-|S|)$ represents the edges in $f_{H_0}(\Tilde{C_1} \backslash S)$ in $M(H_{\alpha'} \Delta H^*)$. In the example, these are the cyan edges, with the bottom right cyan edge ($\in S$) removed. As shown in the example, all these edges have disappeared in $M(\Tilde{H} \Delta H^*)$. This is because these $\frac{1}{2}|\Tilde{C_1}|-|S|$ hyperarcs are elements of $A(H_{\alpha'})$ and not of $A(\Tilde{H})$ and $A(H^*)$. Lastly, we consider the edges in $\Tilde{S}$. As seen in the example, these edges are not present in $M(H_{\alpha'} \Delta H^*)$ and become part of an alternating 2-cycle in $M(\Tilde{H} \Delta H^*)$. This is because these edges are elements of $A(H_{\alpha'})$ and $A(H^*)$, and not of $A(\Tilde{H})$. In $\mathcal{E}(M(H_{\alpha'} \Delta H^*))$, the edges in $S$ counted as 1 edge, and in $\mathcal{E}(M(\Tilde{H} \Delta H^*))$, the edges in $S \cup \Tilde{S}$ count as $\frac{1}{2}\cdot 2=1$ edge. Therefore, for these edges, no term appears.

By construction, $\mathcal{G}(\mathcal{H}^{\textnormal{stub}}_{s}(\vb*{d}))$ contains a path from $H_0$ to $H_{\alpha'}$. By the induction hypothesis, $\mathcal{G}(\mathcal{H}^{\textnormal{stub}}_{s}(\vb*{d}))$ contains a path from $H_{\alpha'}$ to $\Tilde{H}$ and from $\Tilde{H}$ to $H^*$.
\end{proof}

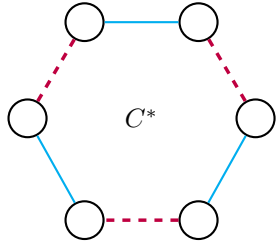
\begin{figure}[tbp]
    \centering
    \begin{tikzpicture}[-,>=Stealth,auto,
                    thick,main node/.style={circle,draw, minimum size=\tikznodesize}, inner sep=1pt]

  \node[main node] (a) at (0,0) {};
  \node[main node] (b) at (2*\tikzscale,0) {};
  \node[main node] (c) at (3*\tikzscale,-1.7*\tikzscale) {};
  \node[main node] (d) at (2*\tikzscale,-3.5*\tikzscale) {};
  \node[main node] (e) at (0,-3.5*\tikzscale) {};
  \node[main node] (f) at (-1*\tikzscale,-1.7*\tikzscale) {};

  \node[font=\bfseries\sffamily] at (1*\tikzscale,-1.7*\tikzscale) {$C^*$};

  \path[every node/.style={font=\sffamily\small}]
    (a) edge[cyan] node {} (b)
    (b) edge[purple, line width=1.3pt, dashed] node {} (c)
    (c) edge[cyan] node {} (d)
    (d) edge[purple, line width=1.3pt, dashed] node {} (e)
    (e) edge[cyan] node {} (f)
    (f) edge[purple, line width=1.3pt, dashed] node {} (a);

\end{tikzpicture}
    \caption{The minimal alternating cycle $C^*$ if for every subtrail $(e_1,e_2)$ of $C^*$ holds $e_1^2 \neq e_2^2$. Legend as in Figure \ref{fig:multiple_cycles_case2_block_vanuit_C1}.}
    \label{fig:single_cycle_no_neighbors_with_same_tail}
\end{figure}

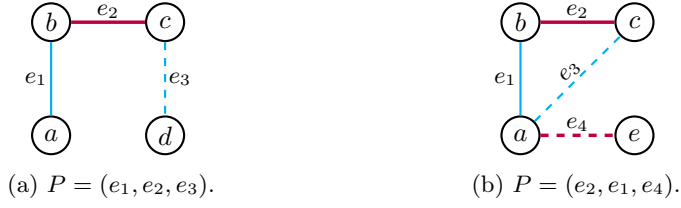
\begin{figure}[tbp]
\centering
\begin{subfigure}{0.4\textwidth}
\centering
    \begin{tikzpicture}[-,>=Stealth,auto,
                    thick,main node/.style={circle,draw, minimum size=\tikznodesize}, inner sep=1pt]

  \node[main node] (a) at (0,0) {$a$};
  \node[main node] (b) at (0,2*\tikzscale) {$b$};
  \node[main node] (c) at (2*\tikzscale,2*\tikzscale) {$c$};
  \node[main node] (d) at (2*\tikzscale,0) {$d$};

  \path[every node/.style={font=\sffamily\small}]
    (a) edge[cyan] node {$\color{black} e_1$} (b)
    (b) edge[purple, line width=1.3pt] node {$\color{black} e_2$} (c)
    (c) edge[cyan, dashed] node {$\color{black} e_3$} (d);

\end{tikzpicture}
    \caption{$P=(e_1,e_2,e_3)$.}
    \label{fig:easy_path_P}
\end{subfigure}
\begin{subfigure}{0.4\textwidth}
\centering
    \begin{tikzpicture}[-,>=Stealth,auto,
                    thick,main node/.style={circle,draw, minimum size=\tikznodesize}, inner sep=1pt]

  \node[main node] (a) at (0,0) {$a$};
  \node[main node] (b) at (0,2*\tikzscale) {$b$};
  \node[main node] (c) at (2*\tikzscale,2*\tikzscale) {$c$};
  \node[main node] (e) at (2*\tikzscale,0) {$e$};

  \path[every node/.style={font=\sffamily\small}]
    (a) edge[cyan] node {$\color{black} e_1$} (b)
    (b) edge[purple, line width=1.3pt] node {$\color{black} e_2$} (c)
    (a) edge[cyan, dashed] node[above,sloped] {$\color{black} e_3$} (c)
    (a) edge[purple, line width=1.3pt, dashed] node {$\color{black} e_4$} (e);

\end{tikzpicture}
    \caption{$P=(e_2,e_1,e_4)$.}
    \label{fig:harder_path_P}
\end{subfigure}
\caption{Constructing the alternating subtrail $P=(\{v,w\}_{.,t}, \{w,x\}_{.,t}, \{x,y\}_{.,.})$ with $v \neq y$. Legend as in Figure \ref{fig:multiple_cycles_case2_block_vanuit_C1}.}
\end{figure}

Next, we consider the case in which $M(H_0 \Delta H^*)$ contains exactly one minimal alternating cycle $C^*$, i.e., $M(H_0 \Delta H^*)=C^*$.
\\
\begin{lemma}
\label{lemma:connected_when_one_cycle}
Let $\vb*{d} \in D$, let $H_0,H^* \in \mathcal{H}^{\textnormal{stub}}_{s}(\vb*{d})$ and $\mathcal{E}(M(H_0 \Delta H^*)) = k+1$. If $M(H_0 \Delta H^*)$ consists of exactly one minimal alternating cycle, then $\mathcal{G}(\mathcal{H}^{\textnormal{stub}}_{s}(\vb*{d}))$ contains a path from $H_0$ to $H^*$.
\end{lemma}
This lemma cannot be proved with the same arguments as Lemma \ref{lemma:connected_when_multiple_cycles}, since inequality (\ref{eq:length cycle leq k}) does not hold. Before proving this lemma, we show that the alternating cycle $C^*$ contains an alternating trail of length 3, of which two adjacent edges have the same tail. Then, we use this trail to construct $\Tilde{H}$.
\\
\begin{claim}
\label{claim:subpath_P}
    If $M(H_0 \Delta H^*)$ contains exactly one minimal alternating cycle $C^*$, then $C^*$ contains an alternating subtrail $P=(\{a,b\}_{.,t}, \{b,c\}_{.,t}, \{c,d\}_{.,.})$ with $a \neq d$.
\end{claim}
\begin{proof}
    If $M(H_0 \Delta H^*)$ contains exactly one minimal alternating cycle $C^*$, then $|C^*|\geq 4$. Indeed, if instead $C^*=(e_1,e_2)$ then $e_1=\{a,b\}_{H_0,t}, e_2=\{a,b\}_{H^*,s}$, with $t \neq s$ (by Property \ref{prop:second_att}), so that $|f_{H_0,t}(M(H_0 \Delta H^*))| > |f_{H^*,t}(M(H_0 \Delta H^*))|$, which contradicts Claim \ref{claim:H_in_H_s_equivalences}.
    
    Suppose that for all alternating subtrails $(e_1,e_2)$ of cycle $C^*$ holds $e_1^2 \neq e_2^2$. Then, all edges in $f_{H_0}(M(H_0 \Delta H^*))$ would have the same tail $t_1$ and all edges in $f_{H^*}(M(H_0 \Delta H^*))$ would have the same tail $t_2$ with $t_1 \neq t_2$ (Figure \ref{fig:single_cycle_no_neighbors_with_same_tail}). Then, $|f_{H_0,\tau}(M(H_0 \Delta H^*))| \neq |f_{H^*,\tau}(M(H_0 \Delta H^*))|$, which contradicts Claim \ref{claim:H_in_H_s_equivalences}. Therefore, we conclude that there exists a subtrail $(e_1,e_2)$ of cycle $C^*$ with $e_1^2 = e_2^2$. Let $e_1=\{a,b\}_{H_0,t}, e_2=\{b,c\}_{H^*,t}$. By Property \ref{prop:second_att}, $(e_1,e_2)$ cannot be an alternating 2-cycle. Therefore, $a \neq c$. Since $M(H_0 \Delta H^*)$ consists of alternating cycles (Lemma \ref{lemma:decomposition_min_alt_cycles}), another edge $e_3 = \{c,d\}_{H_0,.}$ must exist. If $a \neq d$ then the lemma holds for the trail $P=(e_1,e_2,e_3)$ (Figure \ref{fig:easy_path_P}, where $e_1^2=e_2^2=\tau$ and $e_3^2=\sigma$). 
    
    If $a=d$, then another edge $e_4= \{a,e\}_{H^*,.}$ must exist since $M(H_0 \Delta H^*)$ consists of alternating cycles (Lemma \ref{lemma:decomposition_min_alt_cycles}). We obtain that $e \neq b$, otherwise $(e_1,e_4)$ would be an alternating 2-cycle, which cannot be contained in the minimal alternating cycle $C^*$. Similarly, $e \neq c$. The lemma now holds for the trail $P=(e_2,e_1,e_3)$ (Figure \ref{fig:harder_path_P}, where $e_1^2=e_2^2=\tau$ and $e_3^2=\sigma$, $e_4^2=\sigma$).
\end{proof}

Using Claim \ref{claim:subpath_P}, we can prove Lemma \ref{lemma:connected_when_one_cycle}.

\begin{proof}[Proof of Lemma \ref{lemma:connected_when_one_cycle}]
Given a subtrail $P=(e_1,e_2,e_3)$ of cycle $C^*$, where $e_1=\{a,b\}_{.,t},e_2=\{b,c\}_{.,t}, e_3=\{c,d\}_{.,s}$ for some $s,t \in \{\tau,\sigma\}$, as described in Claim \ref{claim:subpath_P}, we create an alternating cycle that can be used to construct $\Tilde{H}$. Let $e_4=\{a,d\}_{.,s}$ and consider the alternating cycle $(e_1,e_2,e_3,e_4)$. We distinguish four cases based on whether the hyperarc $M^{-1}(e_4)$ is present in $H_0$ or $H^*$, which are illustrated in Figure \ref{fig:one_cycle_4cases}:
\\
\begin{enumerate}
    \item $M^{-1}(e_4) \in A(H_0) \cap A(H^*)$;
    \item $M^{-1}(e_4) \notin A(H_0) \cup A(H^*)$;
    \item $M^{-1}(e_4) \in A(H^*) \backslash A(H_0)$, so $e_4 \in f_{H^*,s}(C^*)$;
    \item $M^{-1}(e_4) \in A(H_0) \backslash A(H^*)$, so $e_4 \in f_{H_0,s}(C^*)$. \label{case:case_that_splits_in_4}
\end{enumerate}
Case \ref{case:case_that_splits_in_4} is split into four subcases, which are covered later.

\begin{figure}[tbp]
\centering
\begin{subfigure}{0.4\textwidth}
    \centering
    \begin{tikzpicture}[-,>=Stealth,auto,
                    thick,main node/.style={circle,draw, minimum size=\tikznodesize}, inner sep=1pt]

  \node[main node] (a) at (0,0) {$a$};
  \node[main node] (b) at (0,2*\tikzscale) {$b$};
  \node[main node] (c) at (2*\tikzscale,2*\tikzscale) {$c$};
  \node[main node] (d) at (2*\tikzscale,0) {$d$};

  \path[every node/.style={font=\sffamily\small}]
    (a) edge[cyan] node {$\color{black} e_1$} (b)
    (b) edge[purple, line width=1.3pt] node {$\color{black} e_2$} (c)
    (c) edge[cyan, dashed] node {$\color{black} e_3$} (d)
    (d) edge[lightgray, dashed] node {$\color{black} e_4$} (a);

\end{tikzpicture}
    \caption{Case 1: $M^{-1}(e_4) \in A(H_0) \cap A(H^*)$.}
    \label{fig:one_cycle_case1}
\end{subfigure}
\hspace{1 cm}
\begin{subfigure}{0.4\textwidth}
    \centering
    \begin{tikzpicture}[-,>=Stealth,auto,
                    thick,main node/.style={circle,draw, minimum size=\tikznodesize}, inner sep=1pt]

  \node[main node] (a) at (0,0) {$a$};
  \node[main node] (b) at (0,2*\tikzscale) {$b$};
  \node[main node] (c) at (2*\tikzscale,2*\tikzscale) {$c$};
  \node[main node] (d) at (2*\tikzscale,0) {$d$};

  \path[every node/.style={font=\sffamily\small}]
    (a) edge[cyan] node {$\color{black} e_1$} (b)
    (b) edge[purple, line width=1.3pt] node {$\color{black} e_2$} (c)
    (c) edge[cyan, dashed] node {$\color{black} e_3$} (d);

\end{tikzpicture}
    \caption{Case 2: $M^{-1}(e_4) \notin A(H_0) \cup A(H^*)$.}
    \label{fig:one_cycle_case2}
\end{subfigure}\\ 
\vspace{1cm}
\begin{subfigure}{0.4\textwidth}
    \centering
    \begin{tikzpicture}[-,>=Stealth,auto,
                    thick,main node/.style={circle,draw, minimum size=\tikznodesize}, inner sep=1pt]

  \node[main node] (a) at (0,0) {$a$};
  \node[main node] (b) at (0,2*\tikzscale) {$b$};
  \node[main node] (c) at (2*\tikzscale,2*\tikzscale) {$c$};
  \node[main node] (d) at (2*\tikzscale,0) {$d$};

  \path[every node/.style={font=\sffamily\small}]
    (a) edge[cyan] node {$\color{black} e_1$} (b)
    (b) edge[purple, line width=1.3pt] node {$\color{black} e_2$} (c)
    (c) edge[cyan, dashed] node {$\color{black} e_3$} (d)
    (d) edge[purple, line width=1.3pt, dashed] node {$\color{black} e_4$} (a);

\end{tikzpicture}
    \caption{Case 3: $e_4 \in f_{H^*,s}$.}
    \label{fig:one_cycle_case3}
\end{subfigure}
\hspace{1cm}
\begin{subfigure}{0.4\textwidth}
    \centering
    \begin{tikzpicture}[-,>=Stealth,auto,
                    thick,main node/.style={circle,draw, minimum size=\tikznodesize}, inner sep=1pt]

  \node[main node] (a) at (0,0) {$a$};
  \node[main node] (b) at (0,2*\tikzscale) {$b$};
  \node[main node] (c) at (2*\tikzscale,2*\tikzscale) {$c$};
  \node[main node] (d) at (2*\tikzscale,0) {$d$};

  \path[every node/.style={font=\sffamily\small}]
    (a) edge[cyan] node {$\color{black} e_1$} (b)
    (b) edge[purple, line width=1.3pt] node {$\color{black} e_2$} (c)
    (c) edge[cyan, dashed] node {$\color{black} e_3$} (d)
    (d) edge[cyan, dashed] node {$\color{black} e_4$} (a);

\end{tikzpicture}
    \caption{Case 4: $e_4 \in f_{H_0,s}$.}
    \label{fig:one_cycle_case4}
\end{subfigure}
\caption{Cases 1-4 for the edge $e_4$. Legend as in Figure \ref{fig:multiple_cycles_case2_block_vanuit_C1}.}
\label{fig:one_cycle_4cases}
\end{figure}
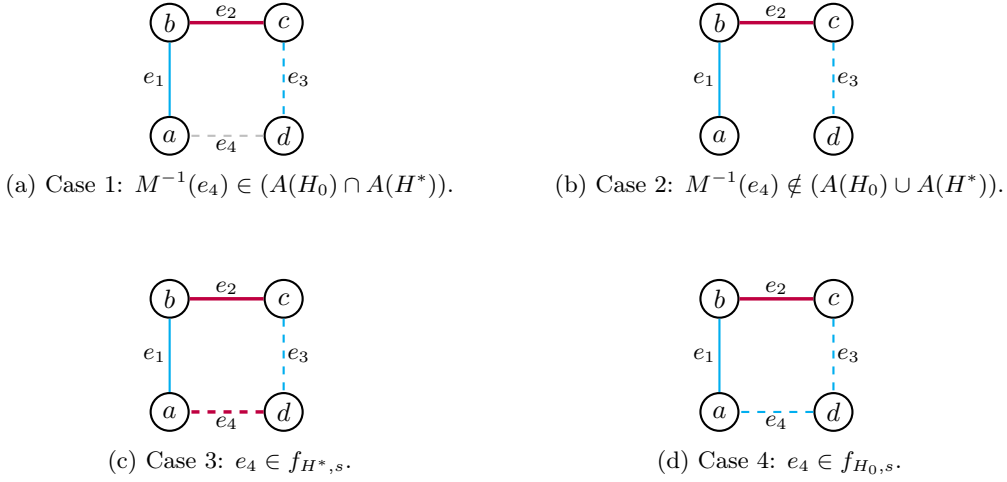

\underline{Case 1}: $M^{-1}(e_4)=(s,\{a,d\}) \in A(H_0) \cap A(H^*)$ (Figure \ref{fig:one_cycle_case1})\\
We construct $\Tilde{H}$ from $H^*$ using the alternating cycle $(e_1,e_2,e_3,e_4)$:
\begin{align*}
    A(\Tilde{H}) = (A(H^*) \backslash \{M^{-1}(e_2), M^{-1}(e_4)\} )\cup \{M^{-1}(e_1),M^{-1}(e_3)\}.
\end{align*}
By Claim \ref{claim:H_in_H_s_equivalences}, $\Tilde{H} \in \mathcal{H}^{\textnormal{stub}}_{s}(\vb*{d})$, as $\Tilde{H}$ contains no degenerate or multi-hyperarcs by construction, $M(\Tilde{H} \Delta H^*) = \{e_1,e_2,e_3,e_4\}$ is an alternating cycle and
\begin{align*}
    |f_{\Tilde{H},\tau}(M(\Tilde{H} \Delta H^*))| &= \mathds{1}_{\{t=\tau\}} + \mathds{1}_{\{s=\tau\}}= |f_{H^*,\tau}(M(\Tilde{H} \Delta H^*))|,
\end{align*}
by construction, where $\mathds{1}$ denotes the indicator function. Moreover,
\begin{align*}
    \mathcal{E}(M(H_0 \Delta \Tilde{H})) &= \mathcal{E}(M(H_0 \Delta H^*)) -3+1 \\
    &= k+1-2 = k-1
\end{align*}
and
\begin{align*}
    \mathcal{E}(M(\Tilde{H} \Delta H^*))= 4 \leq k,
\end{align*}
since $k+1=\mathcal{E}(M(H_0 \Delta H^*)) \geq 5$ so $k \geq 4$. By the induction hypothesis, $\mathcal{G}(\mathcal{H}^{\textnormal{stub}}_{s}(\vb*{d}))$ contains a path from $H_0$ to $\Tilde{H}$ and from $\Tilde{H}$ to $H^*$.

\underline{Case 2}: $M^{-1}(e_4)=(s,\{a,d\}) \notin A(H_0) \cup A(H^*)$ (Figure \ref{fig:one_cycle_case2})\\
We construct $\Tilde{H}$ from $H_0$ using the cycle $(e_1,e_2,e_3,e_4)$:
\begin{align*}
    A(\Tilde{H}) = (A(H_0) \backslash \{M^{-1}(e_1), M^{-1}(e_3)\} )\cup \{M^{-1}(e_2),M^{-1}(e_4)\}.
\end{align*}
With similar arguments as in case 1 of this proof, we can show that $\Tilde{H} \in \mathcal{H}^{\textnormal{stub}}_{s}(\vb*{d})$, $\mathcal{E}(M(H_0 \Delta \Tilde{H})) \leq k$ and $\mathcal{E}(M(\Tilde{H} \Delta H^*)) = k-1$. By the induction hypothesis, $\mathcal{G}(\mathcal{H}^{\textnormal{stub}}_{s}(\vb*{d}))$ contains a path from $H_0$ to $\Tilde{H}$ and from $\Tilde{H}$ to $H^*$.

\underline{Case 3}: $e_4=\{a,d\}_{H^*,s} \in f_{H^*,s}(C^*)$ (Figure \ref{fig:one_cycle_case3})\\
Since $C^*$ is a minimal alternating cycle, we obtain that $C^*=(e_1,e_2,e_3,e_4)$. Thus, $|H_0 \Delta H^*|=|C^*|=4$ and 
\begin{align*}
    |f_{\Tilde{H},\tau}(M(H_0 \Delta H^*))| &= \mathds{1}_{\{t=\tau\}} + \mathds{1}_{\{s=\tau\}}= |f_{H^*,\tau}(M(H_0 \Delta H^*))|.
\end{align*}
By Lemma \ref{lemma:sym_diff_hypershuffle}, $H_0$ and $H^*$ are adjacent in $\mathcal{G}(\mathcal{H}^{\textnormal{stub}}_{s}(\vb*{d}))$.

\underline{Case 4}: $e_4=\{a,d\}_{H_0,s} \in f_{H_0,s}(C^*)$ (Figure \ref{fig:one_cycle_case4})\\
First, we show that there must be an edge $\{a,f\}_{H^*,\cdot}$ in $M(H_0 \Delta H^*)$, for some $f \notin \{a,b,c,d\}$. 

Since $M(H_0 \Delta H^*)$ consists of alternating cycles (Lemma \ref{lemma:decomposition_min_alt_cycles}), there must be at least two edges in $f_{H^*}(M(H_0 \Delta H^*))$ incident to vertex $a$. Let $e_5$ be such an edge. We now show that $e_5$ is incident to a vertex $f \notin \{a,b,c,d\}$.
   
Assume by contradiction that $e_5$ is incident to a vertex in $\{b,c,d\}$. It cannot be incident to $b$, nor $d$, since that creates an alternating 2-cycle $(e_1,e_5)$ or $(e_4,e_5)$, contradicting that $C^*$ is a minimal alternating cycle. Therefore, $e_5$ must be incident to $c$: $e_5=\{a,c\}_{H^*,.}$ (Figure \ref{fig:edge_(a,c)_cannot_exist}).
   
As $M(H_0 \Delta H^*)$ consists of alternating cycles (Lemma \ref{lemma:decomposition_min_alt_cycles}), there must also be an edge in $f_{H^*}(M(H_0 \Delta H^*))$ incident to vertex $d$. Let this edge be $e_i=\{d,.\}_{H^*,.}$. Following the alternating cycle starting from $e_i$, we obtain an alternating trail $Q=(e_i,\hdots,e_j)$ of which none of the edges are in $\{e_1,e_2,e_3,e_4,e_5\}$. Moreover, the trail $Q$ must end at some vertex in $\{a,c,d\}$, to fulfill that there as many edges from $f_{H_0}(M(H_0 \Delta H^*))$ as edges from $f_{H^*}(M(H_0 \Delta H^*))$ incident to any vertex (Claim \ref{claim:H_in_H_s_equivalences}) If the trail ends at $a$, then we obtain $e_j=\{.,a\}_{H^*,.}$, and there exists an alternating subcycle $(Q,e_4)$ (Figure \ref{fig:P_ends_at_a}). This contradicts that $C^*$ is a minimal alternating cycle. If the trail ends at $c$, then we obtain $e_j=\{.,c\}_{H_0,.}$, and there exists an alternating subcycle $(Q,e_5,e_4)$ (Figure \ref{fig:P_ends_at_c}), and if the trail ends at $d$, then we obtain $e_j=\{.,d\}_{H^*}$, and there exists an alternating subcycle $(Q,e_3,e_5,e_4)$ (Figure \ref{fig:P_ends_at_d}). This contradicts that $C^*$ is a minimal alternating cycle. We conclude that the trail $Q$ cannot exist, so $e_5$ must be incident to some vertex $f \notin \{a,b,c,d\}$. Let $e_5=\{a,f\}_{H^*,u}$ for some $u \in \{\tau,\sigma\}$. This is illustrated in Figure \ref{fig:existence_e_5}, for $u=\tau$.

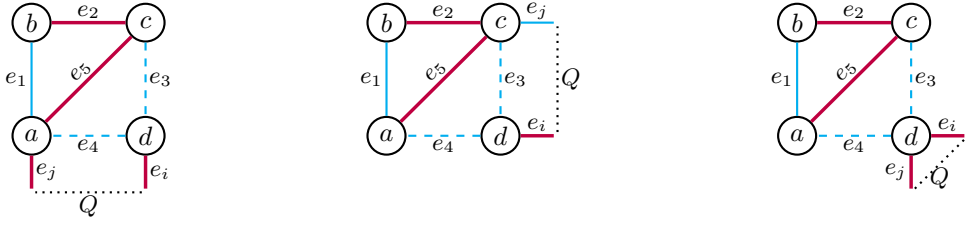
\begin{figure}[tbp]
\centering
\begin{subfigure}[t]{0.3\textwidth}
    \centering
    \begin{tikzpicture}[-,>=Stealth,auto,
                    thick,main node/.style={circle,draw, minimum size=\tikznodesize}, inner sep=1pt]

  \node[main node] (a) at (0,0) {$a$};
  \node[main node] (b) at (0,2*\tikzscale) {$b$};
  \node[main node] (c) at (2*\tikzscale,2*\tikzscale) {$c$};
  \node[main node] (d) at (2*\tikzscale,0) {$d$};

  %invisible nodes used to draw extra paths
  \node (under a) at (0,-1*\tikzscale) {};
  \node (under d) at (2*\tikzscale,-1*\tikzscale) {};

  \path[every node/.style={font=\sffamily\small}]
    (a) edge[cyan] node {$\color{black} e_1$} (b)
    (b) edge[purple, line width=1.3pt] node {$\color{black} e_2$} (c)
    (c) edge[cyan, dashed] node {$\color{black} e_3$} (d)
    (d) edge[cyan, dashed] node {$\color{black} e_4$} (a)
    (a) edge[purple, line width=1.3pt] node[above,sloped] {$\color{black} e_5$} (c)
    (d) edge[purple, line width=1.3pt] node {$\color{black} e_i$} (under d)
    (under d) edge[black, dotted] node {$\color{black} Q$} (under a)
    (a) edge[purple, line width=1.3pt] node {$\color{black} e_j$} (under a);

    % invisible line to make the figure as big as the other figures
    \draw[line width=1pt, color=black, opacity=0.0] (0,0) -- (0,-1.5*\tikzscale);
\end{tikzpicture}
    \caption{Alternating trail $Q$ ends at vertex $a$. Subcycle: $(Q,e_4)$.}
    \label{fig:P_ends_at_a}
\end{subfigure}
\hspace{0.3 cm}
\begin{subfigure}[t]{0.3\textwidth}
    \centering
    \begin{tikzpicture}[-,>=Stealth,auto,
                    thick,main node/.style={circle,draw, minimum size=\tikznodesize}, inner sep=1pt]

  \node[main node] (a) at (0,0) {$a$};
  \node[main node] (b) at (0,2*\tikzscale) {$b$};
  \node[main node] (c) at (2*\tikzscale,2*\tikzscale) {$c$};
  \node[main node] (d) at (2*\tikzscale,0) {$d$};

  %invisible nodes used to draw extra paths
  \node (right of c) at (3*\tikzscale,2*\tikzscale) {};
  \node (right of d) at (3*\tikzscale,0) {};

  \path[every node/.style={font=\sffamily\small}]
    (a) edge[cyan] node {$\color{black} e_1$} (b)
    (b) edge[purple, line width=1.3pt] node {$\color{black} e_2$} (c)
    (c) edge[cyan, dashed] node {$\color{black} e_3$} (d)
    (d) edge[cyan, dashed] node {$\color{black} e_4$} (a)
    (a) edge[purple, line width=1.3pt] node[above,sloped] {$\color{black} e_5$} (c)
    (d) edge[purple, line width=1.3pt] node {$\color{black} e_i$} (right of d)
    (right of d) edge[black, dotted] node[right] {$\color{black} Q$} (right of c)
    (c) edge[cyan] node {$\color{black} e_j$} (right of c);

    % invisible line to make the figure as big as the other figures
    \draw[line width=1pt, color=black, opacity=0.0] (0,0) -- (0,-1.5*\tikzscale);
\end{tikzpicture}
    \caption{Alternating trail $Q$ ends at vertex $c$. Subcycle: $(Q,e_5,e_4)$.}
    \label{fig:P_ends_at_c}
\end{subfigure}
\hspace{0.3 cm}
\begin{subfigure}[t]{0.33\textwidth}
    \centering
    \begin{tikzpicture}[-,>=Stealth,auto,
                    thick,main node/.style={circle,draw, minimum size=\tikznodesize}, inner sep=1pt]

  \node[main node] (a) at (0,0) {$a$};
  \node[main node] (b) at (0,2*\tikzscale) {$b$};
  \node[main node] (c) at (2*\tikzscale,2*\tikzscale) {$c$};
  \node[main node] (d) at (2*\tikzscale,0) {$d$};

  %invisible nodes used to draw extra paths
  \node (right of d) at (3*\tikzscale,0) {};
  \node (below d) at (2*\tikzscale,-1*\tikzscale) {};

  \path[every node/.style={font=\sffamily\small}]
    (a) edge[cyan] node {$\color{black} e_1$} (b)
    (b) edge[purple, line width=1.3pt] node {$\color{black} e_2$} (c)
    (c) edge[cyan, dashed] node {$\color{black} e_3$} (d)
    (d) edge[cyan, dashed] node {$\color{black} e_4$} (a)
    (a) edge[purple, line width=1.3pt] node[above,sloped] {$\color{black} e_5$} (c)
    (d) edge[purple, line width=1.3pt] node {$\color{black} e_i$} (right of d)
    (below d) edge[black, dotted] node[below] {$\color{black} Q$} (right of d)
    (d) edge[purple, line width=1.3pt] node[left] {$\color{black} e_j$} (below d);

    % invisible line to make the figure as big as the other figures
    \draw[line width=1pt, color=black, opacity=0.0] (0,0) -- (0,-1.5*\tikzscale);
\end{tikzpicture}
    \caption{Alternating trail $Q$ ends at vertex $d$. Subcycle: $(Q,e_3,e_5,e_4)$.}
    \label{fig:P_ends_at_d}
\end{subfigure}
\caption{Cases that prove that $e_5=\{a,c\}_{H^*,.}$ cannot exist. Legend as in Figure \ref{fig:multiple_cycles_case2_block_vanuit_C1}.}
\label{fig:edge_(a,c)_cannot_exist}
\end{figure}

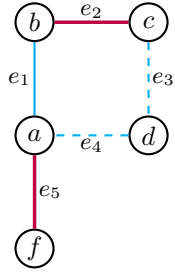
\begin{figure}[tbp]
    \centering
    \begin{tikzpicture}[-,>=Stealth,auto,
                    thick,main node/.style={circle,draw, minimum size=\tikznodesize}, inner sep=1pt]

  \node[main node] (a) at (0,0) {$a$};
  \node[main node] (b) at (0,2*\tikzscale) {$b$};
  \node[main node] (c) at (2*\tikzscale,2*\tikzscale) {$c$};
  \node[main node] (d) at (2*\tikzscale,0) {$d$};
  \node[main node] (f) at (-2*\tikzscale, 0) {$f$};

  \path[every node/.style={font=\sffamily\small}]
    (a) edge[cyan] node {$\color{black} e_1$} (b)
    (b) edge[purple, line width=1.3pt] node {$\color{black} e_2$} (c)
    (c) edge[cyan, dashed] node {$\color{black} e_3$} (d)
    (d) edge[cyan, dashed] node {$\color{black} e_4$} (a)
    (a) edge[purple, line width=1.3pt] node {$\color{black} e_5$} (f);

\end{tikzpicture}
    \caption{The edge $e_5$, connecting $a$ to a vertex $f \notin \{a,b,c,d\}$, must exist. Legend as in Figure \ref{fig:multiple_cycles_case2_block_vanuit_C1}.}
    \label{fig:existence_e_5}
\end{figure}

Let $e_6=\{c,f\}_{.,u}$. We construct $\Tilde{H}$ by using the cycle $(e_5,e_1,e_2,e_6)$. We distinguish four cases based on whether $M^{-1}(e_6)$ is a hyperarc in $H_0$ or $H^*$, which are illustrated in Figure \ref{cases4.1-4.4}:

\begin{enumerate}[label=4.\arabic*.]
    \item $M^{-1}(e_6) \in A(H_0) \cap A(H^*)$;
    \item $M^{-1}(e_6) \notin A(H_0) \cup A(H^*)$;
    \item $M^{-1}(e_6) \in A(H_0) \backslash A(H^*)$, so $e_6=\{f,c\}_{H_0,u} \in f_{H_0,u}(C^*)$;
    \item $M^{-1}(e_6) \in A(H^*) \backslash A(H_0)$, so $e_6=\{f,c\}_{H^*,u} \in f_{H^*,u}(C^*)$.
\end{enumerate}

\begin{figure}[tb]
\centering
\begin{subfigure}{0.4\textwidth}
    \begin{tikzpicture}[-,>=Stealth,auto,
                    thick,main node/.style={circle,draw, minimum size=\tikznodesize}, inner sep=1pt]

  \node[main node] (a) at (0,0) {$a$};
  \node[main node] (b) at (0,2*\tikzscale) {$b$};
  \node[main node] (c) at (2*\tikzscale,2*\tikzscale) {$c$};
  \node[main node] (d) at (2*\tikzscale,0) {$d$};
  \node[main node] (f) at (-2*\tikzscale,0) {$f$};

  \path[every node/.style={font=\sffamily\small}]
    (a) edge[cyan] node[right, pos=0.25] {$\color{black} e_1$} (b)
    (b) edge[purple, line width=1.3pt] node {$\color{black} e_2$} (c)
    (c) edge[cyan, dashed] node {$\color{black} e_3$} (d)
    (d) edge[cyan, dashed] node[below] {$\color{black} e_4$} (a)
    (a) edge[purple, line width=1.3pt] node[below] {$\color{black} e_5$} (f)
    (f) edge[lightgray, dashed] node[above,pos=0.25, sloped] {$\color{black} e_6$} (c);

\end{tikzpicture}
    \centering
    \caption{Case 4.1: $M^{-1}(e_6) \in A(H_0) \cap A(H^*)$.}
    \label{fig:one_cycle_case4.1}
\end{subfigure}
\hspace{1cm}
\begin{subfigure}{0.4\textwidth}
    \begin{tikzpicture}[-,>=Stealth,auto,
                    thick,main node/.style={circle,draw, minimum size=\tikznodesize}, inner sep=1pt]

  \node[main node] (a) at (0,0) {$a$};
  \node[main node] (b) at (0,2*\tikzscale) {$b$};
  \node[main node] (c) at (2*\tikzscale,2*\tikzscale) {$c$};
  \node[main node] (d) at (2*\tikzscale,0) {$d$};
  \node[main node] (f) at (-2*\tikzscale,0) {$f$};

  \path[every node/.style={font=\sffamily\small}]
    (a) edge[cyan] node {$\color{black} e_1$} (b)
    (b) edge[purple, line width=1.3pt] node {$\color{black} e_2$} (c)
    (c) edge[cyan, dashed] node {$\color{black} e_3$} (d)
    (d) edge[cyan, dashed] node {$\color{black} e_4$} (a)
    (a) edge[purple, line width=1.3pt] node[below] {$\color{black} e_5$} (f);

\end{tikzpicture}
    \centering
    \caption{Case 4.2: $M^{-1}(e_6) \notin A(H_0) \cup A(H^*)$.}
    \label{fig:one_cycle_case4.2}
\end{subfigure}\\
\vspace{1cm}
\begin{subfigure}{0.4\textwidth}
    \begin{tikzpicture}[-,>=Stealth,auto,
                    thick,main node/.style={circle,draw, minimum size=\tikznodesize}, inner sep=1pt]

  \node[main node] (a) at (0,0) {$a$};
  \node[main node] (b) at (0,2*\tikzscale) {$b$};
  \node[main node] (c) at (2*\tikzscale,2*\tikzscale) {$c$};
  \node[main node] (d) at (2*\tikzscale,0) {$d$};
  \node[main node] (f) at (-2*\tikzscale,0) {$f$};

  \path[every node/.style={font=\sffamily\small}]
    (a) edge[cyan] node[right, pos=0.25] {$\color{black} e_1$} (b)
    (b) edge[purple, line width=1.3pt] node {$\color{black} e_2$} (c)
    (c) edge[cyan, dashed] node {$\color{black} e_3$} (d)
    (d) edge[cyan, dashed] node[below] {$\color{black} e_4$} (a)
    (a) edge[purple, line width=1.3pt] node[below] {$\color{black} e_5$} (f)
    (f) edge[cyan, dashed] node[above,pos=0.25, sloped] {$\color{black} e_6$} (c);

\end{tikzpicture}
    \centering
    \caption{Case 4.3: $e_6 \in f_{H_0,u}$.}
    \label{fig:one_cycle_case4.3}
\end{subfigure}
\hspace{1cm}
\begin{subfigure}{0.4\textwidth}
    \begin{tikzpicture}[-,>=Stealth,auto,
                    thick,main node/.style={circle,draw, minimum size=\tikznodesize}, inner sep=1pt]

  \node[main node] (a) at (0,0) {$a$};
  \node[main node] (b) at (0,2*\tikzscale) {$b$};
  \node[main node] (c) at (2*\tikzscale,2*\tikzscale) {$c$};
  \node[main node] (d) at (2*\tikzscale,0) {$d$};
  \node[main node] (f) at (-2*\tikzscale,0) {$f$};

  \path[every node/.style={font=\sffamily\small}]
    (a) edge[cyan] node[right, pos=0.25] {$\color{black} e_1$} (b)
    (b) edge[purple, line width=1.3pt] node {$\color{black} e_2$} (c)
    (c) edge[cyan, dashed] node {$\color{black} e_3$} (d)
    (d) edge[cyan, dashed] node[below] {$\color{black} e_4$} (a)
    (a) edge[purple, line width=1.3pt] node[below] {$\color{black} e_5$} (f)
    (f) edge[purple, line width=1.3pt, dashed] node[above, sloped, pos=0.25] {$\color{black} e_6$} (c);

\end{tikzpicture}
    \centering
    \caption{Case 4.4: $e_6 \in f_{H^*,u}$.}
    \label{fig:one_cycle_case4.4}
\end{subfigure}
\caption{Cases 4.1-4.4 for the edge $e_6$. Legend as in Figure \ref{fig:multiple_cycles_case2_block_vanuit_C1}.}
\label{cases4.1-4.4}
\end{figure}
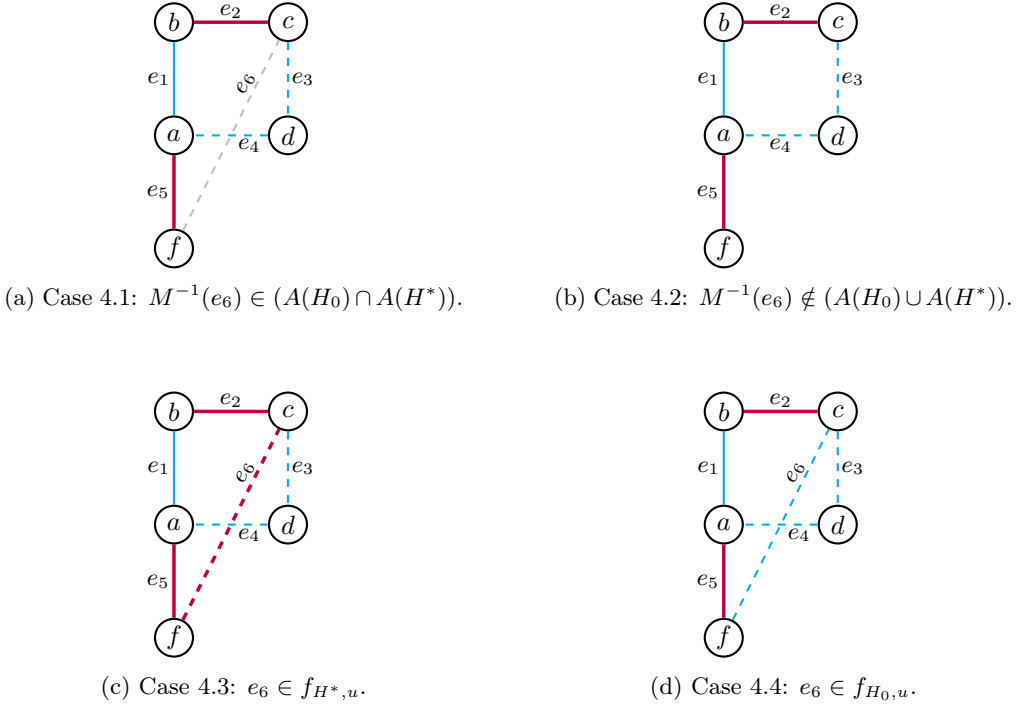

Case 4.1: $M^{-1}(e_6)=(u,\{f,c\}) \in A(H_0) \cap A(H^*)$ (Figure \ref{fig:one_cycle_case4.1})\\
We construct $\Tilde{H}$ from $H_0$ using the alternating cycle $(e_1,e_2,e_6,e_5)$:
\begin{align*}
    A(\Tilde{H}) = (A(H_0) \backslash \{M^{-1}(e_1), M^{-1}(e_6)\} )\cup \{M^{-1}(e_2),M^{-1}(e_5)\}.
\end{align*}
With similar arguments as in case 1 of this proof, we can show that $\Tilde{H} \in \mathcal{H}^{\textnormal{stub}}_{s}(\vb*{d})$, $\mathcal{E}(M(H_0 \Delta \Tilde{H})) \leq k$ and $\mathcal{E}(M(\Tilde{H} \Delta H^*)) = k-1$. By the induction hypothesis, $\mathcal{G}(\mathcal{H}^{\textnormal{stub}}_{s}(\vb*{d}))$ contains a path from $H_0$ to $\Tilde{H}$ and from $\Tilde{H}$ to $H^*$.

Case 4.2: $M^{-1}(e_6)=(u,\{f,c\}) \notin A(H_0) \cup A(H^*)$ (Figure \ref{fig:one_cycle_case4.2})\\
We construct $\Tilde{H}$ from $H^*$ using the alternating cycle $(e_1,e_2,e_6,e_5)$:
\begin{align*}
    A(\Tilde{H}) = (A(H^*) \backslash \{M^{-1}(e_2), M^{-1}(e_5)\} )\cup \{M^{-1}(e_1),M^{-1}(e_6)\}.
\end{align*}
With similar arguments as in case 1 of this proof, we can show that $\Tilde{H} \in \mathcal{H}^{\textnormal{stub}}_{s}(\vb*{d})$, $\mathcal{E}(M(H_0 \Delta \Tilde{H})) =k-1$ and $\mathcal{E}(M(\Tilde{H} \Delta H^*)) \leq k$. By the induction hypothesis, $\mathcal{G}(\mathcal{H}^{\textnormal{stub}}_{s}(\vb*{d}))$ contains a path from $H_0$ to $\Tilde{H}$ and from $\Tilde{H}$ to $H^*$.

\begin{figure}[tbp]
\centering
\begin{subfigure}{0.4\textwidth}
    \centering
    \begin{tikzpicture}[-,>=Stealth,auto,
                    thick,main node/.style={circle,draw, minimum size=\tikznodesize}, inner sep=1pt]

  \node[main node] (a) at (0,0) {$a$};
  \node[main node] (b) at (0,2*\tikzscale) {$b$};
  \node[main node] (c) at (2*\tikzscale,2*\tikzscale) {$c$};
  \node[main node] (d) at (2*\tikzscale,0) {$d$};
  \node[main node] (f) at (-2*\tikzscale,0) {$f$};

   %invisible nodes used to draw extra paths
  \node (right of d) at (3*\tikzscale,0) {};
  \node (below d) at (2*\tikzscale,-\tikzscale) {};

  \path[every node/.style={font=\sffamily\small}]
    (a) edge[cyan] node[right, pos=0.25] {$\color{black} e_1$} (b)
    (b) edge[purple, line width=1.3pt] node {$\color{black} e_2$} (c)
    (c) edge[cyan, dashed] node[left] {$\color{black} e_3$} (d)
    (d) edge[cyan, dashed] node[above] {$\color{black} e_4$} (a)
    (a) edge[purple, line width=1.3pt] node[below] {$\color{black} e_5$} (f)
    (f) edge[purple, line width=1.3pt, dashed] node[above,pos=0.25, sloped] {$\color{black} e_6$} (c)
    (c) edge[cyan] node[above,sloped] {$\color{black} e_i$} (right of d)
    (right of d) edge[black, dotted] node[below, sloped] {$\color{black} Q$} (below d)
    (a) edge[purple, line width=1.3pt] node[below,sloped] {$\color{black} e_j$} (below d);

\end{tikzpicture}
    \caption{Alternating trail $Q$ ends at vertex $a$. Subcycle: $(Q,e_1,e_2)$.}
    \label{fig:P2_ends_at_a}
\end{subfigure}
\hspace{1cm}
\begin{subfigure}{0.4\textwidth}
    \centering
    \begin{tikzpicture}[-,>=Stealth,auto,
                    thick,main node/.style={circle,draw, minimum size=\tikznodesize}, inner sep=1pt]

  \node[main node] (a) at (0,0) {$a$};
  \node[main node] (b) at (0,2*\tikzscale) {$b$};
  \node[main node] (c) at (2*\tikzscale,2*\tikzscale) {$c$};
  \node[main node] (d) at (2*\tikzscale,0) {$d$};
  \node[main node] (f) at (-2*\tikzscale,0) {$f$};

   %invisible nodes used to draw extra paths
  \node (above b) at (0,3*\tikzscale) {};
  \node (above f) at (-2*\tikzscale,2*\tikzscale) {};

  \path[every node/.style={font=\sffamily\small}]
    (a) edge[cyan] node[right, pos=0.25] {$\color{black} e_1$} (b)
    (b) edge[purple, line width=1.3pt] node {$\color{black} e_2$} (c)
    (c) edge[cyan, dashed] node[right] {$\color{black} e_3$} (d)
    (d) edge[cyan, dashed] node[below] {$\color{black} e_4$} (a)
    (a) edge[purple, line width=1.3pt] node[below] {$\color{black} e_5$} (f)
    (f) edge[purple, line width=1.3pt, dashed] node[above,pos=0.25, sloped] {$\color{black} e_6$} (c)
    (c) edge[cyan] node[above, sloped] {$\color{black} e_i$} (above b)
    (above b) edge[black, dotted] node[above, sloped] {$\color{black} Q$} (above f)
    (f) edge[cyan] node[left] {$\color{black} e_j$} (above f);

    % invisible line to make the figure as big as the other figures
    \draw[line width=1pt, color=black, opacity=0.0] (0,-\tikzscale) -- (\tikzscale, \tikzscale);

\end{tikzpicture}
    \caption{Alternating trail $Q$ ends at vertex $f$. Subcycle: $(Q,e_6)$.}
    \label{fig:P2_ends_at_f}
\end{subfigure}\\
\vspace{1cm}
\begin{subfigure}[t]{0.4\textwidth}
    \centering
    \begin{tikzpicture}[-,>=Stealth,auto,
                    thick,main node/.style={circle,draw, minimum size=\tikznodesize}, inner sep=1pt]

  \node[main node] (a) at (0,0) {$a$};
  \node[main node] (b) at (0,2*\tikzscale) {$b$};
  \node[main node] (c) at (2*\tikzscale,2*\tikzscale) {$c$};
  \node[main node] (d) at (2*\tikzscale,0) {$d$};
  \node[main node] (f) at (-2*\tikzscale,0) {$f$};

   %invisible nodes used to draw extra paths
  \node (right of c) at (3*\tikzscale,2*\tikzscale) {};
  \node (right of d) at (3*\tikzscale,0) {};
  \node (below d) at (2*\tikzscale,-\tikzscale) {};
  \node (below a) at (0,-\tikzscale) {};

  \path[every node/.style={font=\sffamily\small}]
    (a) edge[cyan] node[right, pos=0.25] {$\color{black} e_1$} (b)
    (b) edge[purple, line width=1.3pt] node {$\color{black} e_2$} (c)
    (c) edge[cyan, dashed] node {$\color{black} e_3$} (d)
    (d) edge[cyan, dashed] node[below] {$\color{black} e_4$} (a)
    
    (a) edge[purple, line width=1.3pt] node[below] {$\color{black} e_5$} (f)
    (f) edge[purple, line width=1.3pt, dashed] node[above,pos=0.25, sloped] {$\color{black} e_6$} (c)
    
    (c) edge[cyan] node {$\color{black} e_i$} (right of c)
    (right of c) edge[black, dotted] node[right] {$\color{black} Q$} (right of d)
    (d) edge[purple, line width=1.3pt] node {$\color{black} e_j$} (right of d)

    (d) edge[purple, line width=1.3pt] node {$\color{black} e_{i'}$} (below d)
    (below d) edge[black, dotted] node[below] {$\color{black} Q'$} (below a)
    (below a) edge[purple, line width=1.3pt] node {$\color{black} e_{j'}$} (a);
\end{tikzpicture}
    \caption{Alternating trail $Q$ ends at vertex $d$, alternating trail $Q'$ ends at vertex $a$. Subcycle: $(Q',e_4)$.}
    \label{fig:P'_ends_at_a}
\end{subfigure}
\hspace{1cm}
\begin{subfigure}[t]{0.4\textwidth}
    \centering
    \begin{tikzpicture}[-,>=Stealth,auto,
                    thick,main node/.style={circle,draw, minimum size=\tikznodesize}, inner sep=1pt]

  \node[main node] (a) at (0,0) {$a$};
  \node[main node] (b) at (0,2*\tikzscale) {$b$};
  \node[main node] (c) at (2*\tikzscale,2*\tikzscale) {$c$};
  \node[main node] (d) at (2*\tikzscale,0) {$d$};
  \node[main node] (f) at (-2*\tikzscale,0) {$f$};

   %invisible nodes used to draw extra paths
  \node (right of c) at (3*\tikzscale,2*\tikzscale) {};
  \node (right of d) at (3*\tikzscale,0) {};
  \node (below d) at (2*\tikzscale,-\tikzscale) {};
  \node (below f) at (-2*\tikzscale,-\tikzscale) {};

  \path[every node/.style={font=\sffamily\small}]
    (a) edge[cyan] node[right, pos=0.25] {$\color{black} e_1$} (b)
    (b) edge[purple, line width=1.3pt] node {$\color{black} e_2$} (c)
    (c) edge[cyan, dashed] node {$\color{black} e_3$} (d)
    (d) edge[cyan, dashed] node[below] {$\color{black} e_4$} (a)
    
    (a) edge[purple, line width=1.3pt] node[below] {$\color{black} e_5$} (f)
    (f) edge[purple, line width=1.3pt, dashed] node[above,pos=0.25, sloped] {$\color{black} e_6$} (c)
    
    (c) edge[cyan] node {$\color{black} e_i$} (right of c)
    (right of c) edge[black, dotted] node[right] {$\color{black} Q$} (right of d)
    (d) edge[purple, line width=1.3pt] node {$\color{black} e_j$} (right of d)

    (d) edge[purple, line width=1.3pt] node {$\color{black} e_{i'}$} (below d)
    (below d) edge[black, dotted] node[below] {$\color{black} Q'$} (below f)
    (below f) edge[cyan] node {$\color{black} e_{j'}$} (f);

\end{tikzpicture}
    \caption{Alternating trail $Q$ ends at vertex $d$, alternating trail $Q'$ ends at vertex $f$. Subcycle: $(Q',e_6,e_3)$.}
    \label{fig:P'_ends_at_f}
\end{subfigure}
\caption{Cases to show that $e_6=\{f,c\}_{H^*,u} \in f_{H^*,u}(M(H_0,H^*))$ cannot exist. Legend as in Figure \ref{fig:multiple_cycles_case2_block_vanuit_C1}.}
\end{figure}
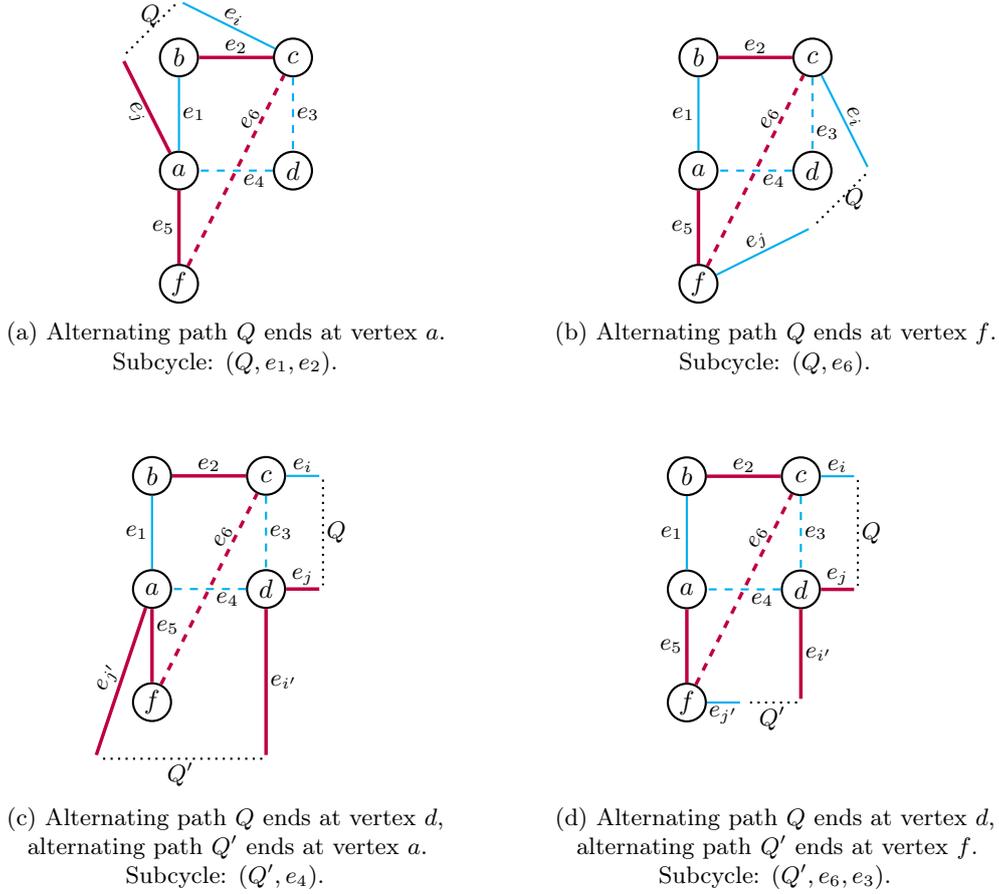

Case 4.3: $e_6=\{f,c\}_{H_0,u} \in f_{H_0,u}(M(H_0 \Delta H^*))$ (Figure \ref{fig:one_cycle_case4.3})\\
This case cannot occur, because there exists an alternating subcycle $\{e_1,e_2,e_6,e_5\}$, which contradicts that $C^*$ is a minimal alternating cycle.

Case 4.4: $e_6=\{f,c\}_{H^*,u} \in f_{H^*,u}(M(H_0 \Delta H^*))$ (Figure \ref{fig:one_cycle_case4.4})\\
We show that this case cannot occur, as it would imply the existence of an alternating subcycle of $C^*$. Since $M(H_0 \Delta H^*)$ consists of alternating cycles (Lemma \ref{lemma:decomposition_min_alt_cycles}), there must be at least one edge in $f_{H_0}(M(H_0 \Delta H^*)) \backslash \{e_3\}$ incident to vertex $c$. Let $e_i=\{c,.\}_{H_0,.}$ be such an edge. Following the alternating cycle starting from $e_i$, we obtain an alternating trail $Q=(e_i,\hdots,e_j)$ of which none of the edges are in $\{e_1,e_2,e_3,e_4,e_5,e_6\}$. Moreover, the trail $Q$ must end at some vertex in $\{a,d,f\}$, to fulfill that there as many edges from $f_{H_0}(M(H_0 \Delta H^*))$ as edges from $f_{H^*}(M(H_0 \Delta H^*))$ incident to any vertex (Lemma \ref{lemma:decomposition_alt_cycles}). If the trail ends at $a$, then we obtain $e_j=\{.,a\}_{H^*,.}$, and there exists an alternating subcycle $(Q,e_1,e_2)$ (Figure \ref{fig:P2_ends_at_a}). This contradicts that $C^*$ is a minimal alternating cycle. If the trail ends at $f$, then we obtain $e_j=\{.,f\}_{H_0,.}$, and there exists an alternating subcycle $(Q,e_6)$ (Figure \ref{fig:P2_ends_at_f}). This contradicts that $C^*$ is a minimal alternating cycle. If the trail ends at $d$, then we obtain $e_j=\{.,d\}_{H^*}$. Since $M(H_0 \Delta H^*)$ consists of alternating cycles (Lemma \ref{lemma:decomposition_min_alt_cycles}), there must be another edge in $f_{H^*}(M(H_0 \Delta H^*))$ incident to $d$. Let $e_{i'}=\{d,.\}_{H^*,.}$ be such an edge. Following the alternating cycle starting from $e_{i'}$, we obtain an alternating trail $Q'=(e_{i'},\hdots,e_{j'})$ of which none of the edges are in $\{e_1,e_2,e_3,e_4,e_6\} \cup Q$. Moreover, the trail $Q'$ must end at some vertex in $\{a,f\}$. If the trail ends at $a$, the we obtain $e_{j'}=\{.,a\}_{H^*,.}$, and there exists an alternating subcycle $(Q',e_4)$ (Figure \ref{fig:P'_ends_at_a}). This contradicts that $C^*$ is a minimal alternating cycle. If the trail ends at $f$, then we obtain $e_{j'}=\{.,f\}_{H_0,.}$, and there exists an alternating subcycle $(Q',e_6,e_3)$ (Figure \ref{fig:P'_ends_at_f}). This contradicts that $C^*$ is a minimal alternating cycle. 

We conclude that such trails $Q$ and $Q'$ cannot exist, so this case cannot occur.
\end{proof}

We are now ready to prove Lemma \ref{lemma:stub_s_connected}.

\begin{proof}[Proof of Lemma \ref{lemma:stub_s_connected}]
    The result follows from Lemmas \ref{lemma:connected_when_multiple_cycles} and \ref{lemma:connected_when_one_cycle}.
\end{proof}

Finally, we prove Theorem \ref{thm:uniform_degrees(2,1)}.
\begin{proof}[Proof of Theorem \ref{thm:uniform_degrees(2,1)}]
    Let $\vb*{d} \in D$. Since $\mathcal{G}(\mathcal{H}^{\textnormal{stub}}_{s}(\vb*{d}))$ is aperiodic (Lemma \ref{lemma:aperiodic}) and strongly connected (Lemma \ref{lemma:stub_s_connected}) for $\vb*{d}$ as in the theorem, random walks on $\mathcal{G}(\mathcal{H}^{\textnormal{stub}}_{s}(\vb*{d}))$ are ergodic. Since $\mathcal{G}(\mathcal{H}^{\textnormal{stub}}_{s}(\vb*{d}))$ is also regular (Lemma \ref{lemma:regular}), $\mathcal{G}(\mathcal{H}^{\textnormal{stub}}_{s}(\vb*{d}))$ has a uniform stationary distribution.
\end{proof}

\section{Vertex labelled spaces: Proof of Theorem \ref{thm:stub_to_vert}}
\label{section:proofThm4}
Lastly, we prove Theorem \ref{thm:stub_to_vert}. To that end, we first show that strong connectivity and aperiodicity of $\mathcal{G}(\mathcal{H}^{\textnormal{vert}}_{x}(\vb*{d}))$ follow directly from strong connectivity and aperiodicity of $\mathcal{G}(\mathcal{H}^{\textnormal{stub}}_{x}(\vb*{d}))$. Second, we show that $\mathcal{G}(\mathcal{H}^{\textnormal{vert}}_{x}(\vb*{d}))$ is a regular graph.

Let us introduce a mapping that relates $\mathcal{G}(\mathcal{H}^{\textnormal{vert}}_{x}(\vb*{d}))$ to $\mathcal{G}(\mathcal{H}^{\textnormal{stub}}_{x}(\vb*{d}))$.
\\
\begin{definition}[Map $g$ (\cite{chodrow2019})]
    Let $g: \mathcal{H}^{\textnormal{stub}}_{x}(\vb*{d}) \rightarrow \mathcal{H}^{\textnormal{vert}}_{x}(\vb*{d})$ be the function that maps a stub-labeled directed hypergraph to its vertex-labeled equivalent. Note that $g$ is not one-to-one. Let $g^{-1}(H)$ denote the set of stub-labeled directed hypergraphs that map to $H$.
\end{definition}

The graph $\mathcal{G}(\mathcal{H}^{\textnormal{vert}}_{x}(\vb*{d}))$ can be considered the same graph as $\mathcal{G}(\mathcal{H}^{\textnormal{stub}}_{x}(\vb*{d}))$, but with vertices $H,H'$ merged if $H$ and $H'$ are equal up to stub-labels, i.e., $H,H' \in g^{-1}(H^*)$ for some $H^* \in \mathcal{H}^{\textnormal{vert}}_{x}(\vb*{d})$. 
\\
\begin{claim}
\label{claim:adjacency_vert_stub}
    Two vertices $v,w$ in $\mathcal{G}(\mathcal{H}^{\textnormal{vert}}_{x}(\vb*{d}))$ are adjacent if and only if there exist two vertices $v' \in g^{-1}(v)$ and $w' \in g^{-1}(w)$ that are adjacent in $\mathcal{G}(\mathcal{H}^{\textnormal{stub}}_{x}(\vb*{d}))$.
\end{claim}
\begin{proof}
    Two vertices $v,w$ in $\mathcal{G}(\mathcal{H}^{\textnormal{vert}}_{x}(\vb*{d}))$ are adjacent if and only if there exists a double hyperarc shuffle that transforms $v$ into $w$. This is true if and only if there exists a double hyperarc shuffle that transforms $v'$ into $w'$, for some vertices $v' \in g^{-1}(v)$ and $w' \in g^{-1}(w)$. This is the case if and only if $v'$ and $w'$ are adjacent in $\mathcal{G}(\mathcal{H}^{\textnormal{stub}}_{x}(\vb*{d}))$.
\end{proof}

Now, we can show when $\mathcal{G}(\mathcal{H}^{\textnormal{vert}}_{x}(\vb*{d}))$ is strongly connected.
\\
\begin{lemma}
\label{lemma:vert_connected}
    Let $x \subseteq \{s,d,m\}$. Then, $\mathcal{G}(\mathcal{H}^{\textnormal{vert}}_{x}(\vb*{d}))$ is a strongly connected graph for some degree sequence $\vb*{d}$, if and only if $\mathcal{G}(\mathcal{H}^{\textnormal{stub}}_{x}(\vb*{d}))$ is a strongly connected graph for degree sequence $\vb*{d}$.
\end{lemma}
\begin{proof}
    The result follows directly from Claim \ref{claim:adjacency_vert_stub}.
\end{proof}

Next, we show that $\mathcal{G}(\mathcal{H}^{\textnormal{vert}}_{x}(\vb*{d}))$ is aperiodic.
\\
\begin{lemma}
\label{lemma:vert_aperiodic}
     $\mathcal{G}(\mathcal{H}^{\textnormal{vert}}_x(\vb*{d}))$ is an aperiodic graph, for any degree sequence $\vb*{d}$ and any $x \subseteq \{s,d,m\}$.
\end{lemma}
\begin{proof}
    In the proof of Lemma \ref{lemma:aperiodic}, we show that $\mathcal{G}(\mathcal{H}^{\textnormal{stub}}_x(\vb*{d}))$ contains self-loops. By Claim \ref{claim:adjacency_vert_stub}, $\mathcal{G}(\mathcal{H}^{\textnormal{vert}}_x(\vb*{d}))$ also contains self-loops and is therefore aperiodic.
\end{proof}

Now, we show that $\mathcal{G}(\mathcal{H}^{\textnormal{vert}}_x(\vb*{d}))$ is a regular graph.
\\
\begin{lemma}
\label{lemma:vert_regular}
    $\mathcal{G}(\mathcal{H}^{\textnormal{vert}}_x(\vb*{d}))$, with transition probabilities $\alpha$ as given in Theorem \ref{thm:stub_to_vert}, is a regular graph, for any degree sequence $\vb*{d}$ and any $x \subseteq \{s,d,m\}$.
\end{lemma}
\begin{proof}
    We show that the transition probability between any two vertices $H_0,H^* $ in $\mathcal{G}(\mathcal{H}^{\textnormal{vert}}_{x}(\vb*{d}))$ equals the transition probability between two vertices $S_0,S^*$ in $\mathcal{G}(\mathcal{H}^{\textnormal{stub}}_{x}(\vb*{d}))$, with $S_0 \in g^{-1}(H_0)$ and $S^* \in g^{-1}(H^*)$.
 
    Let $H_0,H^* \in \mathcal{H}^{\textnormal{vert}}_x(\vb*{d})$ be two vertices in $\mathcal{G}(\mathcal{H}^{\textnormal{vert}}_x(\vb*{d}))$ with $H^* \in S(a,b|H_0)$ and let $p^{\textnormal{vert}}_{a,b}(H^*|H_0)$ denote the transition probability with acceptance rate $\alpha$, when the double hyperarc shuffle $s(a,b|H_0)$ is applied. Similarly, let $p_{a,b}(S^*|S_0)$ denote the transition probability between two stub-labeled directed hypergraphs $S_0$ and $S^*$ when the double hyperarc shuffle $s(a,b|S_0)$ is applied. Let $S_0 \in g^{-1}(H_0)$ and $S^* \in g^{-1}(H^*)$ s.t. $S^* \in S(a,b|S_0)$. Then,
    \begin{align*}
        p^{\textnormal{vert}}_{a,b}(H^*|H_0) &= \sum_{S' \in g^{-1}(H^*)} \alpha(S'|S_0) p_{a,b}(S'|S_0)\\
        &= \alpha(S^*|S_0) \sum_{S' \in g^{-1}(H^*)} p_{a,b}(S'|S_0)\\
        &= \alpha(S^*|S_0) \frac{1}{m_am_b\alpha(S^*|S_0)} p_{a,b}(S^*|S_0)\\
        &= \frac{1}{m_am_b} p_{a,b}(S^*|S_0).
    \end{align*}
    The first equality is by construction. The second equality holds since the value of $\alpha(S'|S_0)$ is equal for all $S' \in g^{-1}(H^*)$ with $p_{a,b}(S'|S_0)>0$. For the third equality, firstly note that $p_{a,b}(S'|S_0)$ is equal for all $S' \in g^{-1}(H^*)$ with $p_{a,b}(S'|S_0)>0$. Secondly, note that the sum contains exactly $\frac{1}{m_a m_b \alpha(S^*|S_0)}$  non-zero terms: 
\begin{itemize}
        \item for every vertex $v$, there are exactly $\binom{m_{\hat{a}^{\textnormal{h}}}(v) + m_{\hat{b}^{\textnormal{h}}}(v)}{m_{\hat{a}^{\textnormal{h}}}(v)}$ ways to rearrange the stubs within the heads of the resulting hyperarc heads;
        \item similarly, for every vertex $v$, there are exactly $ \binom{m_{\hat{a}^{\textnormal{t}}}(v) + m_{\hat{b}^{\textnormal{t}}}(v)}{m_{\hat{a}^{\textnormal{t}}}(v)}$ ways to arrange the stubs within the tails.
\end{itemize}
All of these configurations result in a distinct element in $g^{-1}(H^*)$. Therefore, the sum contains $ \prod_{v \in V} \binom{m_{\hat{a}^{\textnormal{h}}}(v) + m_{\hat{b}^{\textnormal{h}}}(v)}{m_{\hat{a}^{\textnormal{h}}}(v)} \binom{m_{\hat{a}^{\textnormal{t}}}(v) + m_{\hat{b}^{\textnormal{t}}}(v)}{m_{\hat{a}^{\textnormal{t}}}(v)}= \frac{1}{m_a m_b \alpha(S^*|S_0)}$  non-zero terms. Finally, let $p^{\textnormal{vert}}(H^*|H_0)$ denote the transition probability with acceptance rate $\alpha$, for any double hyperarc shuffle. Then,
\begin{align*}
    p^{\textnormal{vert}}(H^*|H_0) &= \binom{|A(H_0)|}{2}^{-1}  \sum_{\hat{a},\hat{b} \in A(H_0)} p^{\textnormal{vert}}_{\hat{a},\hat{b}}(H^*|H_0) \\
    &= \binom{|A(H_0)|}{2}^{-1}  m_a m_b  p^{\textnormal{vert}}_{a,b}(H^*|H_0)\\
    &=  \binom{|A(H_0)|}{2}^{-1}  p_{a,b}(S^*|S_0)\\
    &= p(S^*|S_0).
\end{align*}
The first equality is by construction. The second equality follows from $H^* \in S(a,b|H_0)$. The fourth equality follows from uniqueness of every stub in a stub-labeled directed hypergraph.

The lemma follows from $\mathcal{G}(\mathcal{H}^{\textnormal{stub}}_{x}(\vb*{d}))$ being a regular graph (Lemma \ref{lemma:regular}).
\end{proof}

Finally, we prove Theorem \ref{thm:stub_to_vert}.
\begin{proof}[Proof of Theorem \ref{thm:stub_to_vert}]
For any $x \subseteq \{s,d,m\}$, $\mathcal{G}(\mathcal{H}^{\textnormal{stub}}_{x}(\vb*{d}))$ is regular and aperiodic (Lemmas \ref{lemma:regular} and \ref{lemma:aperiodic}, resp.), so the double hyperarc shuffling method on the space $\mathcal{H}^{\textnormal{stub}}_{x}(\vb*{d})$ samples uniformly from that space if and only if $\mathcal{G}(\mathcal{H}^{\textnormal{stub}}_{x}(\vb*{d}))$ is strongly connected. 

Now, $\mathcal{G}(\mathcal{H}^{\textnormal{vert}}_{x}(\vb*{d}))$ is strongly connected if and only if $\mathcal{G}(\mathcal{H}^{\textnormal{stub}}_{x}(\vb*{d}))$ is strongly connected (Lemma \ref{lemma:vert_connected}), so if and only if the double hyperarc shuffling method on the space $\mathcal{H}^{\textnormal{stub}}_{x}(\vb*{d})$ samples uniformly from $\mathcal{H}^{\textnormal{stub}}_{x}(\vb*{d})$. Therefore, $\mathcal{G}(\mathcal{H}^{\textnormal{vert}}_{x}(\vb*{d}))$ is ergodic if and only if the double hyperarc shuffling method on $\mathcal{H}^{\textnormal{stub}}_{x}(\vb*{d})$ samples uniformly from $\mathcal{H}^{\textnormal{stub}}_{x}(\vb*{d})$ (Lemmas \ref{lemma:vert_aperiodic} and \ref{lemma:vert_connected}). In addition, $\mathcal{G}(\mathcal{H}^{\textnormal{vert}}_{x}(\vb*{d}))$ is regular (Lemma \ref{lemma:vert_regular}), so has a uniform stationary distribution if and only if the double hyperarc shuffling method on $\mathcal{H}^{\textnormal{stub}}_{x}(\vb*{d})$ samples uniformly from $\mathcal{H}^{\textnormal{stub}}_{x}(\vb*{d})$.
\end{proof}

\paragraph{Acknowledgements.} Y.K. and C.S. are funded through NWO M2 grant 0.379.

\bibliographystyle{apalike}
\bibliography{mybib.bib}

\end{document}